\documentclass[12pt]{article}
\usepackage{mathrsfs}
\usepackage{amsmath}
\usepackage{amssymb}
\usepackage{amsfonts}
\usepackage{mathtools}
\usepackage[top=2cm,bottom=2cm,left=2.5cm,right=2cm]{geometry}
\usepackage{graphicx}
\usepackage[colorinlistoftodos]{todonotes}
\setuptodonotes{
  inline,
  backgroundcolor=blue!8!white,
  bordercolor=blue!65!black,
  linecolor=blue!65!black
}
\usepackage{amsthm}
\usepackage{bbm}
\usepackage{enumerate}
\usepackage{enumitem}
\usepackage{verbatim}
\usepackage[
  backend=biber,
  style=alphabetic,
  sorting=nyt,
  maxnames=99,
  giveninits=true,
  doi=true,
  url=false,
  isbn=false,
  eprint=true
]{biblatex}
\usepackage{float}
\usepackage{tikz-cd}
\usepackage[all]{xy}
\usepackage[T1]{fontenc}
\usepackage{xcolor}
\usepackage[skins,breakable]{tcolorbox}
\usepackage{adjustbox}
\usepackage[bottom]{footmisc}

\newtcolorbox{commentbox}[1][]{colback=yellow!10!white, colframe=red!50!black, title=#1}

\usepackage[affil-it]{authblk}
\usepackage{commands}
\usepackage[colorlinks=true, allcolors=blue]{hyperref}

\newif\ifshowcomments%
\showcommentstrue%

\title{The Chabauty--Kim dimension inequality via noncommutative Iwasawa theory}
\author{Sa'ar Zehavi\thanks{saarzehavi@gmail.com}}
\affil{Max Planck Institute for Mathematics}

\begin{document}
\maketitle

\begin{abstract}
We establish Kim's conjectural dimension inequality for smooth
projective geometrically connected curves $C/\Q$ of genus at
least two with a rational basepoint. For every odd prime $p$
of good reduction, we prove that
\[
 \dim H_f^1(G_{\Q,\Sigma},U_n^\et)
 <
 \dim(U_n^{\dR}/F^0)
\]
for all sufficiently large $n$, where $U_n^\et$ and $U_n^\dR$ are the $n$'th
lower central series truncations of the $\Q_p$-pro-unipotent geometric \'etale,
and de Rham fundamental groups of $C$, and $\Sigma$ is any finite set of places
containing $p$, $\infty$, and the primes of bad reduction.
Inspired by Coates and Kim's work on curves with CM Jacobians,
we use noncommutative Iwasawa theory to bound second Galois
cohomology through dimensions of kernels of certain Iwasawa operators,
acting on finite dimensional representations of the ambient generalized symplectic group.
We establish bounds for these kernels in terms of the highest weight
and dimension of the ambient symplectic representation.
No restriction is imposed on the Jacobian's geometric
endomorphism ring or its Tate-module Galois image.
Consequently, the Chabauty--Kim loci are finite at all
sufficiently large levels, recovering the Mordell theorem
over $\Q$ through the Chabauty--Kim method.
\end{abstract}

\tableofcontents
\section{Introduction}
Let $S$ denote a finite set of rational primes, and let $p \not\in S$ denote an odd\footnote{To see where the oddness of $p$ plays a role, see Remark~\ref{rem:very-good-primes}.} auxiliary prime. Fix $\Sigma := S\cup\{p,\infty\}$, and write $G_{\mathbb{Q},\Sigma} := \Gal(\mathbb{Q}_\Sigma/\mathbb{Q})$ for the Galois group of the maximal extension of $\mathbb{Q}$ unramified outside $\Sigma$. Let $\mathcal{A}/\mathbb{Q}$ denote a principally polarized $g$-dimensional abelian variety with good reduction outside $S$. Write $V := T_p\mathcal{A}\otimes_{\Zp}\Qp$ for its (rationalized) $p$-adic Tate module equipped with its Galois-action, and $\omega: V\otimes V\to \Qp(1)$ for its symplectic form, coming from the Weil pairing and the polarization. The Galois action preserves the $\Qp\cdot\omega$-line, hence acts on $(V,\omega)$ through symplectomorphisms. Therefore, the Galois representation $\rho_V: G_{\mathbb{Q},\Sigma}\to \GL(V)$ factors through $\GSp(V,\omega) = \GSp_{2g,\Qp}(\Qp)$. 

The irreducible representations\footnote{Throughout the manuscript, the term irreducible representations should always refer to finite dimensional irreducible representations.} of $\GSp_{2g,\Qp}$ are determined by their restrictions to $\Sp_{2g,\Qp}$ and to $\GSp_{2g,\Qp}$'s center. Recall that the irreducible representations of $\Sp_{2g,\Qp}\subseteq \GSp_{2g,\Qp}$ are in bijection with the highest weight vectors $\lambda\in \mathbb{Z}^g$; these are descending $g$-tuples of non-negative integers,
\[
\lambda_1\ge \dots\ge \lambda_g\ge 0.
\]
For such a highest weight vector, denote by $V_\lambda$ the corresponding finite dimensional rational $p$-adic representation of $\Sp_{2g,\Qp}$. A representation $W$ of $\GSp_{2g,\Qp}$ is said to have weight $s$, when the center, whose $\Qp$-points consist of the scalar matrices $t\cdot I\in \GSp_{2g,\Qp}(\Qp)$, with $I = I_{2g}$ the identity matrix, acts on $W$ via multiplication by $t^s$. We shall denote by $V_{\lambda,s}$ the irreducible representation of $\GSp_{2g,\Qp}$ of weight $s$, whose restriction to $\Sp_{2g,\Qp}$ is isomorphic to $V_\lambda$. Endow $V_{\lambda,s}$ with the Galois structure obtained by the factorization $G_{\mathbb{Q},\Sigma}\to \GSp_{2g,\Qp}(\Qp)\to \GL(V_{\lambda,s})$.

Write 
\[
d_\lambda := \dim_\Qp V_{\lambda,s},\qquad ||\lambda||^2 := \sum_{i=1}^g \lambda_i^2,\qquad D:=\dim \Sp_{2g,\Qp} = 2g^2 + g.
\]
Our main theorem reads as follows.
\begin{theorem}[Main theorem]\label{thm:h2-bound}
There is a finite set of weights, $E\subseteq \Z$, such that for all irreducible representations $V_{\lambda,s}$ of $\GSp_{2g,\Qp}$ with $s\not\in E$, one has
\[
\dim_\Qp H^2(G_{\mathbb{Q},\Sigma}, V_{\lambda,s}) \ll \min\left\{d_\lambda,\dfrac{(1 + ||\lambda||)^{D - g - 1}}{d_\lambda}\right\},
\]
where the bound $\dim_\Qp H^2(G_{\mathbb{Q},\Sigma}, V_{\lambda,s}) \ll d_\lambda$ holds for all $s$.
\end{theorem}
The proof of Theorem~\ref{thm:h2-bound} occurs in \S\ref{ssec:h2-bound}; see \S\ref{sssec:h2-bound-proof-idea} for an overview of the proof.

Our main application is in Chabauty--Kim theory. Namely, we solve a conjecture of Kim~\cite[Conjecture~1, p.~96]{KIM}; see Conjecture~\ref{conj:kim-dimension} below, proving that the Chabauty--Kim method works unconditionally for all smooth projective curves of genus $g\ge 2$ over $\mathbb{Q}$. In particular, Kim's method successfully recovers the Mordell conjecture. Before we state our results more concretely, let us review the basic setup in Chabauty--Kim theory.
\subsection{Motivation and main results}
Let $C/\mathbb{Q}$ denote a smooth quasi-projective curve of genus $g \ge 2$ with good reduction outside $S$. Let $b\in C(\mathbb{Q})$ denote a rational point, and let $U^\et$ and $U^\dR$ denote the pro-unipotent \'etale and de Rham fundamental groups based at $b$, in the sense of Kim~\cite[\S1, p.~97, and \S2, p.~114]{KIM}. In short, these are affine the pro-algebraic $\Qp$-group schemes, which are Tannaka dual to the (neutral) Tannakian categories of iterated extensions of: (i) the trivial $\Qp$-local system on the geometric generic fibre of $C$ in the \'etale realization, and (ii) the trivial bundle with (integrable) connection on $C\otimes \Qp$ in the de Rham realization.

\begin{definition}[The descending central series filtration]
For a (pro-)unipotent group scheme $U$, write
\[
U = U^1 \supseteq U^2\supseteq U^3\supseteq \dots,\qquad U^{n+1} := \overline{[U,U^n]},
\]
for its (closed) descending central series (DCS) filtration, we denote the induced DCS truncations via $U_n := U/U^{n+1}$.
\end{definition}
The fundamental groups $U^\et$ and $U^\dR$, as well as their DCS truncations, have extra structures: a $G_{\mathbb{Q},\Sigma}$-action in the \'etale realization and an ascending Hodge filtration and a Frobenius automorphism in the de Rham realization. Write $F^0 := F^0U^\dR$ for the $\mathrm{Fil}^0$-piece of $U^\dR$'s Hodge filtration, and $U^\dR/F^0$ for its associated quotient, which is representable by an affine scheme as well. Fix a (unipotent) truncation parameter $n\ge 1$. In~\cite[p.~95 and \S2, p.~120]{KIM}, Kim introduced the following diagram:
\begin{equation}\label{eq:kim_diagram}
\begin{tikzcd}
\mathcal{C}(\mathbb{Z}[1/S])\arrow[r]\arrow[d,"j_n^\et"]&\mathcal{C}(\mathbb{Z}_p)\arrow[rd,"\Phi_p"]\arrow[d,"j_{n,p}^\et"]\\
H^1_f(G_{\mathbb{Q},\Sigma}, U_n^\et)(\mathbb{Q}_p)\arrow[r,"\mathrm{loc}_p"]&H^1_f(G_{\mathbb{Q}_p}, U_n^\et)(\mathbb{Q}_p)\arrow[r,"\sim"]& U_n^\dR/F^0(\mathbb{Q}_p).
\end{tikzcd}
\end{equation}
Here, $\mathcal{C}$ is a smooth $\mathbb{Z}[1/S]$-model of $C$. The top row consists of the $\mathbb{Z}[1/S]$ and $\Zp$-points of $\mathcal{C}$. The bottom row consists of the $\mathbb{Q}_p$-points of certain moduli spaces associated to $C$, which are representable by affine schemes. From left to right, these are Kim's global and local Selmer schemes, $H^1_f(G_{\mathbb{Q},\Sigma}, U_n^\et)$ and $H^1_f(G_{\mathbb{Q}_p}, U_n^\et)$, respectively; and the de Rham moduli space $U_n^\dR/F^0$. We defer their definitions to Kim~\cite[Proposition~1, pp.~105--106, and \S2, pp.~118--120]{KIM}. All maps in the bottom row are the maps induced on $\mathbb{Q}_p$-points from algebraic morphisms of affine schemes. The first map is induced by restriction in nonabelian group cohomology, while the second map, called the local Bloch--Kato logarithm, is an isomorphism. Their composition, i.e.~the map $\log_{\mathrm{BK}}: H^1_f(G_{\mathbb{Q},\Sigma}, U_n^\et)\to U_n^\dR/F^0$ is called the (global) Bloch--Kato logarithm. The vertical arrows are the (global and local) Kummer maps; these are set theoretic maps, which associate to $\mathbb{Z}[1/S]$ or $\Zp$-points of $\mathcal{C}$, $U_n^\et$-torsors over $\Spec\,\Qp$, equipped with an equivariant action of $G_{\mathbb{Q},\Sigma}$ or $G_{\mathbb{Q}_p}$, respectively.

The diagonal arrow, $\Phi_p: \mathcal{C}(\Zp)\to U_n^\dR/F^0(\mathbb{Q}_p)$ is Kim's period map (see~\cite[\S1, p.~106]{KIM}). Although, similar to the Kummer maps, it may be viewed set theoretically: it associates $\Zp$-points of $\mathcal{C}$ with admissible $U_n^\dR$-torsors, this map has an analytic structure too. Write $k = \Fp$ and $\mathcal{C}(k)$ for the finite set of points in the special fibre of $\mathcal{C}$. For $\underline{x}\in \mathcal{C}(k)$, denote by $]\underline{x}[$ the tube of $\underline{x}$, also called the residue disk above $\underline{x}$. This is a rigid analytic subspace of the rigid analytification of $C$; its $\Qp$-points are precisely those $\Zp$-points of $\mathcal{C}$, which reduce to $\underline{x}$ in the special fibre. Importantly, the period map restricts to (rigid) analytic maps on each residue disk, and these restrictions are all Zariski dense; see~\cite[Theorem~1, p.~93, and its proof, pp.~111--114]{KIM}.

Kim observed that whenever the dimension of the global Selmer scheme is strictly smaller than the dimension of the de Rham moduli space, i.e.~whenever
\begin{equation}\label{eq:kim_dim_ineq}
\dim H^1_f(G_{\mathbb{Q},\Sigma}, U_n^\et) < \dim U_n^\dR/F^0,
\end{equation}
one recovers an elegant proof of the finiteness of $\#\mathcal{C}(\mathbb{Z}[1/S])$. We repeat Kim's argument~\cite[\S3, pp.~120--121]{KIM} here for the benefit of the reader. Under our assumed dimension inequality, the (Zariski) closure of the image of the (global) Bloch--Kato logarithm, i.e.~$\overline{\mathrm{Im}\left(\log_{\mathrm{BK}}\right)}$, is contained in a proper closed subscheme of $U_n^\dR/F^0$. Since the latter is an affine scheme, there exists some nontrivial global algebraic function $0\neq f\in \Gamma(U_n^\dR/F^0; \mathcal{O}_{U_n^\dR/F^0})$, which vanishes on $\overline{\mathrm{Im}\left(\log_{\mathrm{BK}}\right)}$. The pullback $\Phi_p^*f$ is a Coleman function in the sense of~\cite[Definitions~4.1 and~4.3]{BESSER}, and restricts to a rigid analytic function on each residue disk $]\underline{x}[$ for $x\in \mathcal{C}(k)$. Since the restrictions $\Phi_p|_{]\underline{x}[}$ are Zariski dense, the rigid analytic functions $\Phi_p^*f|_{]\underline{x}[}$ are all nontrivial. Finally, since Kim's diagram~\eqref{eq:kim_diagram} is commutative,
\[
\mathcal{C}(\mathbb{Z}[1/S])\cap ]\underline{x}[(\Qp)\subseteq \left\{\Phi_p^*f|_{]\underline{x}[} = 0\right\},
\]
for all $\underline{x}\in \mathcal{C}(k)$. That is, the set of $S$-integral points within each of the finitely many residue disks is cut out by the zero-set of a nontrivial analytic function on that same disk. Since the residue disks are compact, and a nontrivial analytic function can only have finitely many zeros within a compact domain, the sets $\{\Phi_p^*f|_{]\underline{x}[} = 0\}$ are finite, hence, so are the $\mathcal{C}(\mathbb{Z}[1/S])\cap ]\underline{x}[(\Qp)$. Applying this argument repeatedly within each of the finitely many residue disks, one establishes the finiteness of $\#\mathcal{C}(\mathbb{Z}[1/S])$.

A key advantage of Kim's method is its computational aspect: when~\eqref{eq:kim_dim_ineq} holds, it not only proves finiteness but also suggests how to determine the integral points, by constructing a nontrivial algebraic function on $U^\dR_n/F^0$ vanishing on the global Selmer image.\footnote{The resulting common zero locus may contain additional $p$-adic points. Kim's stronger conjecture predicts that for a large enough unipotent truncation parameter $n\ge 1$, the Chabauty--Kim locus eventually equals the set of integral points; see~\cite[Conjecture~3.1 and \S8.1]{BDKWTateShafarevich}.} Successful computations include $X_{\mathrm{ns}}^+(13)$ by quadratic Chabauty~\cite[Corollary~1.3]{BALAKRISHNAN_DOGRA} and the thrice-punctured line by motivic Chabauty~\cite[Theorem~1.2]{CDC}. It is thus very natural to ask: when does~\eqref{eq:kim_dim_ineq} hold? In~\cite[\S3, Observations~1--2 and Proposition~2]{KIM}, Kim deduced his dimension inequality for all smooth quasi-projective hyperbolic curves over $\mathbb{Z}[1/S]$, assuming either one of the Bloch--Kato or Fontaine--Mazur conjectures. In that same paper, Kim formulated the following conjecture.

\begin{conjecture}[{\cite[Conjecture~1, p.~96]{KIM}}]
\label{conj:kim-dimension}
Notations as above. There exists some $n\ge 1$, possibly large, for which \eqref{eq:kim_dim_ineq} holds.
\end{conjecture}
Kim's conjecture was proved by Kim and others for various special families of hyperbolic curves. In~\cite[\S3]{KIM2}, Kim proved his conjecture for the thrice-punctured line, and in~\cite[\S4, pp.~127--129]{KIM}, for punctured CM elliptic curves of rank one over $\mathbb{Q}$. Kim subsequently removed the rank restriction for CM elliptic curves, using commutative Iwasawa theory, for a certain quotient of the unipotent fundamental group, under certain additional conditions on the auxiliary prime $p$~\cite[Theorem~0.1]{KimCMElliptic}. This result was later extended by Kim and Coates in~\cite[Corollary~0.2]{COATES_KIM} to arbitrary smooth projective curves whose Jacobians are isogeneous to a product of CM abelian varieties. In~\cite[\S9]{EllenbergHastSolvable}, Ellenberg and Hast established Kim's condition for solvable covers of $\mathbb{P}^1$ of genus at least $2$ over $\mathbb{Q}$. In particular, this establishes the validity of Kim's method for superelliptic curves.

\paragraph{Our main results}
We now specialize to the case where $C/\mathbb{Q}$ is smooth projective geometrically connected of genus $g\ge 2$, and with $p$ odd. Denote by $\mathfrak{n} := \Lie(U^\et)$ the pro-nilpotent Lie algebra of $U^\et$, and by $L := \gr\,\mathfrak{n}$ its associated graded Lie algebra. Write $L_n$ for the $n$'th graded piece of $L$, regarded as a finite dimensional $\Qp$-vector space equipped with a $G_{\mathbb{Q},\Sigma}$-action. Throughout, we write $b_n := \dim_\Qp L_n$. We shall prove, using Theorem~\ref{thm:h2-bound}, the following power-saving bound.
\begin{theorem}\label{thm:power_saving}
Notations as above, there is a constant $A := A(p,S,\mathcal{C})$, such that for all $n\ge 1$, one has:
\[
\dim H^2(G_{\mathbb{Q},\Sigma}, L_n) < A\cdot b_n\cdot n^{-\frac{1}{D+1}}.
\]
\end{theorem}
As a corollary, we establish a new case of Kim's conjecture.
\begin{corollary}\label{cor:kim-conj}
The Chabauty--Kim condition holds for all $n$ large enough in the smooth projective geometrically connected case, for any odd auxiliary prime $p$. Hence, Conjecture~\ref{conj:kim-dimension} is true in this setting.
\end{corollary}
In particular, the Chabauty--Kim method recovers the Mordell conjecture (now Faltings' theorem~\cite{FALTINGS}) for smooth projective geometrically connected curves over $\mathbb{Q}$.

Theorem~\ref{thm:power_saving} is proved in~\S\ref{ssec:character-estimates}, and Corollary~\ref{cor:kim-conj} is established in~\S\ref{ssec:kim-conj-proof}.
\subsection{Proof Ideas}
We now discuss the key ingredients involved in the proofs of Theorems~\ref{thm:h2-bound} and~\ref{thm:power_saving}.

\subsubsection{Our bound on \texorpdfstring{$H^2(G_{\mathbb{Q},\Sigma}, V_{\lambda,s})$}{H2(G(Q,Sigma), V(lambda,s))}}\label{sssec:h2-bound-proof-idea}
Theorem~\ref{thm:h2-bound} bounds the dimension of $H^2(G_{\mathbb{Q},\Sigma}, V_{\lambda,s})$, outside a finite set of bad weights, by a constant multiple of
\[
\left\{d_\lambda,\dfrac{(1 + ||\lambda||)^{D-g-1}}{d_\lambda}\right\}.
\]
Our strategy is heavily inspired by the work of Coates and Kim in~\cite{COATES_KIM}, and is closely related to Fu's bounds~\cite[Theorem~1.4]{FU} for the homology of finitely generated Iwasawa modules tensored with algebraic representations, for products of principal congruence subgroups of $\mathrm{SL}_2(\Zp)$.

Let $\mathcal{A}/\mathbb{Q}$ denote our ambient abelian variety, which plays the role of the Jacobian of our curve of interest in our later applications to Chabauty--Kim theory. In~\S\ref{ssec:iwasawa_setup} we establish the ambient Iwasawa theoretic setting. Briefly, for each $n\ge 1$, let $K_n := \mathbb{Q}(\mathcal{A}[p^n])$ denote the smallest extension of $\mathbb{Q}$, containing the coordinates of all $p^n$-torsion points of $\mathcal{A}$. Write $K_{\infty} := \cup_{n\ge 1}\mathbb{Q}(\mathcal{A}[p^n])$ for their colimit. In~\S\ref{ssec:iwasawa_setup} we fix some $c\ge 1$, for which $\Gamma := \Gal(K_\infty/K_c)$ is a uniform pro-$p$ group, which decomposes into a product of the form $H_c\times Z_c$, with
\[
H_c := \Sp_{2g,\Zp}(\Zp)\cap \mathrm{Ker}\left(\Gamma\to \GSp_{2g,\Zp}(\Zp/p^c)\right),\qquad Z_c := (1 + p^c\Zp)\cdot I,
\]
where $I\in \Gamma\subseteq \GSp_{2g,\Zp}(\Zp)$ is the identity element. We then fix the integral and rational Iwasawa algebras
\[
\Omega := \Zp[[\Gamma]],\qquad \Lambda := \Omega[1/p].
\]
It is an immediate consequence of the Hochschild--Serre spectral sequence that for any finite dimensional $p$-adic $G_{\mathbb{Q},\Sigma}$-representation $W$,
\[
\dim_\Qp H^2(G_{\mathbb{Q},\Sigma}, W)\le \dim_\Qp H^2(G_{F,\Sigma_F}, W),
\]
where $F$ is any finite Galois extension of $\mathbb{Q}$ and $\Sigma_F$ is the set of places of $F$ lying above $\Sigma$. We apply this inequality with $F = K_c := \mathbb{Q}(\mathcal{A}[p^c])$, and fix the Iwasawa module $X := H^2_{\mathrm{Iw}}(K_\infty/F; \Zp(1))$ (see Definition~\ref{def:x}). We then prove that $X$ is finitely generated and torsion; see Propositions~\ref{prop:fg_X} and~\ref{prop:torsion_X}.

In Theorem~\ref{thm:galois_iwasawa} of \S\ref{ssec:key_reduction}, we prove that for any finite dimensional $p$-adic $G_{F,\Sigma_F}$-representation $W$, for which the Galois action factors through $\Gamma$, one has
\[
\dim_\Qp H^2(G_{F,\Sigma_F}; W) = \dim_\Qp \Hom_\Gamma(X, W^\vee(1)),
\]
where $W^\vee(1) := \Hom(W, \Qp(1))$, and $\Hom_\Gamma(\bullet, \bullet)$ denotes the subspace of $\Gamma$-equivariant homomorphisms. Applying this with $W = V_{\lambda,s}$, Theorem~\ref{thm:h2-bound} reduces to bounding $\dim_\Qp \Hom_\Gamma(X, V_{\lambda,2-s})$. Here we used that if $W$ is pure of weight $s$, its Tate dual is pure of weight $2-s$, and the fact that representations of the symplectic group are self-dual (due to the symplectic form).

The next idea is also taken from Coates and Kim. Since $X$ is finitely generated we can pick a finite set of generators $x_1,\dots,x_q$, and since it is torsion, we can choose nontrivial left annihilators $0\neq \delta_1,\dots,\delta_q\in \Lambda$, such that $\delta_ix_i = 0$ for every $1\le i\le q$. Evaluation at $x_i$ induces an embedding; see~\eqref{eq:hom-to-ker},
\[
\Hom_\Gamma(X, V_{\lambda,2-s}) \subseteq \bigoplus_{i=1}^q \mathrm{Ker}\,\delta_i|_{V_{\lambda,2-s}}.
\]
Theorem~\ref{thm:h2-bound} is actually a pair of bounds: a trivial bound of order $d_\lambda$, valid for every weight, and a nontrivial bound, which is valid outside a finite set of weights. To see the trivial bound, observe that $\dim_\Qp \mathrm{Ker}\,\delta_i|_{V_{\lambda,2-s}}$ is trivially bounded by $d_\lambda = \dim_\Qp V_{\lambda,2-s}$ for every $i$, which implies that
\[
\dim_\Qp \Hom_\Gamma(X, V_{\lambda,2-s})\le q\cdot d_\lambda.
\]
As for the nontrivial bound, it suffices to bound $\mathrm{Ker}\,\delta|_{V_{\lambda,2-s}} \ll \dfrac{(1 + ||\lambda||)^{D-g-1}}{d_\lambda}$ for an arbitrary nontrivial $\delta$. As our bound is independent of $s$, we may as well treat $\mathrm{Ker}\,\delta|_{V_{\lambda,s}}$. The nontrivial bound is the main focus of~\S\ref{sec:bounds}. The rest of this segment is divided into 5 steps, reflecting the structure of our proof of the nontrivial bound, which is eventually carried out in~\S\ref{sssec:nontriv_proof}.

\textbf{Step I.} In Proposition~\ref{prop:wpt} and Corollary~\ref{cor:red-2}, we argue via a noncommutative Weierstrass preperation-type of an argument, that outside of finitely many values of $s$, $\mathrm{Ker}\,\delta|_{V_{\lambda,s}}$ coincides with $\mathrm{Ker}\,\eta|_{V_{\lambda,s}}$, where $\eta \in \Zp[[H_c]][[Z_c]]$ is a possibly different Iwasawa operator, with the property that $\eta(T) \neq 0$ for any $T\in p\Zp$. Since the central coordinate $T$ acts on $V_{\lambda,s}$ via $t_s := (1+p^c)^s - 1 \in p\Zp$, $\eta$ acts on $V_{\lambda,s}$ via $0\neq \eta(t_s) \in \Zp[[H_c]]$, i.e.~through specializing $T$ to $t_s$.

\textbf{Step II.} Let $\mathfrak{h} := \Lie(\Sp_{2g,\Zp})$, and let $\mathcal{U} := \mathcal{U}(\mathfrak{h})$ denote its (integral) universal enveloping algebra. The universal enveloping algebra is equipped with its ascending PBW filtration $\mathcal{F}_\bullet$, where $\mathcal{F}_d\,\mathcal{U}$ consists of those elements of $\mathcal{U}$, which can be expressed via products of at most $d$ elements of $\mathfrak{h}$. For a $\Zp$-module $M$, write $\widehat{M} := \varprojlim_n M/p^nM$ for its pro-$p$ completion. In~\S\ref{ssec:univ_env} we introduce the Lazard comparison map, 
\[
\Zp[[H_c]]\to \widehat{\mathcal{U}}
\]
that factors $\eta(t_s)$'s action on $V_{\lambda,s}$. Denote by $\mathcal{Z}\subseteq \mathcal{U}$ the center of the universal enveloping algebra. In~\S\ref{ssec:univ_env}, we prove (see Claim~\ref{claim:center_presentation}) that
\[
\mathcal{Z} \cong \Zp[c_1,\dots,c_g],
\]
with $c_i \in \mathcal{F}_{2i}\,\mathcal{U}$. The center of $\mathcal{U}$ acts on $V_{\lambda,s}$ through the Harish--Chandra character (see Definition~\ref{def:central_character}),
\[
\zeta_\lambda: \mathcal{Z}\to \Zp,
\]
and so the action of $\mathcal{U}$ on $V_{\lambda,s}$ factors through the central reduction,
\[
R_\lambda := \mathcal{U}/(c_i - \zeta_\lambda(c_i): 1\le i\le g).
\]
Similarly, the action of $\widehat{U}$ (hence that of $\eta(t_s)$) factors through $\widehat{R_\lambda}$. The affinoid Verma faithfulness theorem of Ardakov and Wadsley as in~\cite[Theorem~5.4]{ArdakovWadsleyVerma}, see Corollary~\ref{cor:eta_reduction}, implies that the composed map 
\[
\Zp[[H_c]]\to \widehat{\mathcal{U}}\to \widehat{R_\lambda}
\]
is injective, and so the image of $0\neq \eta(t_s)$, denoted $\eta_{\lambda,s} \in \widehat{R_\lambda}$, is nontrivial as well.

\textbf{Step III.}
For every $(\lambda,s)$-pair, write $\theta_{\lambda,s}$ for the quotient of $\eta_{\lambda,s}$ by the largest possible power of $p$, such that $p^{-1}\theta_{\lambda,s}\not\in \widehat{R_\lambda}$. In particular, the mod-$p$ reduction of $\theta_{\lambda,s}$ satisfies
\[
0\neq \overline{\theta_{\lambda,s}} \in \overline{R_\lambda} := R_\lambda/pR_\lambda.
\]
In Lemma~\ref{lem:key_lemma}, we establish the existence of an integer $d\ge 1$, independent of $s$ and $\lambda$, such that 
\begin{equation}\label{eq:theta-d}
\overline{\theta_{\lambda,s}} \in \mathcal{F}_d\overline{R_\lambda},  
\end{equation}
where we equip $\overline{R_\lambda}$ with the quotient PBW filtration.

\textbf{Step IV.}
Fix a Chevalley $\Zp$-lattice $\mathcal{V}_{\lambda,s}\subseteq V_{\lambda,s}$, and write 
\[
N_{\lambda} := \mathrm{Im}\left(\mathcal{U}\to \End_{\Zp}(\mathcal{V}_{\lambda,s})\right)
\]
for the image of $\mathcal{U}$ in the endomorphism algebra of $\mathcal{V}_{\lambda,s}$. It follows, as a corollary of Jacobson's density theorem (see Corollary~\ref{cor:jacobson}), that $N_{\lambda}$ is a lattice in $\End_{\Qp}(V_{\lambda,s})$, in the sense that
\[
N_{\lambda}[1/p] = \End_\Qp V_{\lambda,s}.
\]
Write $\theta$ for $\theta_{\lambda,s}$'s image in $N_{\lambda}$. Then,
\[
\rk_\Zp \dfrac{N_\lambda}{N_\lambda\theta} = \dim_\Qp \dfrac{\End_\Qp(V_{\lambda,s})}{\End_\Qp(V_{\lambda,s})\theta} = \dim_\Qp \mathrm{Ker}\,\theta|_{V_{\lambda,s}}\cdot d_\lambda.
\]
By construction, the actions of $\delta$ and $\theta$ on $V_{\lambda,s}$ coincide, up to a nonzero scalar. Therefore, Theorem~\ref{thm:h2-bound} reduces to showing that
\[
\rk_\Zp \dfrac{N_\lambda}{N_\lambda\theta} \ll (1 + ||\lambda||)^{D-g-1}.
\]

As for the latter, in Claim~\ref{claim:jacobson} we show that the surjection $\mathcal{U}[1/p]\twoheadrightarrow \End_\Qp(V_{\lambda,s})$ is already attained at PBW filtration length $2D\lceil ||\lambda||\rceil$. In parciular, since $\mathcal{U}$'s action on $V_{\lambda,s}$ factors through the central reduction, we have a surjection of $\Zp$-modules, 
\[
\mathcal{F}_{2D\lceil ||\lambda||\rceil + d}R_\lambda\twoheadrightarrow N_\lambda.
\]
Writing $\overline{N_\lambda}$ and $\overline{\theta}$ for the mod-$p$ reductions of $N_\lambda$ and $\theta$. By Nakayama's lemma, 
\[
\rk_\Zp \dfrac{N_\lambda}{N_\lambda\theta} \le \dim_\Fp \dfrac{\overline{N_\lambda}}{\overline{N_\lambda}\cdot\overline{\theta}}.
\]
Hence we have a surjection
\[
\mathcal{F}_{2D\lceil ||\lambda||\rceil + d}\overline{R_\lambda}\twoheadrightarrow \overline{N_\lambda} \twoheadrightarrow \dfrac{\overline{N_\lambda}}{\overline{N_\lambda}\cdot\overline{\theta}}.
\]
The kernel of the surjection contains $\mathcal{F}_{2D\lceil ||\lambda||\rceil + d}(\overline{R_\lambda}\cdot \overline{\theta_{\lambda,s}})$, which contains $\mathcal{F}_{2D\lceil ||\lambda||\rceil}(\overline{R_\lambda})\cdot \overline{\theta_{\lambda,s}}$, since $\overline{\theta_{\lambda,s}}\in \mathcal{F}_d\overline{R_\lambda}$ by~\eqref{eq:theta-d}. In Proposition~\ref{prop:r-mu-reduction-domain} we show that $\overline{R_\lambda}$ is a domain, and therefore
\[
\dim_\Fp \mathcal{F}_{2D\lceil ||\lambda||\rceil}(\overline{R_\lambda})\cdot \overline{\theta_{\lambda,s}} \ge \dim_\Fp \mathcal{F}_{2D\lceil ||\lambda||\rceil}\overline{R_\lambda}.
\]
It then follows from the first isomorphism theorem, that
\begin{equation}\label{eq:naka-bound}
\dim_\Fp \dfrac{\overline{N_\lambda}}{\overline{N_\lambda}\overline{\theta}} \le \dim_\Fp\mathcal{F}_{2D\lceil ||\lambda||\rceil + d}\overline{R_\lambda} - \dim_\Fp \mathcal{F}_{2D\lceil ||\lambda||\rceil}\overline{R_\lambda}.
\end{equation}

\textbf{Step V.}
The right-hand side of~\eqref{eq:naka-bound} is precisely $\sum_{m=2D\lceil ||\lambda||\rceil + 1}^{2D\lceil ||\lambda||\rceil + d}\dim_\Fp \gr_m \overline{R_\lambda}$, and a result of Ardakov and Wadsely allows us to identify $\gr \overline{R_\lambda}$ with a certain graded algebra, denoted $S(\mathfrak{h}_k)/(\overline{q_1},\dots,\overline{q_g})$, where $S(\mathfrak{h}_k)$ is the symmetric algebra on $\mathfrak{h}_k$, and the $\overline{q_i}$ are the mod-$p$ reductions of the principal symbols of the $c_i$ in $S(\mathfrak{h})$. Importantly, the algebra $\gr \overline{R_\lambda}$ is graded, generated in degree $1$, and has Krull dimension $D-g$. A basic Hilbert-series argument then shows that
\[
\dim_\Fp \gr_m \overline{R_\lambda} = O(m^{D-g-1}).
\]
The crucial point here is that although the central reduction is dependent on $\lambda$, the dimensions of its graded pieces are independent of $\lambda$. Summing these bounds over the graded pieces with $2D\lceil ||\lambda||\rceil + 1\le m\le 2D\lceil ||\lambda||\rceil + d$, the nontrivial bound follows. This last estimate is carried out in Proposition~\ref{prop:graded_bound}.

\subsubsection{Our bound on \texorpdfstring{$H^2(G_{\mathbb{Q},\Sigma}, L_n)$}{H2(G(Q,Sigma), Ln)}}
Theorem~\ref{thm:power_saving} is proved in \S\ref{sec:ck-dim-ineq}, through averaging. With $L_1 = V := H_1^\et(\overline{C}, \Qp)$, where $\overline{C} = C\otimes \overline{\mathbb{Q}}$ is the geometric generic fibre of $C$. Observe that $L_1$ is pure of weight $1$ since $C$ is smooth and projective. It follows that $L_n$, as a quotient of $L_1^{\otimes n}$, is pure of weight $n$. Decomposing $L_n$ into a direct sum of $\Sp_{2g,\lambda}$-irreducibles, we have
\[
L_n = \bigoplus_{\lambda}V_{\lambda,n}^{\oplus m_n(\lambda)},
\]
where the integer $m_n(\lambda)\in\mathbb{Z}_{\ge 0}$ denotes the multiplicity of the irreducible representation $V_{\lambda,n}$ in $L_n$. In~\S\ref{sec:ck-dim-ineq}, we bound the following representation theoretic invariants.
\begin{enumerate}
  \item In Proposition~\ref{prop:mn-bound} we bound $m_n(\lambda) \ll d_\lambda \cdot b_n\cdot n^{-D/2}$.
  \item In Proposition~\ref{prop:second_moment} we bound $\sum_{\lambda} m_n(\lambda)\cdot d_\lambda\cdot ||\lambda||^2 \ll b_n\cdot n$.
\end{enumerate}
The key idea is to move to a real compact form of the symplectic group, where these estimates may be obtained via analysis, using the character theory of the compact group. The relevant compact group here is $\mathrm{USp}(2g)$, the unitary symplectic group, and the important fact is that $L_n$'s irreducible decomposition persists through extension of scalars to $\mathbb{C}$, and then upon restriction to $\mathrm{USp}(2g)$. Regarding $L_n$ now as a representation of the unitary symplectic group, Filip~\cite[Proposition~4.4, Eq.~(4.8)]{FilipSurfaceLie} (see Theorem~\ref{thm:filip_lie}) provides us with an explicit expression for the character of $L_n$, denoted $\chi_n$. Namely, for all $u\in \mathrm{USp}(2g)$, one has 
\[
\chi_n(u) = \dfrac{1}{n}\sum_{d|n}\mu(n/d)P_d(\chi_V(u^{n/d})),
\]
where $\chi_V$ is the character of the standard representation of $\mathrm{USp}(2g)$, and
\[
P_d(x) = \left(\dfrac{x + \sqrt{x^2 - 4}}{2}\right)^d + \left(\dfrac{x - \sqrt{x^2 - 4}}{2}\right)^d
\]
is a rescaled Chebyshev polynomial. 

\textbf{The bound on $m_n(\lambda)$} is obtained by using the orthogonality of characters, i.e.~through bounding an integral expression for $m_n(\lambda)$ given explicitly as 
\[
\int_{\mathrm{USp}(2g)}\chi_n\cdot \overline{\chi_\lambda}dk,
\]
where $\chi_\lambda$ is the character of $V_{\lambda,n}$ and $dk$ is a probability Haar measure. Weyl's unitary trick provides the bound $|\chi_\lambda(u)| \ll d_\lambda$, and the bound $m_n(\lambda) \ll d_\lambda \cdot b_n\cdot n^{-D/2}$ boils down to the estimate $\int_{\mathrm{USp}(2g)}|\chi_n|dk \ll b_n\cdot n^{-D/2}$, which is carried out within the proof of Proposition~\ref{prop:mn-bound} using elementary analysis. 

\textbf{The second-moment estimate} $\sum_{\lambda} m_n(\lambda)\cdot d_\lambda\cdot ||\lambda||^2 \ll b_n\cdot n$ is obtained by considering the Casimir's action on $\chi_n$ at the identity, observing that it is bounded below by the second-moment sum $\sum_{\lambda} m_n(\lambda)\cdot d_\lambda\cdot ||\lambda||^2$, and bounded above by the absolute values of the action of the Casimir on the summands in Filip's formula, evaluated at the identity.

\textbf{Combining the pieces.} In~\S\ref{sec:ck-dim-ineq} we complete the proof of Theorem~\ref{thm:power_saving} as follows. Assuming $n$ is large enough so that $n$ is not in the set of bad weights $E$ from Theorem~\ref{thm:h2-bound}, applied to $\mathcal{A} = \mathrm{Jac}(C)$, we bound
\[
\dim_\Qp H^2(G_{\mathbb{Q},\Sigma}, L_n) = \sum_\lambda m_n(\lambda)\dim_\Qp H^2(G_{\mathbb{Q},\Sigma}, V_{\lambda,n})
\]
as follows. First, we partition the sum into contributions from $||\lambda|| < R\sqrt{n}$ and $||\lambda|| \ge R\sqrt{n}$, with $R$ an optimization parameter, whose choice is to be fixed in hindsight. As for the small $\lambda$'s, we apply the nontrivial bound,
\[
\dim_\Qp H^2(G_{\mathbb{Q},\Sigma}, V_{\lambda,n}) \ll \dfrac{(1 + ||\lambda||)^{D-g-1}}{d_\lambda}.
\]
Plugging in our bound on $m_n(\lambda)$, the contribution of the small $\lambda$'s is bounded by
\[
b_n\cdot n^{-D/2}\sum_{||\lambda|| < R\sqrt{n}}(1 + ||\lambda||)^{D-g-1}.
\]
The sum being taken over the lattice points within a fixed sector of the $g$-dimensional ball of radius $R\sqrt{n}$, hence a sum over $O((R\sqrt{n})^g)$-summands, each contributes at most a constant multiple of $(R\sqrt{n})^{D-g-1}$, and therefore the small weight contribution is bounded above by $O\left(\dfrac{b_nR^{D-1}}{\sqrt{n}}\right)$.

As for the larger weights, we use the trivial bound, 
\[
\dim_\Qp H^2(G_{\mathbb{Q},\Sigma}, V_{\lambda,n}) \ll d_\lambda.
\]
The large-weight sum is bounded by $\sum_{||\lambda||\ge R\sqrt{n}}m_n(\lambda)d_\lambda$, which we may further bound above, by the well-known Rankin trick, via
\[
\dfrac{1}{R^2n}\sum_{||\lambda||\ge R\sqrt{n}}m_n(\lambda)d_\lambda||\lambda||^2.
\]
Completing the sum and applying our second-moment estimate, we bound the large-weight contribution above by $b_nR^{-2}$.

We therefore obtain the upper bound,
\[
\dim_\Qp H^2(G_{\mathbb{Q},\Sigma}, L_n) \ll \dfrac{b_nR^{D-1}}{\sqrt{n}} + b_nR^{-2},
\]
which is optimal when the contributions of the two terms are balanced. This is achieved by fixing $R = n^{\frac{1}{2(D+1)}}$, which gives rise to the desired upper bound,
\[
\dim_\Qp H^2(G_{\mathbb{Q},\Sigma}, L_n) \ll b_n\cdot n^{-\frac{1}{D+1}}.
\]

\subsection{Plan of Paper}
The paper is structured as follows.
\begin{itemize}
  \item In~\S\ref{sec:galois_iwasawa}, we fix the Iwasawa theoretic setup (in \S\ref{ssec:iwasawa_setup}), carry out the key reduction step in the proof of Theorem~\ref{thm:h2-bound} (in \S\ref{ssec:red2}), reducing the study of Galois cohomology groups to the study of kernels of certain integral Iwasawa operators, and prove the trivial bound on the latter (in \S\ref{ssec:red2}).
  \item In~\S\ref{sec:bounds}, we introduce relevant background material regarding the universal enveloping algebra of the symplectic Lie algebra, and its central reduction (in \S\ref{ssec:univ_env}), we then prove the nontrivial bound on the dimensions of our kernels (in \S\ref{ssec:nontriv}), and complete the proof of Theorem~\ref{thm:h2-bound} in \S\ref{ssec:h2-bound}.
  \item In~\S\ref{sec:ck-dim-ineq}, we prove our power saving bound, i.e.~Theorem~\ref{thm:power_saving}, in~\S\ref{ssec:character-estimates}, and complete the proof of Kim's conjecture in~\S\ref{ssec:kim-conj-proof}.
\end{itemize}

\subsection*{Acknowledgements}
The author wishes to thank Qingyuan Bai, Rok Gregoric, Alexandre Maksoud, Ostap Niezhenskyi, Wyatt Reeves, and Germ\'an Stefanich for general discussions on the ambient theory and on various aspects of the problem. The author extends special thanks to Peter Scholze for discussing a preliminary draft, and for suggesting the use of the Weyl character formula in order to compare the trivial and nontrivial bounds (see Remark~\ref{rem:weyl-bound-comparison}).
\subsection*{AI Declaration}
The general strategy is the author's, or perhaps of Coates and Kim~\cite[Introduction, penultimate paragraph]{COATES_KIM}, and the reduction of dimensions of Galois cohomology groups to dimensions of kernels of Iwasawa operators has been carried out independently by the author about a year prior to this draft. The development of the nontrivial bound of Theorem~\ref{thm:equiv-hom-bound}, as well as its application to the power-saving bound of Theorem~\ref{thm:power_saving} was carried out through a lengthy back-and-forth interaction with GPT 6 Astra (Open AI) on Ultra mode. Subsequently, the author generated a dense shorter draft ($33$-page long) for a careful examination of the argument with the help of Astra. 

Albeit only stated briefly in~\S\ref{sssec:h2-bound-proof-idea}, some of the ideas taking part in the proof of the nontrivial bound, which came up in the author's discussions with Astra, already appear in a paper of Weibo Fu~\cite[Theorem 1.4]{FU} and played a key role in the proof of Theorem~\ref{thm:equiv-hom-bound}, albeit not using Fu's results directly.

This draft was later written by the author without any use of LLMs, including without consulting the older version, except perhaps for discussions pertaining standard ambient theory, and gathering the relevant references. Some aspects in this draft's proof of Theorem~\ref{thm:equiv-hom-bound} are different, and less technical, than those initially carried out in the preliminary draft.

\section{From Galois cohomology to kernels of Iwasawa operators}\label{sec:galois_iwasawa}
This chapter is divided as follows:
\begin{enumerate}
  \item In \S\ref{ssec:iwasawa_setup}, we fix our ongoing Iwasawa theoretic setting.
  \item In \S\ref{ssec:key_reduction}, we define the Iwasawa module $\mathcal{X} := H^2_{\mathrm{Iw}}(K_\infty/F; \Zp(1))$ as well as its rational version, $X := \mathcal{X}[1/p]$, and establish their key properties.
  \item In \S\ref{ssec:red2}, we carry out the key reduction step in the proof of Theorem~\ref{thm:h2-bound}, and prove the trivial bound $\dim_\Qp H^2(G_{\mathbb{Q},\Sigma},V_{\lambda,s}) \ll d_\lambda$.
\end{enumerate}

\subsection{Iwasawa theoretic setup}\label{ssec:iwasawa_setup}
We briefly repeat our setting from the introduction. Let $\mathcal{A}/\mathbb{Q}$ denote an abelian variety with good reduction outside a finite set of primes $S$, and let $p\not\in S$ denote an arbitrary \emph{odd} auxiliary prime. Let $V = T_p\mathcal{A}\otimes_{\mathbb{Z}_p}\mathbb{Q}_p$ denote its (rationalized) $p$-adic Tate module, and let
\[
\rho_{\mathcal{A},p}: G_{\mathbb{Q},\Sigma}\to \GSp_{2g}(\mathbb{Q}_p)
\]
denote its associated Galois representation, where $\Sigma := S\cup\{p,\infty\}$ and $G_{\mathbb{Q},\Sigma}$ is the Galois group of the maximal extension of $\mathbb{Q}$ unramified outside of $\Sigma$. Write 
\[
\Gamma_{\mathcal A,p} := \mathrm{Im}(\rho_{\mathcal{A},p})\subseteq \GSp_{2g}(\mathbb{Z}_p)
\]
for $\rho_{\mathcal{A},p}$'s image, equipped with its $p$-adic Lie group structure. Let
\[
\mathbf{G} := \overline{\Gamma_{\mathcal A,p}},\qquad \mathbf{R} := \mathbf{G}^\circ
\]
denote the Zariski closure of $\Gamma_{\mathcal A,p}\subseteq \GSp_{2g,\mathbb{Q}_p}$, and its connected component of the identity, respectively. 

\begin{definition}
For every $c\ge 1$, define the level $c$ congruence subgroup of $\Gamma_{\mathcal A,p}$, denoted $\Gamma_{\mathcal A,p}(c)$, via:
\[
\Gamma_{\mathcal A,p}(c) := \mathrm{Ker}\left(\Gamma_{\mathcal A,p}\to \GSp_{2g}(\mathbb{Z}_p/p^c)\right).
\]
\end{definition}

Our next theorem will make use of the following pair of facts, due to Bogomolov.
\begin{theorem}[{\cite[Theorem~1.1(c)]{PinkReductionPoint}}]
\label{thm:bogomolov_open}
$\Gamma_{\mathcal{A},p}$ is open in $\mathbf{G}(\mathbb{Q}_p)$ in the $p$-adic topology.
\end{theorem}
\begin{theorem}[{\cite[Theorem~3]{BogomolovFiniteOrder}}]
\label{thm:bogomolov_center}
Set $\mathfrak{r} := \Lie(\mathbf{R})$, then $\mathfrak{r}\supseteq \Qp\cdot I$.
\end{theorem}

\begin{theorem}\label{thm:gamma_decomposition}
Let $H_c := \Gamma_{\mathcal{A},p}(c)\cap \Sp_{2g}(\mathbb{Z}_p)$, and $Z_c := (1 + p^c\Zp)\cdot I$. Then, for $c\gg 1$, one has:
\begin{enumerate}
  \item $\Gamma_{\mathcal{A},p}(c)$ is uniform.
  \item $\Gamma_{\mathcal{A},p}(c) = Z_c\times H_c$.
\end{enumerate}
\end{theorem}

\begin{proof}
\textbf{Uniformity.} 
Since every $p$-adic Lie group contains a uniform open subgroup (see~\cite[Theorems 8.29-8.31]{DDMSAnalyticProPGroups}), and the congruence subgroups $\Gamma_{\mathcal{A},p}(c)$ form a cofinal family of open neighborhoods of the identity, for $c$ large enough the subgroups $\Gamma_{\mathcal{A},p}(c)$ are all uniform. 

\textbf{Decomposition.}
We first argue that for all $c$ large enough,
\[
\Gamma_{\mathcal{A},p}(c)\supseteq Z_c := (1+p^c\Zp)\cdot I.
\]
Write $M_{2g}(\Zp)$ for the ring of $2g$ by $2g$ matrices with entries in $\Zp$, and fix $\mathfrak{r}_{\Zp} := \mathfrak{r}\cap M_{2g}(\Zp)$, where $\mathfrak{r} := \Lie(\mathbf{R})$. By Theorem~\ref{thm:bogomolov_center}, $\mathfrak{r}\supseteq \Qp\cdot I$. Therefore, it suffices to prove that $\log \Gamma_{\mathcal{A},p}(c) = p^{c}\mathfrak{r}_{\Zp}$ for $c$ large enough. We argue via a two-sided inclusion. First, arguing that for large enough $c$, $\log \Gamma_{\mathcal{A},p}(c) \subseteq p^{c}\mathfrak{r}_{\Zp}$. 

Since $\mathbf{R}(\mathbb{Q}_p)$ is open in $\mathbf{G}(\Qp)$, it contains an open neighborhood of the identity. Since the $\{\Gamma_{\mathcal{A},p}(c): c\ge 1\}$ form a cofinal family of neighborhoods of the identity, for all $c$ large enough, one has
\[
\Gamma_{\mathcal{A},p}(c)\subseteq \mathbf{R}(\mathbb{Q}_p).
\]
Taking formal logarithms, while recalling that by construction $\log \Gamma_{\mathcal{A},p}(c)\subseteq p^c\cdot M_{2g}(\Zp)$, the inclusion follows. Next, we seek to show that for large enough $c$, $\log \Gamma_{\mathcal{A},p}(c) \supseteq p^{c}\mathfrak{r}_{\Zp}$.

The formal logarithm map defines an analytic homeomorphism between a sufficiently small neighborhood of the identity in $\mathbf{R}(\mathbb{Q}_p)$ and a sufficiently small neighborhood of $0$ in $\mathfrak{r}$. Since by Theorem~\ref{thm:bogomolov_open}, $\Gamma_{\mathcal{A},p}\subseteq \mathbf{G}(\mathbb{Q}_p)$ is open, there exists a constant $c_1$, such that:
\[
p^{c_1}\mathfrak{r}_{\Zp}\subseteq \log \Gamma_{\mathcal{A},p}.
\]
From here on out, we fix $c$ large enough so that both inclusions are satisfied simultaneously. Hence, we may assume that $\Gamma_{\mathcal{A},p}(c)\supseteq Z_c := (1+p^c\Zp)\cdot I$. 

Our next goal is to establish the decomposition $\Gamma_{\mathcal{A},p}(c) = H_c\times Z_c$, with $H_c := \Gamma_{\mathcal{A},p}(c)\cap \Sp_{2g}(\mathbb{Z}_p)$. We first prove that the subgroup generated by such products generates $\Gamma_{\mathcal{A},p}(c)$. Let 
\[
\omega: V\otimes_{\Qp} V\to \Qp(1)
\]
denote the Weil pairing on $V$, regarded as the symplectic bilinear form for the ambient symplectic group $\Sp_{2g}(\Qp)$. Write $\nu: \GSp_{2g}(\Qp)\to \Qp^\times$ for the similitude character, which is defined, for $g\in \GSp_{2g}(\Zp)$, to be the unique $p$-adic number $\nu(g)\in \Qp^\times$ with the property that
\[
\omega(gu,gv) = 
\nu(g)\cdot \omega(u,v),\qquad \text{for all }u,v\in V.
\]
Let $M\in \Gamma_{\mathcal{A},p}(c)$ be arbitrary, and write $t := \nu(M) \in 1 + p^c\Zp$ for its similitude. We shall prove that $M\in H_c\times Z_c$. Since $p$ is odd, every element of $1 + p^c\Zp$ admits a unique square root, by Hensel's lemma. Let $s\in 1 + p^c\Zp$ denote $t$'s unique square root, we claim that $s^{-1}M \in H_c$. To see this, observe that $H_c$ is precisely the kernel of $\nu|_{\Gamma_{\mathcal{A},p}(c)}$, and that for every element $x\cdot I\in Z_c$, with $x\in 1 + p^c\Zp$, one has $\nu(x\cdot I) = x^2$. It follows that 
\[
\nu(s^{-1}M) = \nu(s^{-1}I)\cdot \nu(M) = s^{-2}\cdot t = 1,
\]
and hence that $s^{-1}M \in H_c$. Hence also that $M = sI\cdot s^{-1}M \in Z_c\cdot H_c$.

To observe that $\Gamma_{\mathcal{A},p}(c)$ indeed admits a decomposition into $Z_c\times H_c$, since $Z_c$ is central, it suffices to observe that $Z_c\cap H_c = \{I\}$, i.e.~that the intersection consists only of the identity. Indeed, let $g\in Z_c\cap H_c$. Since $g\in Z_c$, $g = t\cdot I$, with $t\in 1 + p^c\Zp$. Therefore, 
\[
\nu(g) = \nu(t\cdot I) = t^2.
\]
Since $g\in H_c$ and $\nu(H_c) = 1$, it follows that $t^2 = 1 \in 1 + p^c\Zp$, however, since every element of $1+p^c\Zp$ admits a unique square root and $t = 1$ is such a square root, we must have $g = I$, which completes the proof.
\end{proof}

From this point onward, we fix a $c\ge 1$ as guaranteed in Theorem~\ref{thm:gamma_decomposition}. We abbreviate $\Gamma := \Gamma_{\mathcal{A},p}(c)$ and $Z := Z_c$. 

\begin{notation}
For each $a\ge 1$, write $K_a := \mathbb{Q}(\mathcal{A}[p^a])$ (resp. $K_\infty := \mathbb{Q}(\mathcal{A}[p^\infty])$) for the extension of $\mathbb{Q}$ obtained by appending to $\mathbb{Q}$ the coordinates of all $p^a$ (resp.~all $p$-power) torsion points of $\mathcal{A}$. With our fixed $c$, we write $F = K_c$.
\end{notation}
With these notations we have $\Gamma \cong \Gal(K_\infty/F)$.

\begin{definition}[Iwaswa algebra~{\cite[\S\S2.4 and 3.1]{ArdakovWadsleyVerma}}]
We shall refer to the group rings, $\Omega := \Zp[[\Gamma]]$ and $\Lambda := \Omega[1/p]$ as the integral, and rational Iwasawa algebras.
\end{definition}

\begin{remark}[Noncommutative power-series ring presentation~{\cite[\S4, discussion before Lemma~4.1]{SchneiderTeitelbaumDistributions}}]\label{rmk:iwasawa_generators}
The Iwasawa algebra $\Omega = \Zp[[\Gamma]]$ has a description as a noncommutative power series ring, given as follows. Since $H_c$ is uniform, it is finitely topologically generated, and we write $\gamma_1,\dots,\gamma_{d_H}$ for its topological generators. Similarly for $Z$, write $\gamma_0\in Z$ for a topological generator. For every $1\le i\le d_H$, fix $u_i := \gamma_i - 1$, and write $T := \gamma_0 - 1$. Then every element in $\Omega$ may be written uniquely as
\[
\sum_{\underline{\alpha}\in \mathbb{Z}_{\ge 0}^{d_H+1}}b_{\underline{\alpha}}\cdot u_1^{\alpha_1}\cdot \dots \cdot u_{d_H}^{\alpha_{d_H}}\cdot T^{\alpha_0},\qquad b_{\underline{\alpha}}\in \Zp,
\]
giving $\Omega$ the structure of a noncommutative power-series ring, with $T$ a central coordinate. In particular, $\Omega\cong \Zp[[H_c]][[T]]$. $\Lambda$ may be presented similarly, with the exception that the $b_{\underline{\alpha}}$ are now in $\Qp$, and have a common, finite $p$-power, denominator.
\end{remark}

\begin{proposition}[{\cite[Theorem~A.1(a),(b)]{LimFineSelmer}}]\label{prop:iwasawa_noetherian_domain}
The Iwasawa algebras $\Omega$ and $\Lambda$ are noetherian domains.
\end{proposition}
\begin{proof}
The integral Iwasawa algebra is a (left and right) noetherian domain by Lim~\cite[Theorem~A.1(a),(b)]{LimFineSelmer}. These properties persist under central localization, implying the same holds for $\Lambda$.
\end{proof}

A (left and right) noetherian domain $R$ admits a skew field of fractions, denoted $Q(R)$; see~\cite[\S4, opening paragraph]{LimFineSelmer}. In particular, by Proposition~\ref{prop:iwasawa_noetherian_domain}, 
\begin{definition}
Let $R$ denote a noetherian domain, and let $M$ denote an $R$-module. We define the rank of $M$ as follows:
\[
\operatorname{rank}_R M := \dim_{Q(R)} M\otimes_R Q(R).
\]
An $R$-module $M$ is a torsion $R$-module if it has rank $0$.
\end{definition}

\begin{notation}
For an extension $E$ of $\mathbb{Q}$, write $\Sigma_E$ for the lifts of all places in $\Sigma$ to $E$, and let $G_{E,\Sigma_E}$ denote the Galois group of the maximal extension of $E$ unramified outside of $\Sigma_E$.
\end{notation}

One of the key objects of interest in this paper will be the Iwasawa cohomology group $H^2_{\mathrm{Iw}}(K_\infty/F;\Zp(1))$. Let us recall first the definition of Iwasawa cohomology.

\begin{definition}[Iwasawa cohomology~{\cite[\S2]{LimFineSelmer}}]\label{def:iwasawa_coh}
Let $M$ denote a finitely generated $\Zp$-module, endowed with the $p$-adic topology, equipped with a continuous $\mathbb{Z}_p[G_{F,\Sigma_F}]$-action. For $i\ge 0$, we define the Iwasawa cohomology group,
\[
H^i_{\mathrm{Iw}}(K_\infty/F; M) := \varprojlim_{K_\infty\supseteq E\supseteq F}H^i(G_{E,\Sigma_E}; M),
\]
where the transition maps are given by corestriction, and the cohomology on the right is continuous.
\end{definition}

\begin{remark}[Corestriction in group cohomology]\label{rmk:corestriction}
Let $N\trianglelefteq G$ denote a finite index normal subgroup. Write $i: BN\to BG$ for the essential geometric morphism of $\infty$-topoi of functors $\Fun(BN,\mathcal{D}(\mathbb{Z}))$ and $\Fun(BG,\mathcal{D}(\mathbb{Z}))$, with $\mathcal{D}(\mathbb{Z})$ the stable $\infty$-category of complexes of abelian groups. From $i$, we obtain three natural functors on underlying derived categories of sheaves, $i^*, i_*$, and $i_!$, referred to as restriction, coinduction, and induction, respectively. These functors satisfy the adjunction relations,
\[
i_!\dashv i^*\dashv i_*.
\]
Let $M$ denote a $\mathbb{Z}[G]$-module, regarded as a complex concentrated in degree $0$ in $\Fun(BG, D(\mathbb{Z}))$. The corestriction maps on cohomology are the ones induced by the following sequence of morphisms in the derived category:
\begin{equation}\label{eq:corestriction}
R\Gamma(N; i^*M) \cong R\Gamma(G; i_*i^*M) \cong R\Gamma(G; i_!i^*M) \to R\Gamma(G; M),
\end{equation}
where the last map is the counit of the induction-restriction adjunction, and the equivalence before it is due to the fact that the canonical comparison $i_! \to i_*$ is an equivalence: fibrewise, this is the identification of finite direct sums with finite direct products (see~\cite[Lemma I.1.9]{NikolausScholze2018}). The complex $R\Gamma(N; i^*M)$ is equipped with a homotopy coherent residual action of $G/N$, endowing the cohomology groups $H^i(N; i^*M)$ with a $\mathbb{Z}[G/N]$-module structure.

In the Iwasawa theoretic setting, we work with continuous cochains with $\Zp$-coefficients, as well as the continuous versions of $i^*,i_*,i_!$, however, not much changes beyond it; see~\cite[Lemma~8.1.5 and \S\S~8.2.1-8.2.2]{NekovarSelmerComplexes}. The Iwasawa cohomology groups obtained in this way form a cofiltered limit, equipped with a $\varprojlim_{K_\infty\supseteq E\supseteq F}\Zp[G_{F,\Sigma_F}/G_{E,\Sigma_E}] \cong \Omega$ action.
\end{remark}

\subsection{The Iwasawa module $H^2_{\mathrm{Iw}}(K_\infty/F;\Zp(1))$}\label{ssec:key_reduction}
We are now ready to define the Iwasawa module $\mathcal{X}$, which will play a key role in the proof of our main theorem, as explained in~\S\ref{sssec:h2-bound-proof-idea}.
\begin{definition}\label{def:x}
Denote by $\mathcal{X}$ and $X$ the integral and rational Iwasawa modules 
\[
\mathcal{X} := H^2_{\mathrm{Iw}}(K_\infty/F;\Zp(1)),\qquad X := \mathcal{X}[1/p].
\]
\end{definition}
The rest of this subsection is dedicated to the study of $\mathcal{X}$ and $X$'s properties. Namely, we prove that they are finitely generated and torsion Iwasawa modules.

\begin{proposition}\label{prop:fg_X}
The Iwasawa modules $\mathcal{X}$ and $X$ are finitely generated.
\end{proposition}

\begin{proof}
Since finite generation persists from $\mathcal{X}$ to $X$, it suffices to prove the property for $\mathcal{X}$.
The Iwasawa version of Shaprio's lemma, i.e.~of the equivalence $R\Gamma(N; i^*M) \cong R\Gamma(G; i_!i^*M)$ from~\eqref{eq:corestriction}, gives rise to the isomorphism
\[
\mathcal{X} := \varprojlim_{K_\infty\supseteq E\supseteq F}H^2(G_{E,\Sigma_E}; \Zp(1)) \cong H^2(G_{F,\Sigma_F}; \Zp(1)\widehat \otimes_{\Zp} \Omega);
\]
see~\cite[Lemma 5.3.1]{LimPoitouTate}. Since $\Zp(1)$ is a finitely generated $\Zp$-module, the canonical map $\Zp(1) \otimes_{\Zp} \Omega\to \Zp(1)\widehat \otimes_{\Zp} \Omega$ is an equivalence (see~\cite[Lemma 5.1.5]{LimPoitouTate}), hence we may omit the completed tensors from the notation. Furthermore, since all tensor products below are taken over $\Zp$, we shall simply write $\otimes$ for brevity.

Here, the module $\Zp(1)\otimes \Omega$ has commuting left Galois-action, and right $\Omega$-action, given by
\[
g(t\otimes x)\mapsto gt\otimes \overline{g}x,\qquad (t\otimes x)\cdot y = t\otimes (x\cdot y),
\]
for $x,y\in \Omega$, $t\in \Zp(1)$ and $g\in G_{F,\Sigma_F}$, where $\overline{g}$ denotes $g$'s image in $\Gamma$. In the notations of Remark~\ref{rmk:iwasawa_generators}, write $\mathfrak{m} = (p, u_0,\dots,u_{d_H})$ for the generators of the maximal ideal in $\Omega$, regarded as right-generators of the associated left ideal. For brevity, let us write here $u_{d_H+1} := p$. Write 
\[
d: \bigoplus_{i=0}^{d_H+1}\Omega\to \Omega \cdot \mathfrak{m},\qquad \underline{x} = (x_0,\dots, x_{d_H+1})\mapsto \sum_{i=0}^{d_H+1} x_i\cdot u_i.
\]
Consider the short exact sequences of $\Zp[G_{F,\Sigma_F}]$-modules:
\[
0\to \mathrm{Ker}(d)\to \bigoplus_{i=0}^{d_H+1}\Omega\xrightarrow{d} \mathfrak{m} \to 0,
\]
and
\[
0\to \mathfrak{m} \to \Omega\to \Fp\to 0.
\]

Tensoring with $\Zp(1)$ is exact, and gives rise to long exact sequences in cohomology
\[
\dots\to \bigoplus_{i=0}^{d_H+1}H^2(G_{F,\Sigma_F}, \Zp(1)\otimes \Omega)\xrightarrow{d_*} H^2(G_{F,\Sigma_F}, \Zp(1)\otimes \mathfrak{m})\to 0,
\]
and
\[
\dots \to H^2(G_{F,\Sigma_F}, \Zp(1)\otimes \mathfrak{m})\to H^2(G_{F,\Sigma_F}, \Zp(1)\otimes \Omega) \to H^2(G_{F,\Sigma_F}, \mu_p) \to 0.
\]
Here the long exact sequences truncate at the $H^3$ stage, since the $p$-cohomological dimension of $G_{F,\Sigma_F}$ is at most $2$ for $p$ odd. The long exact sequences give rise to the exact sequence
\begin{equation}\label{eq:les_fg}
\bigoplus_{i=0}^{d_H+1}H^2(G_{F,\Sigma_F}, \Zp(1)\otimes \Omega)\to H^2(G_{F,\Sigma_F}, \Zp(1)\otimes \Omega) \to H^2(G_{F,\Sigma_F}, \mu_p) \to 0.
\end{equation}
It follows that
\[
\operatorname{rank}_{\Omega} H^2(G_{F,\Sigma_F}; \Zp(1)\otimes \Omega) \le \dim_{\Fp}\dfrac{H^2(G_{F,\Sigma_F}; \Zp(1)\otimes \Omega)}{\mathfrak{m}\cdot H^2(G_{F,\Sigma_F}; \Zp(1)\otimes \Omega)}
= \dim_{\Fp} H^2(G_{F,\Sigma_F}; \mu_p).
\]
where we used compact Nakayama's lemma for the first inequality, and~\eqref{eq:les_fg} for the following equality. The finiteness of $\dim_{\Fp} H^2(G_{F,\Sigma_F}; \mu_p)$ now follows from the general fact (see~\cite[I, Corollary 4.15]{MilneADT}) that if $M$ is a finite discrete continuous $G_{F,\Sigma_F}$-module, $\Sigma_F$ contains all archimedean places, and all places above the rational primes dividing $|M|$, then
\[
\#H^2(G_{F,\Sigma_F}; M) < \infty.
\]
\end{proof}

Next, we shall prove that $\mathcal{X}$ and $X$ are torsion Iwasawa modules using a general criterion by Lim; see Theorem~\ref{thm:lim_criterion} below. To state Lim's criterion, we need to define $\Sigma$-admissible extensions. Recall that for a number field $L/\mathbb{Q}$, its $p$-adic cyclotomic extension, denoted $L^\cyc$, is the unique sub-extension of $L(\mu_{p^{\infty}})/L$ whose Galois group is isomorphic to $\Zp$.
\begin{definition}[See~{\cite[discussion before Lemma~5.2]{LimFineSelmer}}]
Let $L/\mathbb{Q}$ denote a number field. An extension $L_\infty$ of $L$ is said to be $\Sigma_L$-admissible if:
\begin{enumerate}
  \item $\Gamma_L := \Gal(L_\infty/L)$ is a uniform pro-$p$ group,
  \item $L_\infty \supseteq L^{\cyc}$, and
  \item $L_\infty$ is contained in the maximal extension of $L$ unramified outside $\Sigma_L$.
\end{enumerate}
\end{definition}

We are ready to state Lim's criterion.
\begin{theorem}[See~{\cite[Lemma 7.1; (a) iff (c)]{LimFineSelmer}}]\label{thm:lim_criterion}
Notations as above. Let $L_\infty$ be a $\Sigma_L$-admissible $p$-adic Lie extension of $L$, and let $T$ denote a free $\Zp$-module of finite rank, equipped with a continuous $G_{L,\Sigma_L}$-action. Then 
\[
H^2_{\mathrm{Iw}}(L_\infty/L; T)\text{ is a torsion } \Zp[[\Gamma_L]]\text{-module}\iff H^2(G_{L_{\infty},\Sigma_{L_{\infty}}}; T \otimes_\Zp \Qp/\Zp) = 0.
\]
\end{theorem}

\begin{proposition}\label{prop:torsion_X}
The Iwasawa modules $\mathcal{X}$ and $X$ are torsion.
\end{proposition}

\begin{proof}
Since $\Omega$ is $p$-torsion free (as a domain), the localization map $\Omega \to \Lambda$ is an inclusion. Hence, if $\mathcal{X}$ is torsion, so is $X$, and it suffices to establish the torsion property for $\mathcal{X}$.

Since $K_\infty/F$ is $\Sigma_F$-admissible, Theorem~\ref{thm:lim_criterion} applies, and it suffices to prove that $H^2(G_{K_{\infty,\Sigma_{K_{\infty}}}}; \mu_{p^{\infty}}) = 0$. Observe that
\[
H^2(G_{K_{\infty,\Sigma_{K_{\infty}}}}; \mu_{p^{\infty}}) = \varinjlim_{m\ge 1}\varinjlim_{K_\infty\supseteq E\supseteq F}H^2(G_{E,\Sigma_E}; \mu_{p^m}),
\]
where the direct limits are taken with respect to: (i) the filtered system of finite intermediate extensions of $F$, with transition maps given by restriction, and (ii) the natural maps induced functorially by the inclusions $\mu_{p^{m}}\hookrightarrow \mu_{p^{m+1}}$.

Fix $m\ge 1$ and a number field $K_\infty\supseteq E\supseteq F$. Write $\Sigma_E^{\mathrm{fin}}$ for the finite places of $E$ within $\Sigma_E$, $U := \Spec\,\mathcal{O}_{E,\Sigma_E^{\mathrm{fin}}}$ for the affine scheme associated to the spectrum of the ring $\mathcal{O}_{E,\Sigma_E^{\mathrm{fin}}}$ of $\Sigma_E^{\mathrm{fin}}$-integers, and $\overline\eta$ for $U$'s geometric point. Since $p$ is invertible in $\mathcal{O}_{E,\Sigma_E^{\mathrm{fin}}}$, by Milne~\cite[II, Proposition~2.9]{MilneADT}, for any constructible sheaf $\mathcal{F}$ of locally constant $p$-primary torsion abelian groups on $U$, we have
\[
H^r_{\et}(U; \mathcal F) \cong H^r(G_{E,\Sigma_E}; \mathcal F_{\overline\eta}),\qquad \text{for every }r.
\]
We apply Milne's proposition with $r=2$, and with $\mathcal F$ the \'etale local system $\mu_{p^m,U}$ on $U$. It follows that
\[
H^2_{\et}(U; \mu_{p^m,U}) = H^2(G_{E,\Sigma_E}; \mu_{p^m}).
\]
The Kummer short exact sequence for $U$ reads as:
\[
1\to \mu_{p^m,U}\to \mathbb{G}_{m,U}\xrightarrow{\cdot p^m} \mathbb{G}_{m,U}\to 1,
\]
and its long exact sequence in cohomology yields the short exact sequence
\[
0\to \mathrm{Cl}(\mathcal{O}_{E,\Sigma_E^{\mathrm{fin}}})/p^m \to H^2_{\et}(U; \mu_{p^m,U})\to \mathrm{Br}(\mathcal{O}_{E,\Sigma_E^{\mathrm{fin}}})[p^m]\to 0,
\]
where $\mathrm{Cl}(\mathcal{O}_{E,\Sigma_E^{\mathrm{fin}}})$ and $\mathrm{Br}(\mathcal{O}_{E,\Sigma_E^{\mathrm{fin}}})$ denote the class group and the Brauer group of $\mathcal{O}_{E,\Sigma_E^{\mathrm{fin}}}$. Since filtered colimits are exact, taking the colimit on $m$ of the above sequence, we get the following short exact sequence
\[
0\to \varinjlim_{m}\mathrm{Cl}(\mathcal{O}_{E,\Sigma_E^{\mathrm{fin}}})/p^m \to H^2_{\et}(U; \mu_{p^\infty,U})\to \mathrm{Br}(\mathcal{O}_{E,\Sigma_E^{\mathrm{fin}}})[p^\infty]\to 0,
\]
where the transition maps on the direct limit $\varinjlim_{m}\mathrm{Cl}(\mathcal{O}_{E,\Sigma_E^{\mathrm{fin}}})/p^m$ are given by multiplication by $p$. The groups $\mathrm{Cl}(\mathcal{O}_{E,\Sigma_E^{\mathrm{fin}}})$ are finite, hence every element of $\mathrm{Cl}(\mathcal{O}_{E,\Sigma_E^{\mathrm{fin}}})/p^m$ is $p$-power torsion, and therefore
\[
\varinjlim_{m}\mathrm{Cl}(\mathcal{O}_{E,\Sigma_E^{\mathrm{fin}}})/p^m = 0.
\]
It follows that,
\[
H^2_{\et}(U; \mu_{p^\infty,U}) \cong \mathrm{Br}(\mathcal{O}_{E,\Sigma_E^{\mathrm{fin}}})[p^\infty].
\]
Our goal is therefore to prove that 
\begin{equation}\label{eq:brauer_vanishing}
\varinjlim_{E}\mathrm{Br}(\mathcal{O}_{E,\Sigma_E^{\mathrm{fin}}})[p^\infty] = 0,
\end{equation}
where the direct limit is taken over the natural restriction maps on Brauer groups, taken along the finite extensions $K_\infty\supseteq E\supseteq F$. According to~\cite[II, Proposition~2.1 and Remark~2.2(a)]{MilneADT}, there is an embedding
\[
0\to \mathrm{Br}(\mathcal{O}_{E,\Sigma_E^{\mathrm{fin}}}) \to \bigoplus_{v\in \Sigma_E}\mathrm{Br}(\mathcal{O}_{E_v}).
\]
Now, observe that the archimedean places do not contribute any $p$-torsion to the right-hand side. Indeed, if $v$ is archimedean, we have $\mathrm{Br}(\mathcal{O}_{E_v})[p^\infty] = 0$ since $\mathrm{Br}(\mathbb{C}) = 0$ and $\mathrm{Br}(\mathbb{R}) = \mathbb{Z}/2\mathbb{Z}$, while $p$ is odd. Therefore, on $p$-torsion, we get the inclusion
\[
\mathrm{Br}(\mathcal{O}_{E,\Sigma_E^{\mathrm{fin}}})[p^\infty]\subseteq \bigoplus_{v\in \Sigma_E^{\mathrm{fin}}}\mathrm{Br}(\mathcal{O}_{E_v})[p^\infty].
\]
By~\cite[III, Theorem~2.1]{MilneCFT}, for every $v\in \Sigma_E^{\mathrm{fin}}$, there are isomorphisms,
\[
\mathrm{inv}_{E_v}: \mathrm{Br}(\mathcal{O}_{E_v})[p^\infty] \cong \Qp/\Zp,
\]
and for every local field extension $L_w\supseteq E_v$, there is a commutative diagram
\[
\begin{tikzcd}
\mathrm{Br}(\mathcal{O}_{E_v})[p^\infty]\arrow[r, "\mathrm{inv}_{E_v}"]\arrow[d, "\mathrm{res}"]& \Qp/\Zp\arrow[d, "{\cdot [L_w:E_v]}"]\\
\mathrm{Br}(\mathcal{O}_{L_w})[p^\infty]\arrow[r, "\mathrm{inv}_{L_w}"]& \Qp/\Zp.
\end{tikzcd}
\]
Here the horizontal arrows are isomorphisms, the left vertical arrow is the local restriction map, while the right vertical arrow is multiplication by the extension degree. Since every class in $\mathrm{Br}(\mathcal{O}_{E,\Sigma_E^{\mathrm{fin}}})[p^\infty]$ is $p$-power torsion, we can establish the vanishing of the direct limit~\eqref{eq:brauer_vanishing} if we can construct, for every finite extension $K_\infty\supset E\supseteq F$ and every place $v\in \Sigma_E^{\mathrm{fin}}$, an intermediate finite extension $K_\infty \supset L\supseteq E$ and a place $w$ of $L$ lying above $v$, such that
\[
p^m|[L_w:E_v].
\]
We shall construct $L$ from $E$ by appending sufficiently high $p$-power roots of unity, which is justified since $E^\cyc = \bigcup_{n\ge 0} E(\mu_{p^n})\subseteq K_\infty$. We study a by case analysis, with respect to whether $v|p$ or not.

\textbf{The case where $v|p$.} Observe that,
\[
[E_v(\mu_{p^n}): E_v]\cdot [E_v : \Qp] = [E_v(\mu_{p^n}): \Qp] = [E_v(\mu_{p^n}): \Qp(\mu_{p^n})]\cdot [\Qp(\mu_{p^n}): \Qp].
\]
Therefore, $[\Qp(\mu_{p^n}): \Qp]$ divides $[E_v(\mu_{p^n}): E_v]\cdot [E_v : \Qp]$ and since $[E_v : \Qp]$ is a constant, the $p$-adic valuations of the sequence of extension degrees $([E_v(\mu_{p^n}): E_v])_{n\ge 1}$ are unbounded if those of the sequence $([\Qp(\mu_{p^n}): \Qp])_{n\ge 1}$ are unbounded. As for the latter, by~\cite[Proposition~6.2(a),(c)]{MilneANT}, $\Qp(\mu_{p^n})$ is a totally ramified extension of $\Qp$ of degree $\phi(p^n) = p^{n-1}(p-1)$, and, in particular, is unbounded.

\textbf{The case where $v\nmid p$.} Denote $D_v^\cyc := \Gal(E_v(\mu_{p^\infty})/E_v)$, and let $\chi_\cyc: D_v^\cyc \to \Zp^\times$ denote the $p$-adic cyclotomic character. Since $E_v\supseteq \mu_p$, $\chi_\cyc$ factors through $1 + p\Zp\subseteq \Zp^\times$. We argue that $\chi_\cyc$ is injective on $D_v^\cyc$. Indeed, if $\sigma \in D_v^\cyc$, and $\zeta\in \mu_{p^n}$ is a primitive $p^n$'th root of unity, then $\sigma(\zeta) = \zeta^{\chi_\cyc(\sigma)}$. Hence, if $\chi_\cyc(\sigma) = 1$, $\sigma$ acts trivially on all $p$-power roots of unity, hence $\sigma = 1$ in $D_v^\cyc$. It follows that $D_v^\cyc \cong \chi_\cyc(D_v^\cyc)$ is a continuous bijection from $D_v^\cyc$ onto its image in $1 + p\Zp$. Since the former is compact and the latter is Hausdorff, the bijection $D_v^\cyc \cong \chi_\cyc(D_v^\cyc)$ is a homeomorphism. In particular,
\[
D_v^\cyc = \overline{\langle \Frob_v\rangle} \cong \chi_\cyc(\overline{\langle \Frob_v\rangle}) = \overline{\chi_\cyc(\langle \Frob_v\rangle)}.
\]
Recall that
\[
\chi_\cyc(\Frob_v) = \#k(v) \in 1 + p\Zp.
\]
Since the group $(1 + p\Zp,\times)$ is torsion free, and $\#k(v) > 1$, $\chi_\cyc(\Frob_v)$ generates an infinite cyclic subgroup of $(1 + p\Zp,\times)$. Since $p$ is odd, the $p$-adic logarithm is a homeomorphism,
\[
\log_p: (1 + p\Zp,\times)\to (p\Zp,+),
\]
and under $\log_p\circ \chi_\cyc$, $D_v^\cyc$ becomes isomorphic with an infinite closed subgroup of $p\Zp$, hence of the form $p^c\Zp\subseteq p\Zp$ for some $c\ge 1$. In particular, $D_v^\cyc$ is infinite pro-$p$, and has arbitrarily large $p$-power degree subextensions.
\end{proof}

\subsection{From Galois cohomology groups to kernels of Iwasawa operators}\label{ssec:red2}
This subsection is dedicated to the reduction step in the proof of Theorem~\ref{thm:h2-bound}, and to the proof of the trivial bound $\dim_\Qp H^2(G_{\mathbb{Q},\Sigma}, V_{\lambda,s})\ll d_\lambda$, as explained in~\S\ref{sssec:h2-bound-proof-idea}. The reduction itself shall be carried out in three steps.

We begin with the first reduction step.
We are now ready to carry out the first reduction step in our proof of .
\begin{theorem}\label{thm:galois_iwasawa}
Let $W$ be a finite dimensional $\mathbb{Q}_p$-vector space equipped with a $G_{F,\Sigma_F}$-action, which factors through $\Gamma$. Then there is an isomorphism,
\[
H^2(G_{F,\Sigma_F}; W)^\vee \cong \Hom_\Gamma(X,W^\vee(1)).
\]
\end{theorem}
\begin{proof}
Let $\mathcal{W}\subset W$ denote a Galois stable $\Zp$-lattice. Write $\mathcal{Y}$ for the integral Iwasawa module,
\[
\mathcal{Y} := H^2_{\mathrm{Iw}}(K_\infty/F; \mathcal{W}).
\]
We argue that for every $m\ge 1$, one has
\begin{equation}\label{eq:y-identity}
\mathcal{Y}/p^m \cong H^2_{\mathrm{Iw}}(K_\infty/F;\mathcal{W}/p^m).  
\end{equation}
Indeed, the Iwasawa version of Shapiro's lemma gives an isomorphism $\mathcal{Y} \cong H^2(G_{F,\Sigma_F}; \mathcal{W}\otimes\Omega)$. Since the coefficient $\mathcal{W}\otimes\Omega$ is $p$-adic and $p$ is odd, we have an isomorphism,
\[
\mathcal{Y}/p^m \cong H^2(G_{F,\Sigma_F}; \mathcal{W}\otimes\Omega)/p^m \cong H^2(G_{F,\Sigma_F}; (\mathcal{W}/p^m)\otimes\Omega).
\]
Applying Shapiro's lemma again but in the other direction, we obtain~\eqref{eq:y-identity}.

Next, we argue that
\begin{equation}\label{eq:xy-identity}
H^2_{\mathrm{Iw}}(K_\infty/F;\mathcal{W}/p^m) \cong \mathcal{X}/p^m\otimes \mathcal{W}(-1)/p^m.
\end{equation}
For this, write $\mathcal{W}/p^m = \mu_{p^m}\otimes \mathcal{W}(-1)/p^m$. Since $\mathcal{W}(-1)/p^m$ is a finite Galois module whose Galois action factors through $\Gamma$\footnote{This is true for $\mathcal{W}$ by our assumption. As for $\Zp(-1)$, this follows from the fact that $K_\infty\supseteq F(\mu_{p^\infty})$.}, there is a finite extension $K_\infty\supseteq F_m\supseteq F$, with the property that $G_{F_m,\Sigma_{F_m}}$ acts trivially on $\mathcal{W}(-1)/p^m$. Since the directed system of subextensions $K_\infty\supseteq E\supseteq F_m$ is cofinal in $K_\infty\supseteq E\supseteq F$, we have the equivalences
\[
H^2_{\mathrm{Iw}}(K_\infty/F;\mu_{p^m}\otimes \mathcal{W}(-1)/p^m) := \varprojlim_{K_\infty\supseteq E\supseteq F}H^2(G_{E,\Sigma_E}; \mu_{p^m}\otimes \mathcal{W}(-1)/p^m) 
\]
\[
\cong \varprojlim_{K_\infty\supseteq E\supseteq F_m}H^2(G_{E,\Sigma_E}; \mu_{p^m}\otimes \mathcal{W}(-1)/p^m) \cong \varprojlim_{K_\infty\supseteq E\supseteq F_m}\left(H^2(G_{E,\Sigma_E}; \mu_{p^m})\otimes \mathcal{W}(-1)/p^m\right).
\]
Since $\mathcal{W}(-1)/p^m$ is a free $\Zp/p^m$-module, tensoring with it commutes with inverse limits. In particular,
\[
\varprojlim_{K_\infty\supseteq E\supseteq F_m}\left(H^2(G_{E,\Sigma_E}; \mu_{p^m})\otimes \mathcal{W}(-1)/p^m\right)\cong \left(\varprojlim_{K_\infty\supseteq E\supseteq F_m} H^2(G_{E,\Sigma_E}; \mu_{p^m})\right)\otimes \mathcal{W}(-1)/p^m.
\]
Identifying $\varprojlim_{K_\infty\supseteq E\supseteq F_m} H^2(G_{E,\Sigma_E}; \mu_{p^m})$ with $\mathcal{X}/p^m$, we obtain~\eqref{eq:xy-identity}. Combining~\eqref{eq:y-identity} and~\eqref{eq:xy-identity} together, we obtain the identity
\[
\mathcal{Y}/p^m \cong \mathcal{X}/p^m\otimes \mathcal{W}(-1)/p^m \cong \mathcal{X}/p^m\otimes \mathcal{W}(-1).
\]
Taking inverse limits with respect to $m\ge 1$, we obtain the equivalence 
\begin{equation}\label{eq:y-cong}
\mathcal{Y}\cong \mathcal{X} \otimes \mathcal{W}(-1).
\end{equation}
Indeed, since tensoring with $\mathcal{W}(-1)$ commutes with inverse limits, it suffices to prove that $\varprojlim_m \mathcal{X}/p^m = \mathcal{X}$ and similarly for $\mathcal{Y}$. We shall argue for $\mathcal{X}$ in detail, the analogous result for $\mathcal{Y}$ is then proven similarly. Since 
\[
\mathcal{X}/p^m\cong H^2_{\mathrm{Iw}}(K_\infty/F; \mu_{p^m})\cong \varprojlim_{E}H^2(G_{E,\Sigma_E}; \mu_{p^m}),
\]
and since limits commute with limits, it suffices to show that $\varprojlim_m H^2(G_{E,\Sigma_E}; \mu_{p^m}) \cong H^2(G_{E,\Sigma_E}; \Zp(1))$. As for the latter, the Grothendieck spectral sequence for the pair of functors $\varprojlim_m$ and $(-)^{G_{E,\Sigma_E}}$ gives rise to the short exact sequence (see~\cite[Eq.~(2.1) and Theorem~2.2, p.~215]{JannsenContinuousEtale}):
\[
0\to \varprojlim^1_{m\ge 1}H^1(G_{E,\Sigma_E}; \mu_{p^m})\to H^2(G_{E,\Sigma_E}; \Zp(1))\to \varprojlim_{m\ge 1}H^2(G_{E,\Sigma_E}; \mu_{p^m})\to 0.
\]
Since the groups $H^1(G_{E,\Sigma_E}; \mu_{p^m})$ are finite, the sequence $(H^1(G_{E,\Sigma_E}; \mu_{p^m}))_{m\ge 1}$ is Mittag--Leffler, and therefore 
\[
\varprojlim^1_{m\ge 1}H^1(G_{E,\Sigma_E}; \mu_{p^m}) = 0.
\]
In particular, $H^2(G_{E,\Sigma_E}; \Zp(1))\cong \varprojlim_{m\ge 1}H^2(G_{E,\Sigma_E}; \mu_{p^m})$ holds, establishing 
\[
\varprojlim_m \mathcal{X}/p^m = \mathcal{X},
\]
and, as explained above, also establishes~\eqref{eq:y-cong}.

Taking continuous $\Gamma$-coinvariants,~\eqref{eq:y-cong} becomes
\begin{equation}\label{eq:coinvs}
\mathcal{Y}_\Gamma^{\mathrm{cont}} \cong \left(\mathcal{X} \otimes \mathcal{W}(-1)\right)_\Gamma^{\mathrm{cont}}.  
\end{equation}
We wish to argue that, in this particular setting, continuous $\Gamma$-coinvariants and ordinary $\Gamma$-coinvariants coincide.

For this, let $(\gamma_i)_{i=0}^{d_H}$ denote the topological generators of $\Gamma$ from Remark~\ref{rmk:iwasawa_generators}. The $\gamma_i$ act on $M := \mathcal{X} \otimes \mathcal{W}(-1)$ through the natural inclusion $\Gamma\subseteq \Zp[[\Gamma]] =: \Omega$. Denote by $M^{\mathrm{alg}}_\Gamma$ the ordinary $\Gamma$-coinvariatns of $M$, and by $D\subseteq M$ the image of the continuous map $M^{d_H+1}\to M$, sending $(m_0,\dots,m_{d_H})$ to $\sum_{i=0}^{d_H}m_i\cdot (\gamma_i-1)$. Clearly,
\[
D\subseteq \mathrm{Ker}\left(M\to M_\Gamma^{\mathrm{alg}}\right).
\]
We seek to show the inverse inclusion. Since $M$ is compact Hausdorff (by Proposition~\ref{prop:fg_X}, and the fact that $\Omega$ is noetherian), so is the image of $M^{d_H+1}\to M$, and thus $D\subseteq M$ is closed. Observe that $D$ contains $m(\gamma-1)$ for every $\gamma\in \Gamma$ and $m\in M$, where $\gamma$ is a word in the $\gamma_i^{\pm 1}$. Indeed, this follows by induction on the word length using the identities
\[
m(ab-1) = (ma)(b-1) + m(a-1),\qquad m(a^{-1}-1) = (-ma^{-1})(a-1).
\]
Now, fix $m\in M$ and consider the continuous map $f_m: \Gamma\to M$, sending $\gamma$ to $m(\gamma-1)$. Since $D$ is closed, and the map sends a dense subset of $\Gamma$ to $D$, $D$ contains the image of $f_m$. Since this is true for every $m\in M$, we obtain the inverse inclusion, hence
\[
D = \mathrm{Ker}\left(M\to M_\Gamma^{\mathrm{alg}}\right).
\]
Since continuous coinvariants are then given by $M/\overline{D}$ and $D$ is closed, it follows that
\[
M_\Gamma^{\mathrm{cont}} = M_\Gamma^{\mathrm{alg}}.
\]
Since~\eqref{eq:y-cong} is an equivalence of $\Omega$, hence also of $\Gamma$-modules, the same ordinary-continuous coinvariants coincidece also holds for the left hand side of~\eqref{eq:coinvs}, i.e.~for $\mathcal{Y}$. In particular, since there is no difference between ordinary and continuous coinvariants in this setting, we may denote both by $(\bullet)_\Gamma$, i.e.~omitting the superscripts ``$\mathrm{cont}$'' or ``$\mathrm{alg}$''.

We proceed to analyze the left hand side of~\eqref{eq:coinvs}. Write $\mathcal I\subseteq \Zp[\Gamma]$ for the ideal in $\Gamma$, which is left-generated by the $(\gamma_i - 1)_{i=0}^{d_H}$. The previous argument shows that taking the continuous $\Gamma$ coinvariants of $\mathcal{Y}$ coincides with $\mathcal{Y}/\mathcal{Y}\mathcal{I}$. Applying Shapiro's lemma to $\mathcal{Y} := H^2_{\mathrm{Iw}}(K_\infty/F;\mathcal{W})$, we obtain 
\[
\mathcal{Y} \cong H^2(G_{F,\Sigma_F}; \mathcal{W}\otimes\Omega),
\]
and therefore,
\[
\mathcal{Y}_\Gamma \cong H^2(G_{F,\Sigma_F}; \mathcal{W}\otimes\Omega)/\mathcal{I}.
\]
Consider the right exact sequence of $G_{F,\Sigma_F}$-modules,
\[
\bigoplus_{i=0}^{d_H}\mathcal{W}\otimes\Omega\xrightarrow{\underline{x}\mapsto \sum x_i\cdot (\gamma_i-1)}\mathcal{W}\otimes\Omega\to \mathcal{W}\to 0.
\]
It induces the right exact sequence in cohomology:
\[
\bigoplus_{i=0}^{d_H}H^2(G_{F,\Sigma_F}; \mathcal{W}\otimes\Omega)\xrightarrow{\underline{c}\mapsto \sum c_i\cdot (\gamma_i-1)}H^2(G_{F,\Sigma_F}; \mathcal{W}\otimes\Omega)\to H^2(G_{F,\Sigma_F}; \mathcal{W})\to 0,
\]
giving rise to the isomorphism
\[
H^2(G_{F,\Sigma_F}; \mathcal{W}\otimes\Omega)/\mathcal{I}\cong H^2(G_{F,\Sigma_F}; \mathcal{W}).
\]
Plugging back to~\eqref{eq:coinvs}, we get
\begin{equation}\label{eq:dual-of-final}
H^2(G_{F,\Sigma_F}; \mathcal{W})\cong (\mathcal{X}\otimes \mathcal{W}(-1))_\Gamma.
\end{equation}
Taking linear duals, i.e.~$N^\vee:= \Hom(N,\Zp)$, we obtain
\[
H^2(G_{F,\Sigma_F}; \mathcal{W})^\vee \cong \Hom((\mathcal{X}\otimes \mathcal{W}(-1))_\Gamma,\Zp)\cong \Hom_\Gamma(\mathcal{X}\otimes \mathcal{W}(-1),\Zp).
\]
The usual $\otimes$-$\Hom$-adjunction gives
\[
\Hom_\Gamma(\mathcal{X}\otimes \mathcal{W}(-1),\Zp) \cong \Hom_\Gamma(\mathcal{X},\mathcal{W}^{\vee}(1)),
\]
and combining the pieces together, we obtain the isomorphism
\[
H^2(G_{F,\Sigma_F}; \mathcal{W})^\vee \cong \Hom_\Gamma(\mathcal{X},\mathcal{W}^{\vee}(1)).
\]
The theorem now follows by inverting $p$.
\end{proof}

\begin{corollary}\label{cor:galois-iwasawa}
Let $W$ be a $G_{\mathbb{Q},\Sigma}$-module, for which the Galois action, when restricted to $G_{F,\Sigma_F}$ factors through $\Gamma$. Then
\[
\dim_\Qp H^2(G_{\mathbb{Q},\Sigma}; W)\le \dim_\Qp H^2(G_{F,\Sigma_F}; W) = \dim_\Qp \Hom_\Gamma(X, W^\vee(1)).
\]
\end{corollary}

\begin{proof}
The Hochschield--Serre spectral sequence shows that 
\[
H^2(G_{\mathbb{Q},\Sigma}; W) = H^2(G_{F,\Sigma_F}; W)^{G_{\mathbb{Q},\Sigma}/G_{F,\Sigma_F}}\subseteq H^2(G_{F,\Sigma_F}; W),
\]
hence the inequality. The subsequent equality follows from Theorem~\ref{thm:galois_iwasawa}.
\end{proof}

Theorem~\ref{thm:galois_iwasawa} provides us with the first reduction step towards Theorem~\ref{thm:h2-bound}, whose goal is to bound $\dim_\Qp H^2(G_{F,\Sigma_F}; V_{\lambda,s})$ above. Namely, it reduces the problem to the study of $\dim_\Qp\Hom_\Gamma(X,V_{\lambda,2-s})$ by means of Corollary~\ref{cor:galois-iwasawa}. We now carry out the second reduction step, following Coates and Kim as in~\cite{COATES_KIM}, reducing the study of such equivariant $\Hom$-spaces to the study of kernels of Iwasawa operators.

Recall that by Propositions~\ref{prop:fg_X} and~\ref{prop:torsion_X}, $X$ is finitely generated and torsion. Let $x_1,\dots,x_q$ denote an arbitrary set of generators, and let $\delta_1,\dots,\delta_q\in \Lambda$ denote a corresponding set of left annihilators of the $x_i$. By this we mean that for each $1\le i\le q$, one has
\[
\delta_i\cdot x_i = 0,\qquad \delta_i\neq 0.
\]
A $\Gamma$-equivariant morphism naturally extends into a $\Lambda$-linear one, and it follows that any $\Lambda$-linear morphism from $X$ to an arbitrary $\Lambda$-module $W$ must send $x_i$ into $\mathrm{Ker}\delta_i|_{W}$. In particular, evaluation yields a natural embedding,
\begin{equation}\label{eq:hom-to-ker}
\Hom_\Gamma(X, W)\hookrightarrow \bigoplus_{i=1}^q \mathrm{Ker}\,\delta_i|_W.  
\end{equation}

As explained in~\S\ref{sssec:h2-bound-proof-idea}, Theorem~\ref{thm:h2-bound} is an amalgam of two bounds, a trivial and a nontrivial bound. The trivial bound follows immediately from the inclusion~\eqref{eq:hom-to-ker}.
\begin{proposition}[The trivial bound]\label{prop:triv_bound}
For any finite dimensional $\Lambda$-module $W$ of dimension $d$, one has
\[
\dim_\Qp \Hom_\Gamma(X, W) \ll d.
\]
\end{proposition}
\begin{proof}
By~\eqref{eq:hom-to-ker}, one has
\[
\dim_\Qp \Hom_\Gamma(X, W)\le \sum_{i=1}^q \dim_\Qp \mathrm{Ker}\delta_i|_W \le q\cdot d.
\]
\end{proof}

As a first step towards the nontrivial bound, and as the third and final reduction step towards the proof of Theorem~\ref{thm:h2-bound}, we seek to replace the $\delta_i\in \Lambda$ with possibly different, integral Iwasawa-operators, $0\neq \eta_i\in \Omega = \Zp[[H_c]][[T]]$, with the special property that:
\[
0\neq \eta_i(t)\in \Zp[[H_c]],\qquad t\in p\Zp,
\]
and without changing the kernel of $\delta_i$, in the sense that
\[
\mathrm{Ker}\,\delta_i|_{V_{\lambda,s}} = \mathrm{Ker}\,\eta_i|_{V_{\lambda,s}}.
\]
As it turns out, this is possible for any Iwasawa operator $\delta\in \Lambda$, outside a finite set of weights $E = E(\delta)\subseteq \Z$. Possibly replacing $p$ by a $p$-power multiple, we may replace $0\neq \delta\in \Lambda$ with a nonzero element of the integral Iwasawa algebra without altering $\delta$'s kernel. Hence, we may assume that $0\neq \delta\in \Omega$ without loss of generality.

Towards our final reduction step, we first prove:
\begin{proposition}[Noncommutative Weierstrass preparation]\label{prop:wpt}
Let $0\neq \delta\in \Zp[[H_c]][[T]]$ be arbitrary. There is a decomposition
\[
\delta = P\cdot \eta
\]
with $P\in \Zp[T]$, and $\eta \in \Zp[[H_c]][[T]]$, such that $0\neq \eta(t) \in \Zp[[H_c]]$ for every $t\in p\Zp$.
\end{proposition}
\begin{proof}
For every $t_1\in p\Zp$ such that $\delta(t_1) = 0$, we may apply polynomial division in the central coordinate, to rewrite $\delta(T) = (T-t_1)\delta_0(T)$. Indeed, write $\delta(T) = \sum_{n\ge 0} a_n T^n$, with $a_n\in \Zp[[H_c]]$. If $\delta(t_1) = 0$, then $\sum_{n\ge 0}a_nt_1^n = 0$, hence 
\[
\delta(T) = \delta(T) - \delta(t_1) = \sum_{n\ge 0}a_n(T^n - t_1^n) = (T-t_0)\sum_{n\ge 1}a_n\sum_{k=0}^{n-1}T^{n-k-1}t_1^k
\]
\[
= (T-t_1)\sum_{n\ge 0}T^n\sum_{k\ge 1}a_{n+k}t_1^{k-1}
\]
and we set $\delta_0(T) := \sum_{n\ge 0}T^n\sum_{k\ge 1}a_{n+k}t_1^{k-1}$. We observe that the coefficient of $T^n$ is in $\Zp[[H_c]]$. Indeed, this holds since the $a_n$'s are in $\Zp[[H_c]]$, and $t_1\in p\Zp$ has positive $p$-adic valuation. Splitting out linear factors iteratively we obtain a decomposition of the form $P_k\cdot \delta_k$, with 
\[
P_k(T) := (T-t_1)\cdot\cdots\cdot(T-t_k),\qquad \text{and }\delta_k\in \Zp[[H_c]][[T]].
\]
We claim that this process eventually halts. Indeed, let $f: \Zp[[H_c]]\to \Zp$ denote a matrix coefficient, with the property that $f(\delta(T)) := \sum_{n\ge 0}f(a_n)T^n \neq 0$, and divide by a common $p$-power factor, if there is such, so that 
\[
0\neq \overline{f(\delta(T))} = \sum_{n\ge 0}\overline{f(a_n)}T^n \in \Fp[[T]].
\]
Write $\overline{f(\delta(T))} = T^r\cdot \overline{\delta'}(T)$, with $\overline{\delta'}(T)\in \Fp[[T]]$ indivisible by $T$. This is possible since $\overline{f(\delta(T))} \neq 0$ implies that there is some index $n_0\ge 0$, for which $\overline{f(a_{n_0})}\neq 0$, while
\[
\overline{f(a_{n_0})}\neq 0\implies r \le n_0.
\]
Since all roots $t_i$ of $P_k$ have positive $p$-adic valuation, $P_k$ reduces to $T^k$ mod $p$. Furthermore, since $P_k(T) | \delta(T)$, we have $\overline{f(P_k(T))}|\overline{f(\delta(T))}$, and as $f(P_k(T)) = P_k(T)$, it follows that $T^k|\overline{f(\delta(T))}$. In particular, $k\le r$, and therefore splitting linear factors out of $\delta$ must eventually halt. Write $\delta = P\cdot \eta$ for the resulting factorization. By construction, $\eta$ is indivisible by $(T-t)$ for any $t\in p\Zp$. Therefore,
\[
\eta(t)\neq 0,\qquad \text{for all }t\in p\Zp,
\]
which completes the proof.
\end{proof}

We now complete the final reduction step.
\begin{corollary}\label{cor:red-2}
Let $0\neq \delta\in \Omega$ be arbitrary, and let $\delta = P\eta$ denote the decomposition guaranteed by Proposition~\ref{prop:wpt}. Denote by $E = E(\delta) \subseteq \mathbb{Z}$ the set of bad weights, given explicitly via
\[
E := \left\{s\in \Z: P((1+p^c)^s - 1) = 0\right\}.
\]
Then $E$ is finite, and for any $\Lambda$-module $W_s$ of pure weight $s\not\in E$, one has
\[
\mathrm{Ker}\,\delta|_{W_{s}} = \mathrm{Ker}\,\eta|_{W_{s}}.
\]
\end{corollary}

\begin{proof}
The central coordinate $T$ of $\Lambda$ acts on any $W_s$ of pure weight $s\in \Z$ via multiplication by 
\[
t_s := (1 + p^c)^s - 1 \in p\Zp.
\]
Furthermore, since the zero set of $P$ is finite and the map $s\mapsto t_s$ is injective, the subset $E\subseteq \Z$ is finite. Now, assume $s\not\in E$, then $P(t_s)\neq 0$, hence
\[
w\in \mathrm{Ker}\,\delta|_{W_{s}}\iff \delta(t_s)w = 0
\]
\[
\iff P(t_s)\eta(t_s)w = 0\iff \eta(t_s)w = 0 \iff w\in \mathrm{Ker}\,\eta|_{W_{s}}.
\]
\end{proof}

\section{Uniform bounds on Iwasawa operator kernels}\label{sec:bounds}
We maintain the notations of the previous chapters. This chapter is dedicated to the proof of Theorem~\ref{thm:h2-bound}. As explained in~\S\ref{sssec:h2-bound-proof-idea}, Theorem~\ref{thm:h2-bound} is an amalgam of two bounds: the trivial bound, established in Proposition~\ref{prop:triv_bound} above, and the nontrivial bound, given as:
\begin{theorem}[The nontrivial bound]\label{thm:equiv-hom-bound}
Let $0\neq \eta\in \Zp[[H_c]][[T]]$ have the property that $0\neq \eta(t) \in \Zp[[H_c]]$ for every $t\in p\Zp$. Then, for any $(\lambda,s)$-pair,
\[
\dim_\Qp \mathrm{Ker}\,\eta|_{V_{\lambda,s}} \ll \dfrac{(1 + ||\lambda||)^{D-g-1}}{d_\lambda}.
\]
\end{theorem}

This chapter is structured as follows.
\begin{enumerate}
  \item In~\S\ref{ssec:univ_env}, we recall relevant background on the universal enveloping algebra of the symplectic Lie algebra, and define its central reduction.
  \item In~\S\ref{ssec:nontriv}, we establish the nontrivial bound; see Theorem~\ref{thm:equiv-hom-bound}.
  \item In~\S\ref{ssec:h2-bound}, we complete the proof of Theorem~\ref{thm:h2-bound}, which bounds $\dim_\Qp H^2(G_{\mathbb{Q},\Sigma}, V_{\lambda,s})$.
\end{enumerate}


\begin{remark}\label{rem:weyl-bound-comparison}
To appreciate the nontriviality of this bound, we wish to determine when, if at all, it improves on the trivial bound. Up to constants, this occurs whenever
\[
d_\lambda > \dfrac{(1 + ||\lambda||)^{D-g-1}}{d_\lambda}\iff d_\lambda > (1 + ||\lambda||)^{\frac{D-g-1}{2}} = (1 + ||\lambda||)^{g^2 - 1/2}.
\]
Let $e_1,\dots,e_g \in \mathfrak{h}$ denote the standard coordinates of the diagonal maximal torus $T\subset \Sp_{2g,\Qp}$, and denote by $\Phi^+$ the positive roots of $\Sp_{2g,\Qp}$. In the standard coordinates, these are given by
\[
\Phi^+ := \{2e_i:1\le i\le g\} \sqcup \{e_i-e_j: 1\le i < j \le g\}\sqcup \{e_i+e_j: 1\le i < j \le g\}.
\]
Write $\rho := \dfrac{1}{2}\sum_{\alpha \in \Phi^+}\alpha = (g,g-1,\dots, 1)\in \Z^g$. The Weyl dimension formula for the unitary symplectic group $\mathrm{USp}(2g)$ (see~\cite[Corollary~24.6, p.~402]{FULTON_HARRIS}) implies that for any irreducible rational representation of $\Sp_{2g,\Qp}$ of highest weight $\lambda$, one has
\[
d_\lambda = \dim_\Qp V_{\lambda} = \prod_{\alpha \in \Phi^+}\left(1 + \dfrac{\langle\lambda, \alpha\rangle}{\langle\rho,\alpha \rangle}\right).
\]
A by-case analysis on the different positive root types: $2e_i, e_i-e_j$, and $e_i+e_j$, shows that
\[
\dfrac{\langle\lambda, \alpha\rangle}{\langle\rho,\alpha \rangle} \le \lambda_1\le ||\lambda||,\qquad \text{for all }\alpha\in \Phi^+.
\]
Plugging back, we obtain the bound
\[
d_\lambda \le (1 + ||\lambda||)^{g^2}.
\]
The bound derived from Weyl's dimension formula is asymptotically tight, as may be observed by considering $\lambda := m\cdot \rho$ for any $m\ge 1$. In this case,
\[
d_{m\rho} = (1 + m)^{g^2} \gg (1 + ||m\rho||)^{g^2}.
\]
In particular, we observe that the nontrivial bound is asymptotically superior over the trivial bound for the family $V_{m\rho}$, as
\[
(1 + ||m\rho||)^{g^2 - 1/2} = o(d_{m\rho}).
\]
\end{remark}

\subsection{The universal enveloping algebra and its central reduction}\label{ssec:univ_env}
\begin{definition}[The universal enveloping algebra; see~{\cite[I, \S1, p.~18]{MilneLieAlgebras}}]
For a Lie algebra $\mathfrak{g}$, we denote by $\mathcal{U}(\mathfrak{g})$ its universal enveloping algebra. This is the unital associative Lie algebra, given explicitly as the quotient
\[
\mathcal{U}(\mathfrak{g}) := T(\mathfrak{g})/\langle x\otimes y - y\otimes x - [x,y] : x,y\in\mathfrak{g}\rangle,
\]
where $T(\mathfrak{g})$ is the tensor algebra, defined as
\[
T(\mathfrak{g}) := \Zp\oplus \mathfrak{g}\oplus \mathfrak{g}\otimes_\Zp\mathfrak{g}\oplus\dots.
\]
\end{definition}

\begin{notation}
Throughout, let $\mathfrak{h} := \Lie(\Sp_{2g,\Zp})$ denote the Lie algebra of $\Sp_{2g,\Zp}$, and let $\mathcal{U} := \mathcal{U}(\mathfrak{h})$ denote its associated universal enveloping algebra.
\end{notation}

\begin{definition}[The Poincaré--Birkhoff--Witt filtration;
  {\cite[I, \S1, pp.~22--23]{MilneLieAlgebras}}]
The universal enveloping algebra $\mathcal{U}$ is equipped with an increasing, exhaustive filtration, $\mathcal{F}_\bullet$, indexed by the nonnegative integers. For $d\ge 0$, the $\Zp$-submodule $\mathcal{F}_d\,\mathcal{U}$ consists of those elements $u\in \mathcal{U}$, which may be expressed as finite sums of products of at most $d$ elements of $\mathfrak{h}$, where $\mathcal{F}_0\,\mathcal{U} = \Zp$.
\end{definition}

The following claim is an immediate corollary of the PBW Theorem.
\begin{claim}[PBW; {\cite[\S2.3]{ArdakovWadsleyDModulesI}}]\label{claim:pbw}
Write $\mathrm{gr}\,\mathcal{U}$ for the associated graded algebra, then
\[
\mathrm{gr}\,\mathcal{U} \cong S(\mathfrak{h}),
\]
where $S(\mathfrak{h})\cong \Zp[x_1,\dots,x_D]$ denotes the symmetric algebra on $\mathfrak{h}$.
\end{claim}

\begin{definition}[Principal symbol]
For an element $u\in \mathcal{F}_d\,\mathcal{U}$ of PBW degree $d$, we denote by 
\[
\sigma(u) = \sigma_d(u) = u + \mathcal{F}_{d-1}\,\mathcal{U} \in \mathrm{gr}_d\,\mathcal{U} = S(\mathfrak{h})_d
\]
the image of $u$ under the projection from $\mathcal{F}_{d}\,\mathcal{U}$ to $\mathrm{gr}_d\,\mathcal{U}$.
\end{definition}

We record, for later use, the $p$-saturatedness of the PBW filtration on $\mathcal{U}$. First, we recall the definition of $p$-saturation.
\begin{definition}[$p$-saturation; see~{\cite[Ch.~II, \S2, no.~4, Proposition~10(iii) and the paragraph following its proof]{BourbakiCommutativeAlgebra}}]\label{def:p-saturation}
Let $N\subseteq M$ be a pair of $\Zp$-modules. We say that $N$ is $p$-saturated in $M$, if for every $m\in M$, 
\[
pm \in N \implies m\in N.
\]
\end{definition}

\begin{claim}\label{claim:pbw-is-saturated}
The PBW filtration of $\mathcal{U}$ is $p$-saturated.
\end{claim}
\begin{proof}
Assume otherwise, then there is some $d\ge 0$, and a $0\neq u\in \mathcal{U}$, such that
\[
u\not\in \mathcal{F}_d\,\mathcal{U},\qquad\text{while } pu \in \mathcal{F}_d\,\mathcal{U}.
\]
Since the PBW filtration is exhaustive, there is some $e > d$ such that $u\in \mathcal{F}_{e}\,\mathcal{U}\setminus \mathcal{F}_{e-1}\,\mathcal{U}$, and such that 
\begin{equation}\label{eq:saturation}
pu\in \mathcal{F}_{d}\,\mathcal{U} \subset \mathcal{F}_{e-1}\,\mathcal{U}.  
\end{equation}
Consider the principal symbol of $u$,
\[
\sigma_{e}(u) = u + \mathcal{F}_{e-1}\,\mathcal{U} \in \mathrm{gr}_{e}\,\mathcal{U} = S(\mathfrak{h})_{e}.
\]
Since $u\in \mathcal{F}_{e}\,\mathcal{U}\setminus \mathcal{F}_{e-1}\,\mathcal{U}$, $\sigma_{e}(u)\neq 0$. However, as the principal symbol is $\Zp$-linear and $S(\mathfrak{h})_{e}$ is a free $\Zp$-module, we have on the one hand: 
\[
\sigma_e(pu) = p\sigma_e(u) \neq 0 \text{ in } S(\mathfrak{h})_{e},
\]
yet, on the other hand, by~\eqref{eq:saturation},
\[
pu\in \mathcal{F}_{e-1}\,\mathcal{U}\implies \sigma_e(pu) = 0 \text{ in } S(\mathfrak{h})_{e},
\]
a contradiction.
\end{proof}

Denote the center of $\mathcal{U}$ by $\mathcal{Z} := \mathcal{Z}(\mathcal{U})$. Our next goal is to construct a presentation of the algebra $\mathcal{Z}$. This is accomplished in Claim~\ref{claim:center_presentation} below. The next few statements build our path towards it. Write $\mathbf{H} := \Sp_{2g,\Zp}$. Then the adjoint action of $\mathbf{H}$ on $\mathfrak{h}$ extends to give a well-defined action of $\mathbf{H}$ on $\mathcal{U}$. Deriving that action, we get an action of $\mathfrak{h}$ on $\mathcal{U}$.
\begin{proposition}\label{prop:u_center}
One has
\[
\mathcal{Z} = \mathcal{U}^\mathfrak{h} = \mathcal{U}^\mathbf{H}.
\]
\end{proposition}
\begin{proof}
Recall that $\mathbf{H}$ acts on $\mathfrak{h}$ through the adjoint action, $\mathrm{Ad}: \mathbf{H}\to \GL(\mathfrak{h})$. Differentiating this action, one obtains the adjoint action of the Lie algebra on itself, $\mathrm{ad}: \mathfrak{h}\to \End(\mathfrak{h})$, actings as follows. Let $X,Y\in \mathfrak{h}$, then $\mathrm{ad}(X)(Y) := [X,Y]$. This action extends to $\mathcal{U}$ as a derivation, in the sense that it is linear, and for $X_1\cdot \dots\cdot X_r\in\mathcal{U}$, one has
\[
\mathrm{ad}(X)(X_1\cdot \dots\cdot X_r) = \sum_{i=1}^r X_1\cdot \dots \cdot X_{i-1} \cdot [X,X_i] \cdot X_{i+1}\cdot \dots \cdot X_r.
\]
The $\mathfrak{h}$-invariants of $\mathcal{U}$ are those $u\in \mathcal{U}$, for which $[h,u] = hu - uh = 0$ for every $h\in \mathfrak{h}$. Since $\mathfrak{h}$ generates $\mathcal{U}$, $\mathcal{Z} = \mathcal{U}^\mathfrak{h}$. It remains to prove that $\mathcal{U}^\mathfrak{h} = \mathcal{U}^\mathbf{H}$. Observing that there is a natural inclusion $\mathcal{U}^\mathbf{H}\subseteq \mathcal{U}^\mathfrak{h}$ obtained by differentiation, it remains to establish the reverse inclusion.

Write $\mathcal{U}_\Qp$, $\mathfrak{h}_{\Qp}$ and $\mathbf{H}_\Qp$ for the associated base changes to $\Qp$. According to~\cite[II, discussion following~(70), p.~132]{MilneLAG}, since $\mathbf{H}_\Qp$ is connected over a field of characteristic zero, we have a rational version of our desired identification:
\[
\mathcal{U}_{\Qp}^{\mathfrak{h}_{\Qp}} = \mathcal{U}_{\Qp}^{\mathbf{H}_{\Qp}}.
\]
Note that the $\mathbf{H}$ action on $\mathcal{U}$ is encoded by the coaction map, 
\[
\delta: \mathcal{U} \to \mathcal{U}\otimes \mathcal{O}(\mathbf{H}),
\]
and we have
\[
u\in \mathcal{U}^{\mathbf{H}}\iff \delta(u) - u\otimes 1 = 0.
\]
Let $u\in \mathcal{U}^\mathfrak{h}$, our goal is to show that $u\in \mathcal{U}^\mathbf{H}$, i.e.~that $\delta(u) - u\otimes 1 = 0$. Since $\mathcal{U}_{\Qp}^{\mathfrak{h}_{\Qp}} = \mathcal{U}_{\Qp}^{\mathbf{H}_{\Qp}}$, and the $\Qp$-coaction map is the $\Qp$-basechange of the $\Zp$-coaction map, we know that
\[
\delta(u) - u\otimes 1 = 0 \in \mathcal{U}_\Qp\otimes \mathcal{O}(\mathbf{H}_{\Qp}),
\] 
hence, there is some power of $p$, say $m\ge 0$, such that
\begin{equation}\label{eq:p-torsion}
p^m(\delta(u) - u\otimes 1) = 0 \in \mathcal{U}\otimes \mathcal{O}(\mathbf{H}).  
\end{equation}
By the integral PBW theorem, $\mathcal{U}$ is $p$-torsion free. Consider the multiplication-by-$p$ left exact sequence,
\[
0\to \mathcal{U}\xrightarrow{\cdot p}\mathcal{U},
\]
Since $\mathbf{H}$ is smooth over $\Zp$, $\mathcal{O}(\mathbf{H})$ is flat. Therefore, tensoring the above sequence with $\mathcal{O}(\mathbf{H})$ is exact, and we have a left exact sequence
\[
0\to \mathcal{U}\otimes\mathcal{O}(\mathbf{H}) \xrightarrow{\cdot p}\mathcal{U}\otimes\mathcal{O}(\mathbf{H}).
\]
In particular, $\mathcal{U}\otimes\mathcal{O}(\mathbf{H})$ is $p$-torsion free, and therefore~\eqref{eq:p-torsion} implies that
\[
\delta(u) - u\otimes 1 = 0 \in \mathcal{U}\otimes \mathcal{O}(\mathbf{H}),
\]
and hence, that $u\in \mathcal{U}^\mathbf{H}$, which proves the reverse inclusion, and completes the proof.
\end{proof}
The following discussion closely follows~\cite[\S4.10, \S6.9]{ArdakovWadsleyIrreducible}. Write 
\begin{equation}\label{eq:triangular_decomp}
\mathfrak{h} = \mathfrak{n}^-\oplus \mathfrak{t}\oplus \mathfrak{n}^+  
\end{equation}
for the triangular/root-space decomposition of $\mathfrak{h}$, with $\mathfrak{n}^-$ and $\mathfrak{n}^+$ spanned by the negative and positive roots, respectively, and with $\mathfrak{t}$ the Cartan Lie subalgebra. The enveloping algebra admits a description as
\[
\mathcal{U} = \mathcal{U}(\mathfrak{n}^-)\otimes \mathcal{U}(\mathfrak{t})\otimes \mathcal{U}(\mathfrak{n}^+),
\]
giving rise to a direct sum decomposition
\begin{equation}\label{eq:direct_sum_decomp}
\mathcal{U} = \mathcal{U}(\mathfrak{t})\oplus \left(\mathfrak{n}^-\cdot \mathcal{U}\oplus \mathcal{U}\cdot \mathfrak{n}^+\right).
\end{equation}
We also have a version of~\eqref{eq:direct_sum_decomp} for the associated gradeds,
giving rise to the following decomposition of the symmetric algebra; cf.~Claim~\ref{claim:pbw}:
\begin{equation}\label{eq:direct_sum_decomp2}
S(\mathfrak{h}) = S(\mathfrak{t})\oplus \left(\mathfrak{n}^-\cdot S(\mathfrak{h})\oplus S(\mathfrak{h})\cdot \mathfrak{n}^+\right).
\end{equation}
Note that the Weyl group of $\mathbf{H}$, denoted $\mathcal{W}$, acts on $S(\mathfrak{t}) = \Zp[t_1,\dots,t_g]$ by signed permutations. Write,
\begin{equation}\label{eq:hc_projection}
\mathrm{pr}: \mathcal{U}\twoheadrightarrow \mathcal{U}(\mathfrak{t}),\qquad \mathrm{gr}\,\mathrm{pr}: S(\mathfrak{h})\to S(\mathfrak{t})
\end{equation}
for the $\Zp$-module surjections induced from the direct sum decompositions of~\eqref{eq:direct_sum_decomp} and~\eqref{eq:direct_sum_decomp2}. 

Since taking invariants is left exact, there is a natural inclusion 
\[
\iota: \mathrm{gr}(\mathcal{U}^\mathbf{H})\hookrightarrow S(\mathfrak{h})^{\mathbf{H}}.
\]
Denote the restriction of $\mathrm{gr}\,\mathrm{pr}$ to $S(\mathfrak{h})^\mathbf{H}$ by $\psi$. By~\cite[\S6.9]{ArdakovWadsleyIrreducible}, the map
\[
\psi: S(\mathfrak{h})^\mathbf{H} \to S(\mathfrak{t})
\]
is injective, and factors through $S(\mathfrak{t})^\mathcal{W}$. 

\begin{theorem}[{\cite[\S6.9, Proposition in p.~507]{ArdakovWadsleyIrreducible}}]\label{thm:aw_main}
The rows of the diagram
\begin{equation}\label{eq:diagram_aw}
\begin{tikzcd}
0\arrow[r]&\mathrm{gr}(\mathcal{U}^\mathbf{H})\arrow[r,"\cdot p"]\arrow[d,"\iota"]&\mathrm{gr}(\mathcal{U}^\mathbf{H})\arrow[r]\arrow[d, "\iota"]&\mathrm{gr}(\mathcal{U}_k^{\mathbf{H}_k})\arrow[r]\arrow[d,"\iota_k"]&0\\
0\arrow[r]&S(\mathfrak{h})^\mathbf{H}\arrow[r,"\cdot p"]\arrow[d,"\psi"]&S(\mathfrak{h})^\mathbf{H}\arrow[r]\arrow[d, "\psi"]&S(\mathfrak{g}_k)^{\mathbf{H}_k}\arrow[r]\arrow[d,"\psi_k"]&0\\
0\arrow[r]&S(\mathfrak{t})^\mathcal{W}\arrow[r,"\cdot p"]&S(\mathfrak{t})^\mathcal{W}\arrow[r]&S(\mathfrak{t}_k)^{\mathcal{W}_k}\arrow[r]&0,
\end{tikzcd}
\end{equation}
are exact, and each vertical map is an isomorphism.
\end{theorem}
\begin{remark}[Very good primes]\label{rem:very-good-primes}
Theorem~\ref{thm:aw_main} holds under the assumption that the auxiliary prime $p$ is \emph{very good} for the root system of $\mathbf{H}$. Since $\mathbf{H} = \Sp_{2g,\Zp}$ has root system of type $C_g$, every odd prime is very good; see~\cite[\S6.8]{ArdakovWadsleyIrreducible}. Since by assumption $p$ is odd, the very good prime condition is superfluous to us.
\end{remark}

The left vertical column of isomorphisms in Theorem~\ref{thm:aw_main} reads as,
\[
\mathrm{gr}\mathcal{Z} \cong S(\mathfrak{h})^\mathbf{H}\cong S(\mathfrak{t})^\mathcal{W}.
\]
Furthermore, since $S(\mathfrak{t}) \cong \Zp[t_1,\dots,t_g]$ and $\mathcal{W}$ acts via signed permutations, $S(\mathfrak{t})^\mathcal{W}$ is the $\Zp$-subalgebra generated by the symmetric polynomials in the $t_i^2$, and therefore has a unique generator in each even degree $i\in 2,\dots,2g$. Since $S(\mathfrak{h})^\mathbf{H}\cong S(\mathfrak{t})^\mathcal{W}$, there are generators $q_1,\dots, q_{g}\in S(\mathfrak{h})^\mathbf{H}$, with $q_i\in S(\mathfrak{h})_{2i}$, such that
\begin{equation}\label{eq:invariant-polynomial-generators}
S(\mathfrak{h})^\mathbf{H} = \Zp[q_1,\dots q_g].
\end{equation}

\begin{claim}\label{claim:center_presentation}
There are central generators $c_1,\dots, c_{g}\in \mathcal{Z}$ with $c_i \in \mathcal{F}_{2i}\,\mathcal{U}$ such that
\[
\mathcal{Z} = \Zp[c_1,\dots,c_g],
\]
and for each $i$, one has $\sigma(c_i) = q_i$.
\end{claim}
\begin{proof}
The isomorphism $\iota$ of~\eqref{eq:diagram_aw} reads as
\[
\mathrm{gr}\mathcal{Z}\cong S(\mathfrak{h})^\mathbf{H} = \Zp[q_1,\dots q_g].
\]
Choose, for each $1\le i\le g$, a lift $c_i \in \mathcal{Z}\cap \mathcal{F}_{2i}\,\mathcal{U}$ of $q_i$, in the sense that $\sigma(c_i) = q_i$. By construction, $\langle c_1,\dots,c_g\rangle\subseteq \mathcal{Z}$. We argue that $\mathcal{Z} = \Zp[c_1,\dots,c_g]$ in two steps, first showing that $\langle c_1,\dots,c_g\rangle$ generates $\mathcal{Z}$, and then establishing that there are no nontrivial relations between the $c_i$.

\textbf{Generation.} Let $c\in \mathcal{Z}$ be arbitrary, and of PBW degree $d$. Since $\sigma(c) \in S(\mathfrak{h})^{\mathbf{H}}_d$ and $S(\mathfrak{h})^{\mathbf{H}}\cong \Zp[q_1,\dots,q_g]$, there is a homogeneous polynomial $P := P(x_1,\dots,x_g)$ of degree $d$, with the property that $P(q_1,\dots,q_g) = \sigma(c)$. Consider $c' := P(c_1,\dots,c_g)\in \langle c_1,\dots,c_g\rangle$. By construction, $\sigma(c') = \sigma(c)$. Hence, $c-c'\in \mathcal{Z}\cap \mathcal{F}_{d-1}\mathcal{U}$. Continuing this process recursively, establishes that $\langle c_1,\dots,c_g\rangle = \mathcal{Z}$.

\textbf{Independence.} Assume towards contradiction that there is a non-trivial algebraic relation amongst the $c_i$'s. The principal symbol of this relation would yield a non-trivial relation amongst the $q_i$'s, which is impossible. Therefore, $\mathcal{Z}\cong \Zp[c_1,\dots,c_g]$, which completes the proof.
\end{proof}

Next, we want to define the central reductions $\mathcal{U}\twoheadrightarrow R_\mu$. Before doing so, we require the definition of the Harish-Chandra homomorphism, and the central character.
\begin{definition}[The Harish-Chandra homomorphism; see~{\cite[\S4.10]{ArdakovWadsleyIrreducible}}]\label{def:harish-chandra-map}
The Harish-Chandra homomorphism is the restriction of the projection $\mathrm{pr}: \mathcal{U}\to \mathcal{U}(\mathfrak{t})$ of~\eqref{eq:hc_projection}, to $\mathcal{Z}$, i.e.~$\mathrm{HC} := \mathrm{pr}|_{\mathcal{Z}}$,
\[
\mathrm{HC}: \mathcal{Z}\to \mathcal{U}(\mathfrak{t}).
\]
\end{definition}

\begin{remark}[The Harish-Chandra homomorphism]
While a priori only a $\Zp$-module map, the Harish-Chandra homomorphism is a homomorphism of $\Zp$-algebras; see~\cite[especially the lemma at the end of \S4.10]{ArdakovWadsleyIrreducible}.
\end{remark}

\begin{definition}[The central character]\label{def:central_character}
Let $\mu\in \mathbb{Z}_p^g$ denote a highest weight for $\Sp_{2g,\Zp}$, and let $\mathrm{ev}_\mu: S(\mathfrak{t})\to \Zp$ denote the evaluation map. In torus coordinates, $S(\mathfrak{t}) = \Zp[t_1,\dots,t_g]$, this is the map that sends $t_i$ to $\mu_i$. We define the central character associated to the highest weight $\mu$, denoted $\zeta_\mu$, as the composition,
\[
\zeta_\mu = \mathrm{ev}_\mu\circ \mathrm{HC}.
\]
\end{definition}

\begin{definition}[The central reduction]
Let $\mu\in \mathbb{Z}_p^g$ denote a highest weight vector. We define the central reduction of $\mathcal{U}$ (at the highest weight $\mu$) to be the algebra,
\[
R_\mu := \mathcal{U}/(c_i - \zeta_\mu(c_i): 1\le i\le g).
\]
\end{definition}

For later use, we record the following fact about the mod-$p$ reduction of $R_\mu$.
\begin{proposition}\label{prop:r-mu-reduction-domain}
$R_\mu/pR_\mu$ is a domain.
\end{proposition}
\begin{proof}
According to Ardakov and Wadsley in~\cite[\S4.6, Lemma~(b)]{ArdakovWadsleyVerma}, the algebra
\[
S(\mathfrak{h}_k)/(\overline{q_1},\dots,\overline{q_g})
\]
is a domain, which implies that so is $R_\mu/pR_\mu$. Indeed, if $x,y\in R_\mu/pR_\mu$ are a pair of non-zero elements with necessarily non-zero principal symbols $\sigma(x)$ and $\sigma(y)$ in $S(\mathfrak{h}_k)/(\overline{q_1},\dots,\overline{q_g})$. Then $\sigma(xy) = \sigma(x)\cdot \sigma(y)$ is non-zero since $S(\mathfrak{h}_k)/(\overline{q_1},\dots,\overline{q_g})$ is a domain, and therefore $xy\neq 0$, since its principal symbol does not vanish, which completes the proof.
\end{proof}

\begin{notation}
For a $\Zp$-module, or an algebra $A$, we denote by $\widehat{A} := \varprojlim_n A/p^n$ its pro-$p$ completion.
\end{notation}
We shall be particularly interested in the pro-$p$ completions of $\mathcal{U}$ and $R_\mu$, denoted $\widehat{\mathcal{U}}$ and $\widehat{R}_\mu$. The former is directly related to the Iwasawa algebra through the Lazard comparison map.
\begin{definition}[Lazard's comparison map; see~{\cite[proof of Theorem~10.4]{ArdakovWadsleyIrreducible}}]
Let $a\ge 1$, and write $H^{\mathrm{amb}}_a := \mathrm{Ker}\left(\Sp_{2g}(\Zp)\to \Sp_{2g}(\Zp/p^a)\right)$ for the level $a$ congruence subgroup of the ambient symplectic group. There is a well-defined map, which we shall refer to as the Lazard comparison,
\[
\iota_a: \Zp[[H^{\mathrm{amb}}_a]]\to \widehat{\mathcal{U}}_{\Qp},\qquad \gamma\mapsto \exp(\log \gamma).
\]
\end{definition}

\begin{proposition}
For every $a\ge 1$, the Lazard map $\iota_a$ factors through the integral completed enveloping algebra, $\widehat{\mathcal{U}}$.
\end{proposition}
\begin{proof}
First, we observe that since $p$ is odd and $a\ge 1$, the matrix $p$-adic logarithm gives a bijection between $H^{\mathrm{amb}}_a$ and $p^a\mathfrak{h}$. In particular, 
\[
\log(H^{\mathrm{amb}}_a) = p^a\mathfrak{h}.
\]
The logarithm on the Iwasawa algebra $\Zp[[H^{\mathrm{amb}}_a]]$ is extended continuously and $\Zp$-linearly from $\Zp[H^{\mathrm{amb}}_a]$, and we wish to argue that $\log(\Zp[[H^{\mathrm{amb}}_a]]) = p^a\mathfrak{h}$ as well. We have the inclusions:
\[
p^a\mathfrak{h} = \log(\Zp[H^{\mathrm{amb}}_a])\subseteq \log(\Zp[[H^{\mathrm{amb}}_a]]) \subseteq \overline{\log(\Zp[H^{\mathrm{amb}}_a])} = p^a\mathfrak{h},
\]
where the last inclusion follows from the fact that the logarithm sends convergent sequences to convergent sequences, and therefore
\[
\log(\Zp[[H^{\mathrm{amb}}_a]]) = \log(\overline{\Zp[H^{\mathrm{amb}}_a]}) \subseteq \overline{\log(\Zp[H^{\mathrm{amb}}_a])}.
\]
Since $\mathrm{Im}(\iota_a)$ coincides with $\exp(p^a\mathfrak{h})$, it suffices to verify that $\exp: p^a\mathfrak{h}\to \widehat{\mathcal{U}}_{\Qp}$ factors through $\widehat{\mathcal{U}}$. Let $X\in p^a\mathfrak{h}$. By assumption, there exists a $Y \in \mathfrak{h}$ such that $X = pY$. Then
\[
\exp(X) := \sum_{n\ge 0}\dfrac{1}{n!}\cdot X^n = \sum_{n\ge 0}\dfrac{p^n}{n!}Y^n \in \widehat{\mathcal{U}}_\Qp.
\]
To verify that $\exp(X)\in \widehat{\mathcal{U}}$, it suffices to verify that $v_p\left(\dfrac{p^n}{n!}\right) = n - v_p(n!) \ge 0$ for every $n\ge 0$. This obviously holds for $n = 0$. For $n\ge 1$ we have, by Legendre's formula,
\[
v_p(n!) = \sum_{k\ge 1}\left[\dfrac{n}{p^k}\right] = \dfrac{n - s_p(n)}{p-1},
\]
where $s_p(n)$ is the sum of the base-$p$ digits of $n$. In particular, since $n\ge 1$, $s_p(n) \ge 1$, hence $v_p(n!) \le \dfrac{n-1}{p-1}$. Plugging back, we obtain for all $n\ge 1$:
\[
v_p\left(\dfrac{p^n}{n!}\right) \ge n - \dfrac{n-1}{p-1} \ge 1,
\]
which completes the proof.
\end{proof}

\begin{remark}
Since $\Zp[[H_c]]$ naturally embeds into $\Zp[[H_c^{\mathrm{amb}}]]$, we get a map from our integral Iwasawa algebra into the integral, $p$-completed, universal enveloping algebra $\widehat{\mathcal{U}}$.
\end{remark}

Towards studying the kernel of $0\neq \eta \in \Zp[[H_c]]$ acting on $V_{\lambda,s}$, we require the affinoid Verma faithfulness theorem of Ardakov and Wadsley; see Theorem~\ref{thm:affinoid_verma_faithfulness} below. To state it, we need the definition of the affinoid Verma module of highest weight $\mu\in\mathbb{Z}_p^g$.

\begin{definition}[Affinoid Verma module; see~{\cite[\S4.5]{ArdakovWadsleyVerma}}]\label{def:affinoid_verma_module}
Let $\mathfrak{h} = \mathfrak{n}^-\oplus \mathfrak{t}\oplus \mathfrak{n}^+$ denote the triangular decomposition of~\eqref{eq:triangular_decomp}, and fix $\mathfrak{b}^+ := \mathfrak{t}\oplus \mathfrak{n}^+$ for the upper triangular Borel. Let $\mu\in\mathbb{Z}_p^g$ denote a $g$-tuple of $p$-adic integers, which we shall refer to also as a highest weight vector. Regard $\mu$ as a functional, $\mu: \mathfrak{t}\to \Zp$, and define the rank one $\mathfrak{b}^+$-module $\mathbf{Z}_{p,\mu}$ by
\[
t\cdot 1 = \mu(t),\qquad x\cdot 1 = 0\quad (t\in \mathfrak{t}, x\in \mathfrak{n}^+).
\]
This extends to a $\mathcal{\mathcal{U}(\mathfrak{b}^+)}$-module naturally. The integral affinoid Verma module (of highest weight $\mu$) is:
\[
\widehat{V^\mu} := \widehat{\mathcal{U}(\mathfrak{h})}\widehat{\otimes}_{\widehat{\mathcal{U}(\mathfrak{b}^+)}}\mathbf{Z}_{p,\lambda}.
\]
\end{definition}

The affinoid Verma module has a natural action of $\widehat{\mathcal{U}} = \widehat{\mathcal{U}(\mathfrak{h})}$ by left translation, and equips an action of the Iwasawa algebra through the Lazard comparison (see~\cite[proof of Theorem~5.3]{ArdakovWadsleyVerma}). More specifically, with $\iota_a: \Zp[[H_a^{\mathrm{amb}}]]\to \widehat{\mathcal{U}}$ the Lazard comparison of level $a\ge 1$, and $u\widehat{\otimes} 1\in \widehat{V^\mu}$. One has, for every $\nu\in \Zp[[H_a^{\mathrm{amb}}]]$,
\[
\nu\cdot (u\widehat{\otimes} 1) = (\iota_a(\nu)u)\widehat{\otimes} 1.
\]

\begin{theorem}[Affinoid Verma faithfulness; see~{\cite[Theorem~5.4]{ArdakovWadsleyVerma}}]\label{thm:affinoid_verma_faithfulness}
For every $a\ge 1$, and highest weight $\mu\in \mathbb{Z}_p^g$, the affinoid Verma module $\widehat{V^\mu}$ is faithful as a $\Zp[[H_a^{\mathrm{amb}}]]$-module.
\end{theorem}

\begin{proposition}\label{prop:affinoid_center}
The center $\mathcal{Z}$ of $\mathcal{U}$ acts on $\widehat{V^\mu}$ through the central character $\zeta_\mu$.
\end{proposition}
\begin{proof}
By~\cite[\S4.6]{ArdakovWadsleyVerma}, the center, $\mathcal{Z}$ of $\mathcal{U}$ acts on $\widehat{V^\mu}$ through a character $\chi_\mu : \mathcal{Z}\to \Qp$. Our aim is to verify that this character $\chi_\mu$ is indeed the central character $\zeta_\mu$ of Definition~\ref{def:central_character}. The Harish-Chandra homomorphism, $\mathrm{HC}$ of Definition~\ref{def:harish-chandra-map} is the restriction of the canonical projection, $\mathcal{U}\to \mathcal{U}(\mathfrak{t})$ to the center $\mathcal{Z}$. Let $z\in \mathcal{Z}$, and decompose it with respect to the direct sum decomposition of~\eqref{eq:direct_sum_decomp}:
\[
z = \mathrm{HC}(z) + n^-u_1 + u_2n^+,
\]
where $u_1,u_2\in\mathcal{U}$, and $n^-$ and $n^+$ belong to $\mathfrak{n}^-$ and $\mathfrak{n}^+$, respectively. Now, consider the action of $z$ on the highest weight vector $v_\mu := 1\widehat{\otimes} 1\in \widehat{V^\mu}$. By assumption, $zv_\mu = \chi_{\mu}(z)v_\mu$, while by the definition of the affinoid Verma module, $\mathfrak{n}^+$ annihilates $v_\mu$, and $\mathrm{HC}(z)\in \mathcal{U}(\mathfrak{t}) = S(\mathfrak{t})$ acts through the central character $\zeta_\mu := \mathrm{ev}_\mu\circ \mathrm{HC}$. In particular, we have 
\[
\mathrm{HC}(z)v_\mu = \zeta_\mu(z)v_\mu.
\]
Combining with the above, it follows that
\begin{equation}\label{eq:v-mu-equality}
\chi_\mu(z)v_\mu = \zeta_\mu(z)v_\mu + n^-u_1v_\mu.  
\end{equation}
Observe that, since $n^-u_1v_\mu$ has strictly smaller weight than that of $v_\mu$, the only way equality~\eqref{eq:v-mu-equality} can hold is if $n^-u_1v_\mu = 0$. Consequentially, $\chi_\mu(z) = \zeta_\mu(z)$, and since this holds for every $z\in \mathcal{Z}$, this completes the proof.
\end{proof}

We obtain the following key corollary.
\begin{corollary}\label{cor:eta_reduction}
For any level parameter $a\ge 1$ and highest weight parameter $\mu\in \mathbb{Z}_p^g$, the composition
\[
\pi_\mu: \Z_p[[H_a^{\mathrm{amb}}]]\xrightarrow{\iota_a}\widehat{\mathcal{U}}\twoheadrightarrow \widehat{R}_\mu
\]
is injective.
\end{corollary}
\begin{proof}
Ardakov--Wadsely's affinoid Verma faithfulness~\ref{thm:affinoid_verma_faithfulness} implies that the Lazard comparison is injective. Indeed, this follows from the fact that the ambient Iwasawa algebra acts faithfully on the affinoid Verma module $\widehat{V^\mu}$. Since, by Proposition~\ref{prop:affinoid_center}, this action factors through the central reduction, the composed map, $\Z_p[[H_a^{\mathrm{amb}}]]\to \widehat{R}_\mu$ is injective as well.
\end{proof}

\subsection{The nontrivial bound}\label{ssec:nontriv}
We now work our way towards bounding $\dim_\Qp\mathrm{Ker}\eta|_{V_{\lambda,s}}$. This subsection is structured as follows.
\begin{enumerate}
  \item In~\S\ref{sssec:gradeds}, we bound $\dim_\Qp\mathrm{Ker}\eta|_{V_{\lambda,s}}$ in terms of the dimensions of the associated gradeds of the mod $p$ reduction of $R_\lambda$; see Proposition~\ref{prop:key_prop}.
  \item In~\S\ref{sssec:nontriv_proof}, we bound the dimensions of these associated gradeds using Hilbert polynomials, and deduce Theorem~\ref{thm:equiv-hom-bound}.
\end{enumerate} 

\subsubsection{\texorpdfstring{Bounding $\dim_\Qp\mathrm{Ker}\eta|_{V_{\lambda,s}}$ in terms of the $\Fp$-dimensions of $\mathrm{gr}_m R_\lambda/pR_\lambda$}{Bounding kernel dimensions in terms of dimensions of graded pieces modulo p}}\label{sssec:gradeds}
For $t\in p\Zp$ and $\mu\in \mathbb{Z}_p^g$, write $\eta_{\mu, t} := \pi_\mu(\eta(t)) \in \widehat{R_\mu}$ for the central reduction of $\eta(t)$ at the highest weight parameter $\mu\in \mathbb{Z}^g$. According to Corollary~\ref{cor:eta_reduction}, $\eta_{\mu,t} \neq 0$.

\begin{notation}[Weight space]
We shall refer to the product $\mathbb{Z}_p^g\times p\Zp$ as the weight space. It parameterizes tuples $(\mu,t)$ with $t\in p\mathbb{Z}_p$ and $\mu\in \mathbb{Z}_p^g$ a highest weight vector.
\end{notation}

The following key lemma allows us to bound the PBW filtration degree of the mod-$p$ reduction of $\eta_{\mu,t}$, upon clearing out its maximal $p$-power divisor.
\begin{lemma}\label{lem:key_lemma}
Let $(\mu,t)\in \mathbb{Z}_p^g\times p\Zp$ be a point of the weight space, and let
\[
0\neq \eta_{\mu,t} := \pi_\mu(\eta(t)) \in \widehat{R}_\mu
\]
denote the corresponding central reduction. Let $v := v(\mu,t) \ge 0$ be the largest integer with the property that 
\[
\theta_{\mu,t} := p^{-v}\eta_{\mu,t} \in \widehat{R}_\mu.
\]
Write $\overline{\theta_{\mu,t}}\in \overline{R_\mu}$ for $\theta_{\mu,t}$'s reduction mod-$p$. Equip $\overline{R_\mu}$ with its quotient PBW filtration from $\mathcal{U}$. Then, there exists an integer $d\ge 0$ independent of $(\mu,t)$, such that
\[
\overline{\theta_{\mu,t}} \in \mathcal{F}_d\overline{R_\mu}.
\]
\end{lemma}

\begin{proof}
Denote by $E_n \subseteq \mathbb{Z}_p^g\times p\Zp$ the weight subspace given explicitly as
\[
E_n := \left\{(\mu,t)\in \mathbb{Z}_p^g\times p\Zp: p^{-n}\eta_{\mu,t}\not\in \widehat{R}_{\mu}\right\}.
\]
We argue that the $E_n$ form an ascending chain of open subsets of the weight space in the $p$-adic topology, which cover the weight space, in the sense that $\cup_n E_n = \mathbb{Z}_p^g\times p\Zp$.

\textbf{Ascension.} Equivalently, it suffices to prove that $(\mu,t)\not\in E_{n+1}\implies (\mu,t)\not\in E_n$. Indeed, by definition, $(\mu,t)\not\in E_{n+1}$ implies that $p^{-n-1}\eta_{\mu,t}\in \widehat{R}_{\mu}$. However, since $\widehat{R}_{\mu}$ is a $\Zp$-module,
\[
p^{-n-1}\eta_{\mu,t}\in \widehat{R}_{\mu} \implies p^{-n}\eta_{\mu,t}\in \widehat{R}_{\mu},
\]
hence $(\mu,t)\not\in E_{n}$, which shows that $E_{n+1}\supseteq E_n$.

\textbf{Covering.} A $\Zp$-module $M$ is said to be $p$-separated if $\cap_{m\ge 1}p^mM = 0$. The pro-$p$ completion of a module is always $p$-separated, hence $\widehat{R}_\mu$ is such. Let $(\mu,t)$ be a point in the weight space, and assume $(\mu,t)\not\in \cup_{n\ge 1}E_n$. Then, by definition of the $E_n$, 
\[
(\mu,t)\not\in \cup_{n\ge 1}E_n\iff \eta_{\mu,t} \in \cap_{n\ge 1}p^n\widehat{R}_{\mu} = 0,
\]
however, by Corollary~\ref{cor:eta_reduction}, $\eta_{\mu,t}\neq 0$, a contradiction. Therefore, $\cup_{n\ge 1}E_n = \mathbb{Z}_p^g\times p\Zp$.

\textbf{Openness.} Let $(\mu,t) \in E_n$. Our goal is to find an open neighborhood $(\mu,t)\in \mathcal{O}\subseteq E_n$. By definition, since $(\mu,t)\in E_n$, we have $\eta_{\mu,t}\not\in p^n\widehat{R}_\mu$, hence
\[
0\neq \eta_{\mu,t} + p^n\widehat{R}_\mu \in \widehat{R}_\mu/p^n\widehat{R}_\mu = R_\mu/p^nR_\mu.
\]
In particular, since $\widehat{R}_\mu$ is a quotient of $\widehat{\mathcal{U}}$,
\[
0\neq \iota_c(\eta(t)) + p^n\widehat{\mathcal{U}} \in \widehat{\mathcal{U}}/p^n\widehat{\mathcal{U}} = \mathcal{U}/p^n\mathcal{U}.
\]
Since $\iota_c(\eta(T)): p\Zp \to \widehat{\mathcal{U}}$ is continuous, and the reduction $\widehat{\mathcal{U}}\to \widehat{\mathcal{U}}/p^n\widehat{\mathcal{U}}$ is continuous, the composed function is a continuous function from a compact domain to a discrete domain, hence locally constant. Therefore, there is some open $D\subseteq p\Zp$ for which
\[
\iota_c(\eta(T)) = \iota_c(\eta(t)) \mod p^n\widehat{\mathcal{U}},\qquad \text{for all }T\in D.
\]
Next, since $R_\mu$ is obtained from $\mathcal{U}$ through a quotient by the ideal 
\[
I(\mu) := (c_i - \zeta_\mu(c_i): 1\le i\le g)\subseteq \mathcal{U},
\]
and since the central character is a continuous function of the weight, $\mu$. There exists a sufficiently small neighborhood $N\subseteq \mathbb{Z}_p^g$ of $\mu$, such that
\[
I(\mu') + p^n\widehat{\mathcal{U}} = I(\mu) \mod p^n\widehat{\mathcal{U}},\qquad \mu'\in N.
\]
In particular, for such $\mu'\in N$, the quotients $R_\mu/p^nR_\mu$ and $R_{\mu'}/p^nR_{\mu'}$ are isomorphic. Let $(\mu',t')\in \mathcal{O} := N\times D$. Under the isomorphism $R_\mu/p^nR_\mu\cong R_{\mu'}/p^nR_{\mu'}$, one has $\eta_{\mu,t} = \eta_{\mu',t'}$ for every $(\mu',t')\in \mathcal{O}$, as both elements are reductions of the same element $\iota_c(\eta(t)) \mod p^n\widehat{\mathcal{U}}$ along the same ideal $I(\mu) \mod p^n\widehat{\mathcal{U}}$. In particular, it follows that for every $(\mu',t')\in \mathcal{O}$,
\[
0\neq \eta_{\mu',t'} \mod p^n\widehat{R}_{\mu'} \in R_{\mu'}/p^nR_{\mu'},
\]
which implies that $p^{-n}\eta_{\mu',t'}\not\in \widehat{R}_{\mu',t'}$, hence $(\mu',t')\in E_n$, and therefore $\mathcal{O}\subseteq E_n$.

The compactness of the weight space now implies that there is some $N\ge 1$, for which $E_N = \mathbb{Z}_p^g\times p\Zp$. Fix such $N\ge 1$, and consider the map $p\Zp\to \mathcal{U}/p^N\mathcal{U}$, given by
\begin{equation}\label{eq:the_above_map}
T\mapsto \iota_c(\eta(T)) \mod p^N\widehat{\mathcal{U}}.  
\end{equation}
Since this map is continuous from a compact domain to a discrete domain, it is locally constant, hence takes finitely many values. Choose representatives of these values, i.e.~choose 
\[
u_1,\dots, u_r\in \mathcal{U},
\]
such that their reductions mod $p^N\mathcal{U}$ cover the image of the map~\eqref{eq:the_above_map}. For each $1\le i\le r$, let $d_i$ denote the PBW length of $u_i$, and fix $d = \max_i d_i$.

Now, for each $(\mu,t)\in \mathbb{Z}_p^g\times p\Zp$, denote by $0\le v = v(\mu,t) < N$ the maximal integer power of $p$ for which $p^{-v}\eta_{\mu,t} \in \widehat{R}_\mu$. Fix $\theta_{\mu,t} := p^{-v}\eta_{\mu,t}$. Since $p^{-v-1}\eta_{\mu,t}\not\in \widehat{R}_\mu$, 
\[
0\neq \overline{\theta_{\mu,t}} \in \overline{R_\mu} := R_\mu/pR_\mu.
\]
Since the map~\eqref{eq:the_above_map} takes the finitely many values represented by the $(u_i)_{i=1}^r$, there is some $1\le j\le r$, for which
\[
\iota_c(\eta(t)) = u_j \mod p^N\widehat{\mathcal{U}}.
\]
Hence, $p^{-v}\iota_c(\eta(t)) = p^{-v}u_j \mod p^{N-v}\widehat{\mathcal{U}}$, and therefore
\[
0\neq \overline{\theta_{\mu,t}} = \overline{\pi_\mu(p^{-v}u_j)} \in \overline{R_\mu}.
\]
We argue that the lemma now reduces to the $p$-saturation (see Definition~\ref{def:p-saturation}) of the PBW filtration on $\mathcal{U}$. Indeed, if the filtration is $p$-saturated, then 
\[
u_j\in \mathcal{F}_d\,\mathcal{U}\implies p^{-v}u_j \in \mathcal{F}_d\,\mathcal{U}\implies \overline{\theta_{\mu,t}} \in \mathcal{F}_d\overline{R_\mu},
\]
where $\overline{R_\mu}$ is equipped with the quotient filtration. The lemma now follows from Claim~\ref{claim:pbw-is-saturated}.
\end{proof}

Let $V_{\lambda,s}$ denote a pure weight $s$ representation of $\GSp_{2g,\Qp}$ of highest weight $\lambda \in \mathbb{Z}^g$. Fixing
\[
t = t_s := (1+p^c)^s-1 \in p\Zp,
\]
we have
\[
\mathrm{Ker}\,\eta|_{V_{\lambda,s}} = \mathrm{Ker}\,\theta_{\mu,t}|_{V_{\lambda,s}},
\]
where $\theta_{\mu,t}$ is as in Lemma~\ref{lem:key_lemma}. We are ready to state the main proposition of this segment.
\begin{proposition}\label{prop:key_prop}
Notations as above, one has
\[
\dim_\Qp \mathrm{Ker} \theta_{\mu,t}|_{V_{\mu,s}} \le d_\mu^{-1}\left(\dim_\Fp \mathcal{F}_{\lceil 2D||\lambda||\rceil + d}\overline{R_\mu} - \dim_\Fp \mathcal{F}_{\lceil 2D||\lambda||\rceil}\overline{R_\mu}\right).
\]
\end{proposition}

Before we may carry on with the proof of Proposition~\ref{prop:key_prop}, we shall require a few facts about the map $\mathcal{U}_\Qp \to \End_\Qp V_{\mu,s}$, derived from the Jacobson density theorem.
\begin{theorem}[Jacobson's Density Theorem; see~{\cite[Theorem XVII.3.2]{LangAlgebra}}]\label{thm:jacobson}
Let $A$ be a unital, associative, and distributive ring (not necessarily commutative), and let $M$ be a finitely generated simple $A$-module. Write $D := \End_A M$, then the natural map,
\[
A\to \End_D M,
\]
is surjective.
\end{theorem}

We shall apply the Jacobson density theorem in the following way.
\begin{corollary}\label{cor:jacobson}
Notations as above, the natural map $\mathcal{U}_\Qp \to \End_\Qp V_{\mu,s}$ is surjective.
\end{corollary}
\begin{proof}
Since the underlying base field $\Qp$ has characteristic zero, $V_{\mu,s}$ is a simple $\mathfrak{h}_\Qp$-module, hence also a simple $\mathcal{U}_\Qp$ module. Let $v_\mu\in V_{\mu,s}$ denote a highest weight vector, we argue that a $\mathcal{U}_\Qp$-equivariant endomorphism, $\phi\in \End_{\mathcal{U}_\Qp} V_{\mu,s}$ must preserve the highest weight line $\Qp\cdot v_\mu$. Indeed, such an endomorphism must commute with the Cartan algebra $\mathcal{U}(\mathfrak{t})\subseteq \mathcal{U}$, which acts through the weight character, 
\[
t\cdot v_\mu = \mathrm{ev}_\mu(t)v_\mu,\qquad t\in \mathcal{U}(\mathfrak{t}).
\]
If $\phi(v_{\mu}) \not\in \Qp v_\mu$, then $t\cdot \phi(v_\mu) = \phi(t\cdot v_\mu) = \mathrm{ev}_\mu(t)\cdot \phi(v_\mu)$, which implies that $\phi(v_\mu)\in \Qp v_\mu$, a contradiction.

Since $\phi(v_\mu)\in \Qp v_\mu$, let us denote the constant by $c\in \Qp$ the constant for which $\phi(v_\mu) = c v_\mu$. We claim that $\phi(v) = cv$ for any $v\in V_{\mu,s}$. Indeed, since $V_{\mu,s}$ is a simple $\mathcal{U}_\Qp$-module, $\mathcal{U}_{\Qp}v_\mu = V_{\mu,s}$. Therefore, for any $v\in V_{\mu,s}$ there is some $u\in \mathcal{U}_{\Qp}$, such that $uv_\mu = v$. In particular,
\[
\phi(v) = \phi(uv_\mu) = u\phi(v_\mu) = u\cdot cv_\mu = cv.
\]
Where we used, in the second equality, that $\phi$ is $\mathcal{U}_\Qp$-linear.

We conclude that $\End_{\mathcal{U}_\Qp} V_{\mu,s} = \Qp$. The corollary now follows directly from Jacobson's Density Theorem~\ref{thm:jacobson}
\end{proof}

In fact, the surjectivity statement of Corollary~\ref{cor:jacobson} can be strenghtened.
\begin{claim}\label{claim:jacobson}
The restriction of the surjection $\mathcal{U}_{\Qp}\twoheadrightarrow \End_\Qp V_{\mu,s}$ to $\mathcal{F}_{\lceil 2D||\lambda||\rceil}\mathcal{U}_{\Qp}$ is also surjective.
\end{claim}
\begin{proof}
By the PBW theorem, every element of $\mathcal{U}_\Qp$ may be written as a finite linear combination of (ordered) products of the form 
\[
\sum_{\underline{n}}a_{\underline{n}}X_1^{n_1}\cdot\dots\cdot X_D^{n_D},\qquad \underline{n} = (n_1,\dots,n_D)\in \Z_{\ge 0}^D,
\]
where the $(X_i)_{i=1}^D$ form a basis of $\mathfrak{h}$, regarded as the generators of $\mathcal{U}$ lying in degree $1$, and the $a_{\underline{n}}\in \Qp$ are $p$-adic rational numbers having a common denominator of bounded degree. 

We argue that the exponents $n_i$ can be taken to be at most $2||\lambda||$. Indeed, choose a basis $(X_i)_{i=1}^D$ of $\mathfrak{h}$ consisting of root vectors and Cartan elements. 

\textbf{Root vectors.}
For $X_\alpha$ a root vector, shifting weight via $\nu\mapsto \nu+\alpha$, we observe that $\alpha\in \Z^g$ has at least one coordinate of absolute value at least $1$, and that the weight polytope of $V_{\mu} := \mathrm{Res}_{\Sp_{2g,\Qp}}V_{\mu,s}$ is contained in $[-\mu_1,\mu_1]^g$. For any element $v\in V_{\mu,s}$ of pure weight $\nu\in [-\mu_1,\mu_1]^g$, $X_\alpha^{2\mu_1+1}v = 0$. Since such vectors span $V_{\mu,s}$, the image of $\mathcal{U}_\Qp$ is generated by those $u$ in which root basis vectors occur with exponent at most $2\mu_1 \le 2||\mu||$. 

\textbf{Cartan elements.}
If $X$ is a Cartan element, then it acts on a pure weight-$\nu$ vector $v\in V_{\mu}$ via the weight character, $\mathrm{ev}_\nu(X)$. Now, fix a basis of Cartan elements, $X_1,\dots,X_g \in \mathfrak{t}$, where $X_iv = \nu_iv$ for every pure weight $v$ of weight $\nu$. A basis of pure weight $v$'s diagonalizes the $X_i$'s simultaneously, and their eigenvalues are all bounded by $[-\mu_1,\mu_1]$. Therefore, the minimal polynomial of each $X_i$ divides $P(X) := \prod_{m=-\mu_1}^{\mu_1}(X - m)$. Since $P(X)$ has degree $2\mu_1 + 1\le 2||\lambda|| + 1$, it allows us to express the image in $\End(V_{\mu,s})$ of any monomial in the $X_i$'s as a sum of monomials of degree at most $2||\lambda||$. 

We conclude that the image of $\mathcal{U}_\Qp$ in $\End_\Qp(V_{\mu,s})$ is generated by sums of products of $D$ monomials, each of degree at most $2||\lambda||$, which completes the proof.
\end{proof}

We are now ready to complete the proof of Proposition~\ref{prop:key_prop}.
\begin{proof}[Proof of Proposition~\ref{prop:key_prop}]
Write $\theta$ for the image of $\theta_{\mu,t}$ in $\End_\Qp(V_{\mu,s})$. As a first step, we observe that 
\[
\dim_\Qp \dfrac{\End_\Qp(V_{\mu,s})}{\End_\Qp(V_{\mu,s})\theta} = \dim_\Qp \End_\Qp(V_{\mu,s}) - \dim_\Qp \End_\Qp(V_{\mu,s})\theta = d_\mu^2 - d_\mu\cdot\dim_\Qp\mathrm{Im}\theta,
\]
which, by the first isomorphism theorem, is the same thing as $d_\mu\cdot\dim_\Qp\mathrm{Ker}\theta$. Our statement therefore reduces to showing that
\begin{equation}\label{eq:rk-bound-goal}
\dim_\Qp \dfrac{\End(V_{\mu,s})}{\End_\Qp(V_{\mu,s})\theta} \le \dim_\Fp \mathcal{F}_{\lceil 2D||\lambda||\rceil + d}\overline{R_\mu} - \dim_\Fp \mathcal{F}_{\lceil 2D||\lambda||\rceil}\overline{R_\mu}.  
\end{equation}
Next, we move to an integral setting. Let $\mathcal{V}_{\mu,s}\subseteq V_{\mu,s}$ denote a Chevalley lattice. By this we mean a $\Zp$-submodule of $V_{\mu,s}$ of full rank, stable under $\mathcal{U}$. We note that $\End_\Zp(\mathcal{V}_{\mu,s})$ is a free $\Zp$-module of finite rank, and that $\End_\Zp(\mathcal{V}_{\mu,s})[1/p] = \End_\Qp(V_{\mu,s})$. Set
\[
N_\mu := \mathrm{Im}\left(\mathcal{U}\to \End_\Zp(\mathcal{V}_{\mu,s})\right).
\]
The proof of Claim~\ref{claim:jacobson} replaces monomials of degree $\lceil2||\mu||\rceil + 1$ by smaller degree monomials using monic relations. Therefore, we have a surjection
\[
\mathcal{F}_{\lceil 2D||\mu||\rceil}\mathcal{U}\twoheadrightarrow N_\mu.
\]
Since $\End_\Zp(\mathcal{V}_{\mu,s})$ is a free $\Zp$ module of finite rank, the image $N_\mu \subseteq \End_\Zp(\mathcal{V}_{\mu,s})$ is a submodule of a free module hence also free, and therefore $p$-adically complete. It follows that the map $\mathcal{U}\to N_\mu$ factors through the $p$-completion, $\mathcal{U}\to \widehat{\mathcal{U}}\to N_\mu$. In particular, $\mathcal{U}$ surjects onto $\widehat{\mathcal{U}}$'s image, and this includes the image of $\theta_{\mu,t}$, which we also denote by $\theta$. We obtain a surjection:
\begin{equation}\label{eq:u-surjection}
\mathcal{F}_{\lceil 2D||\mu||\rceil}\mathcal{U}\twoheadrightarrow \dfrac{N_\mu}{N_\mu\theta}.  
\end{equation}
By Corollary~\ref{cor:jacobson}, $N_\mu[1/p] = \End_\Qp (V_{\mu,s})$. Since inverting $p$ is exact, we have
\[
\mathrm{rank}_\Zp\dfrac{N_\mu}{N_\mu\theta} = \dim_\Qp \dfrac{N_\mu[1/p]}{N_\mu[1/p]\cdot \theta} = \dim_\Qp\dfrac{\End_\Qp (V_{\mu,s})}{\End_\Qp (V_{\mu,s})\theta}.
\]
Therefore,~\eqref{eq:rk-bound-goal} becomes
\begin{equation}\label{eq:rk-bound-goal2}
\mathrm{rank}_\Zp\dfrac{N_\mu}{N_\mu\theta} \le \dim_\Fp \mathcal{F}_{\lceil 2D||\lambda||\rceil + d}\overline{R_\mu} - \dim_\Fp \mathcal{F}_{\lceil 2D||\lambda||\rceil}\overline{R_\mu}.
\end{equation}
Since by Proposition~\ref{prop:affinoid_center}, the map $\mathcal{U}\to \End_\Zp(\mathcal{V}_{\mu,s})$ factors through $R_\mu$, combining with~\eqref{eq:u-surjection}, we obtain a surjection:
\begin{equation}\label{eq:r-surjection}
\mathcal{F}_{\lceil 2D||\mu||\rceil}R_\mu\twoheadrightarrow \dfrac{N_\mu}{N_\mu\theta}.
\end{equation}

By Nakayma's lemma,
\[
\mathrm{rank}_\Zp\dfrac{N_\mu}{N_\mu\theta} \le \dim_\Fp \dfrac{N_\mu}{N_\mu\theta + pN_\mu}.
\]
As before, maintain the convention of writing $\overline{R_\mu}, \overline{N}_\mu, \overline{\theta_{\mu,t}}$ and $\overline{\theta}$ for mod-$p$ reductions. Since tensor products are right exact,~\eqref{eq:r-surjection} gives rise to the surjection
\[
s: \mathcal{F}_{\lceil 2D||\mu||\rceil + d}\overline{R_\mu} \twoheadrightarrow \dfrac{\overline{N_\mu}}{\overline{N_\mu}\cdot \overline{\theta}}.
\]
Here $d \ge 0$ is as in Proposition~\ref{prop:key_prop}, an integer with the property that $\theta_{\mu,t} + pR_\mu \in \mathcal{F}_d\overline{R_\mu}$. Technically speaking, the surjection also holds without increasing the PBW filtered piece from $\lceil 2D||\mu||\rceil$ to $\lceil 2D||\mu||\rceil + d$, but this makes the rest of the argument slightly smoother. We observe that the kernel of the surjection above contains $\mathcal{F}_{\lceil 2D||\mu||\rceil + d}\left(\overline{R_\mu}\cdot\overline{\theta}\right)$. Since, by Proposition~\ref{prop:key_prop} $0\neq \overline{\theta} \in \mathcal{F}_d\overline{R_\mu}$,
\[
\mathrm{Ker}(s) \supseteq \left(\mathcal{F}_{\lceil 2D||\mu||\rceil}\overline{R_\mu}\right)\overline{\theta},
\]
and since by Proposition~\ref{prop:r-mu-reduction-domain} $\overline{R_\mu}$ is a domain,
\[
\dim_\Fp \mathrm{Ker}(s) \ge \dim_\Fp \mathcal{F}_{\lceil 2D||\mu||\rceil}\overline{R_\mu}.
\]
\eqref{eq:rk-bound-goal2} now follows from the first isomorphism theorem.
\end{proof}

\subsubsection{Proof of the nontrivial bound}\label{sssec:nontriv_proof}
To complete the proof of Theorem~\ref{thm:equiv-hom-bound}, we seek to bound the $\Fp$-dimensions of the associated gradeds $\dim_\Fp \mathrm{gr}_k \overline{R_\mu}$. Indeed, this is due to the identity
\begin{equation}\label{eq:graded-bound-target}
\dim_\Fp \mathcal{F}_{\lceil 2D||\mu||\rceil + d}\overline{R_\mu} - \dim_\Fp \mathcal{F}_{\lceil 2D||\mu||\rceil}\overline{R_\mu} = \sum_{m = \lceil 2D||\mu||\rceil + 1}^{\lceil 2D||\mu||\rceil + d}\dim_\Fp\mathrm{gr}_m \overline{R_\mu}.
\end{equation}

We shall prove the following bound, using a basic Hilbert polynomial argument.
\begin{proposition}\label{prop:graded_bound}
One has
\[
\dim_\Fp \mathrm{gr}_m\overline{R_\mu} \ll_g m^{D-g-1}.
\]
\end{proposition}
Before we carry out the proof of Proposition~\ref{prop:graded_bound}, we shall require a more concrete presentation of the associated graded algebra $\mathrm{gr}\overline{R}_\mu$. We start out by obtaining a presentation of $\mathrm{gr}R_\mu$.

Since $R_\mu := \mathcal{U}/(c_i-\zeta_\mu(c_i): 1\le i\le g)$ is equipped with the quotient filtration, and taking associated gradeds is exact, we have
\begin{equation}\label{eq:gr_R}
\mathrm{gr} R_\mu \cong \mathrm{gr}\,\mathcal{U}/\mathrm{gr}(c_i-\zeta_\mu(c_i): 1\le i\le g).
\end{equation}
Here $(c_i-\zeta_\mu(c_i): 1\le i\le g)$ denotes the two-sided ideal generated by the elements $c_i-\zeta_\mu(c_i)$. By the notation $\langle c_i-\zeta_\mu(c_i): 1\le i\le g\rangle$, we shall mean the $\mathcal{U}$ subalgebra generated by those same elements.

\begin{proposition}\label{prop:gradeds_presentation}
$\mathrm{gr}\langle c_i-\zeta_\mu(c_i): 1\le i\le g\rangle = \mathrm{gr}\mathcal{Z}$ in $\mathrm{gr}\,\mathcal{U} = S(\mathfrak{h})$.
\end{proposition}
\begin{proof}
Write $A(\mu) := \langle c_i-\zeta_\mu(c_i): 1\le i\le g\rangle$ in short, let $a\in \mathcal{F}_n A(\mu)$, and let $P\in \Zp[x_1,\dots,x_g]$ be a polynomial with the property that $P(c_1-\zeta_\mu(c_1),\dots,c_g-\zeta_\mu(c_g)) = a$. Regard the algebra $\Zp[x_1,\dots,x_g]$ as weight-graded, with $x_i$ having degree $2i$. We claim that 
\[
P\in \Zp[x_1,\dots,x_g]_{\le n}
\]
Indeed, assume otherwise, and let $P^{\mathrm{top}}$ denote the maximal degree homogeneous piece of $P$. Under the assumption that $P\not\in \Zp[x_1,\dots,x_g]_{\le n}$, $P^{\mathrm{top}} \neq 0$ is (weighted-)homogeneous of degree $m > n$, then 
\[
\sigma(i) = \sigma_m(P^{\mathrm{top}}(c_1-\zeta_\mu(c_1),\dots,c_g-\zeta_\mu(c_g))) 
\]
\[
= P^{\mathrm{top}}(\sigma_2(c_1-\zeta_\mu(c_1)),\dots,\sigma_{2g}(c_g-\zeta_\mu(c_g))) = P^{\mathrm{top}}(q_1,\dots,q_g).
\]
Since, by assumption, $i$ has degree $n < m$, $\sigma_m(i) = 0$, which implies that $P^{\mathrm{top}}(q_1,\dots,q_g) = 0$ is a non-trivial algebraic relation between the $q_i$, contradicting their algebraic independence; see~\eqref{eq:invariant-polynomial-generators}. It follows that weighted homogeneous polynomials in the $c_i-\zeta_\mu(c_i)$ of degree at most $n$ generate $\mathcal{F}_nA(\mu)$. In particular, principal symbols of weighted homogeneous polynomials of degree exactly $n$ generate $\mathrm{gr}_nA(\mu)$. Let $Q$ be such a polynomial of degree $n$. Then
\[
\sigma(Q(c_1-\zeta_\mu(c_1),\dots,c_g-\zeta_\mu(c_g))) = Q(q_1,\dots,q_g) \in \mathrm{gr}_n\mathcal{Z}.
\]
It follows that there is a natural inclusion
\[
\mathrm{gr}_nA(\mu) \subseteq \mathrm{gr}_n\mathcal{Z}.
\]
Since $\mathrm{gr}\mathcal{Z} = \langle q_1,\dots,q_g\rangle$, and $q_1,\dots,q_g\in \mathrm{gr}A(\mu)$, we obtain the reverse inclusion, which completes the proof.
\end{proof}
Now, since
\[
(c_i-\zeta_\mu(c_i): 1\le i\le g) = \langle c_i-\zeta_\mu(c_i): 1\le i\le g\rangle\cdot\mathcal{U},
\]
Proposition~\ref{prop:gradeds_presentation} implies that
\[
\mathrm{gr}(c_i-\zeta_\mu(c_i): 1\le i\le g) = \mathrm{gr}\mathcal{Z}\cdot \mathcal{U} = (q_1,\dots,q_g).
\]
Plugging back into~\eqref{eq:gr_R}, we obtain
\begin{equation}\label{eq:graded_rmu}
\mathrm{gr} R_\mu \cong \dfrac{S(\mathfrak{h})}{(q_1,\dots,q_g)}.
\end{equation}
We are now ready to prove Proposition~\ref{prop:graded_bound}.
\begin{proof}[Proof of Proposition~\ref{prop:graded_bound}]
Reducing~\eqref{eq:graded_rmu} mod $p$, we obtain
\[
\overline{\mathrm{gr} R_\mu} = \dfrac{S(\mathfrak{h}_k)}{(\overline{q_1},\dots,\overline{q_g})}.
\]
Since $S(\mathfrak{h}_k) \cong \Fp[x_1,\dots,x_D]$ is a free algebra on $D$-generators, and $\overline{\mathrm{gr} R_\mu}$ is presented as its quotient by an ideal generated by $g$ elements, we have
\begin{equation}\label{eq:krull-lb}
\dim_{\mathrm{Krull}}\overline{\mathrm{gr} R_\mu} \ge D-g.
\end{equation}
To obtain the converse dimension inequality, observe that by Theorem~\ref{thm:aw_main},
\[
\overline{\mathrm{gr} R_\mu}\cong \dfrac{S(\mathfrak{h}_k)}{S(\mathfrak{h}_k)^{\mathbf{H}_k}_+}.
\]
We shall require the following Theorem regarding $\mathcal{N}^*:=\Spec\,\dfrac{S(\mathfrak{h}_k)}{S(\mathfrak{h}_k)^{\mathbf{H}_k}_+}$, which is often referred to as the nilpotent cone in the literature. According to~\cite[\S4.6, Proof of Lemma (c)]{ArdakovWadsleyVerma}, there is a dominant map
\[
\mathbf{H}_{\overline{k}}\times^{B_{\overline{k}}}\mathfrak{b}_{\overline{k}}^\perp\to \mathcal{N}^*_{\overline{k}}.
\]
Here $B_{\overline{k}}$ is a Borel subgroup of $\mathbf{H}_{\overline{k}}$, and $\mathfrak{b}_{\overline{k}}^\perp\subseteq \mathfrak{h}^*_{\overline{k}}$ is the annihilator of $\mathfrak{b}_{\overline{k}}\subseteq \mathfrak{h}_{\overline{k}}$ within the Lie co-algebra $\mathfrak{h}^*_{\overline{k}}$. It follows that
\[
\dim \mathbf{H}_{\overline{k}}\times^{B_{\overline{k}}}\mathfrak{b}_{\overline{k}}^\perp\ge \dim \mathcal{N}^*_{\overline{k}} = \dim \mathcal{N}^*.
\]
The left-hand side has dimension $\dim \mathbf{H}_{\overline{k}}/\mathbf{B}_{\overline{k}} + \dim \mathfrak{b}_{\overline{k}}^\perp = 2\dim\mathbf{H}_{\overline{k}} - 2\dim B_{\overline{k}}$. Note that $\dim \mathbf{H}_{\overline{k}} = D = 2g^2 + g$, while $\dim B_{\overline{k}} = g^2 + g$, hence 
\[
\dim \mathcal{N}^* \le 2g^2 = D-g,
\]
which establishes the converse to~\eqref{eq:krull-lb}. It follows that $\dim_{\mathrm{Krull}} \mathrm{gr}\overline{R_\mu} = D-g$. Recall that a graded algebra $R$ over a field $k'$, such that
\begin{enumerate}
  \item $R$ is generated in degree $1$ by finitely many elements, and
  \item $R_0 = k'$,
\end{enumerate}
has a Hilbert series $P(T) \in \mathbb{Q}(T)$, with the property that for all $n$ large enough,
\[
\dim R_n = P(n)
\]
where $\deg P = \dim_{\mathrm{Krull}}R - 1$; see~\cite[Tag 00P5]{STACKS-PROJECT}. Since $\mathrm{gr}\mathcal{Z}$ is a graded ideal, $\mathrm{gr}\overline{R_\mu} = S(\mathfrak{h}_k)/\mathrm{gr}\overline{\mathcal{Z}}$ is a graded algebra. Its degree $0$ piece is exacly $k = \Fp$, and it is generated in degree $1$ over $k$ by finitely many elements as a quotient of $S(\mathfrak{h}_k) = \Fp[x_1,\dots,x_D]$. Finally, since $\dim_{\mathrm{Krull}} \mathrm{gr}\overline{R_\mu} = D-g$, the Hilbert polynomial of $\mathrm{gr}\overline{R_\mu}$ has degree $D-g-1$, which implies that
\[
\dim_\Fp \mathrm{gr}_m\overline{R_\mu} \ll m^{D-g-1},
\]
and completes the proof.
\end{proof}
We are now ready to complete the proof of the nontrivial bound.
\begin{proof}[Proof of Theorem~\ref{thm:equiv-hom-bound}]
Under the assumptions of Theorem~\ref{thm:equiv-hom-bound}, Proposition~\ref{prop:key_prop} and~\eqref{eq:graded-bound-target} show that
\[
\dim_\Qp \mathrm{Ker}\,\eta|_{V_{\lambda,s}} \le d_\lambda^{-1}\sum_{m = \lceil 2D||\lambda||\rceil + 1}^{\lceil 2D||\lambda||\rceil + d}\dim_\Fp \mathrm{gr}_{m}\overline{R_\lambda}.
\]
Plugging back the bound $\dim_\Fp \mathrm{gr}_m\overline{R_\lambda} \ll_g m^{D-g-1}$ of Proposition~\ref{prop:graded_bound}, we obtain
\[
\dim_\Qp \mathrm{Ker}\,\eta|_{V_{\lambda,s}} \ll d_\lambda^{-1}\sum_{m = \lceil 2D||\lambda||\rceil + 1}^{\lceil 2D||\lambda||\rceil + d}m^{D-g-1} \ll d_{\lambda}^{-1}(1 + ||\lambda||)^{D-g-1},
\]
which completes the proof of Theorem~\ref{thm:equiv-hom-bound}.
\end{proof}

\subsection{Proof of Theorem~\ref{thm:h2-bound}}\label{ssec:h2-bound}
We may now complete the proof of Theorem~\ref{thm:h2-bound}. Recall that this theorem states that for all but finitely many bad weights $E\subseteq \Z$, one has for any irreducible representation $V_{\lambda,s}$ with $s\not\in E$,
\[
\dim_\Qp H^2(G_{\mathbb{Q},\Sigma},V_{\lambda,s}) \ll \min\left\{d_\lambda, \dfrac{(1+||\lambda||)^{D-g-1}}{d_\lambda}\right\},
\]
where the bound $d_\lambda$ holds for all $s$.

\begin{proof}[Proof of Theorem~\ref{thm:h2-bound}]
Recall that by Corollary~\ref{cor:galois-iwasawa}, we have
\[
\dim_\Qp H^2(G_{\mathbb{Q},\Sigma},V_{\lambda,s}) \le \dim_\Qp \Hom_\Gamma(X,V_{\lambda,2-s}),
\]
and that by~\eqref{eq:hom-to-ker}, one has
\[
\dim_\Qp \Hom_\Gamma(X,V_{\lambda,2-s})\le \sum_{i=1}^q\dim_\Qp\mathrm{Ker}\,\delta_i|_{V_{\lambda,2-s}},
\]
where the $(\delta_i)_{i=1}^q\in \Lambda$ are a finite set of annihilators of the Iwasawa module $X$. The trivial bound
\[
\dim_\Qp\mathrm{Ker}\,\delta_i|_{V_{\lambda,2-s}} \ll d_\lambda
\]
holds for all $i$ and all $s$, as was already established in Proposition~\ref{prop:triv_bound}. The last reduction step of \S\ref{ssec:red2} (see Corollary~\ref{cor:red-2}) shows that outside of a finite set of weights $E_i = E_i(\delta_i)\subseteq \Z$, the operators $\delta_i$ may be replaced by operators $0\neq \eta_i\in \Omega$, with 
\[
0\neq \eta_i(t)\in \Zp[[H_c]],\qquad t\in p\Zp,
\]
such that for all $s\not\in E_i$,
\[
\mathrm{Ker}\,\delta_i|_{V_{\lambda,2-s}} = \mathrm{Ker}\,\eta_i|_{V_{\lambda,2-s}}.
\]
Theorem~\ref{thm:equiv-hom-bound} now applies, and provides the nontrivial bound
\[
\dim_\Qp \mathrm{Ker}\,\eta_i|_{V_{\lambda,2-s}} \ll \dfrac{(1+ ||\lambda||)^{D-g-1}}{d_\lambda}.
\]
Set $E := \cup_{i=1}^q E_i$, then each $E_i$ is finite, hence so is $E$. It follows that for all $s\not\in E$, both the trivial and nontrivial bounds hold simultaneously for all $1\le i\le q$, and we get
\[
\dim_\Qp H^2(G_{\mathbb{Q},\Sigma},V_{\lambda,s}) \le \sum_{i=1}^q \dim_\Qp \mathrm{Ker}\,\eta_i|_{V_{\lambda,2-s}} \ll \min\left\{d_\lambda,\dfrac{(1+ ||\lambda||)^{D-g-1}}{d_\lambda}\right\},
\]
which completes the proof.
\end{proof}

\section{The Chabauty--Kim dimension inequality via character sum cancellations}\label{sec:ck-dim-ineq}
This chapter is structured as follows:
\begin{enumerate}
  \item In~\S\ref{ssec:character-estimates}, we bound various representation theoretic quantities associated with the graded pieces $L_n$ of the associated graded Lie algebra $L$ of $\mathfrak{n} := U^{\et}$, and complete the proof of Theorem~\ref{thm:power_saving}.
  \item In~\S\ref{ssec:kim-conj-proof}, we complete the proof of Kim's conjecture (see Corollary~\ref{cor:kim-conj}).
\end{enumerate}

\subsection{Auxiliary character estimates}\label{ssec:character-estimates}
As in the introduction, write $L$ for the associated graded Lie algebra of $\mathfrak{n}$, the pro-nilpotent Lie algebra of the pro-unipotent \'etale fundamental group of our fixed smooth projective, geometrically connected, hyperbolic curve $C/\mathbb{Q}$. In this case $L$ would be a quadratically presented Lie algebra, given explicitly as the quotient of the free Lie algebra on $V := H_1^\et(\overline{C};\Qp)$ by the single relation coming from the image of the cup co-product, i.e.~the map $\cup^\vee$ dual to the Weil pairing:
\[
\cup: H^1_\et(\overline{C};\Qp)\otimes H^1_\et(\overline{C};\Qp)\to \Qp(-1);
\]
see~\cite[Fact~3.14]{CORWIN_ZEHAVI}. Our goal in this subsection is to compute asymptotics and bounds of certain representation theoretic invariants associated with $L$. For each $n\ge 1$, denote by $L_n$ the $n$'th graded piece of $L$, and decompose $L$ into a direct sum of irreducible representations of $\Sp_{2g,\Qp}$, where $g\ge 2$ is the genus of $C/\mathbb{Q}$. For each highest weight vector $\lambda \in \Sp_{2g,\Qp}$, write $V_\lambda$ for the associated simple representation, and denote by $m_n(\lambda)$ the multiplicity of $V_\lambda$ in $L_n$, such that we have a decomposition of the form:
\begin{equation}\label{eq:ln-decomp}
L_n = \bigoplus_\lambda V_\lambda^{\oplus m_n(\lambda)}.
\end{equation}
We shall be interested in the asymptotics of $b_n = \dim_\Qp L_n$, and in bounding $m_n(\lambda)$ and the second-moment sum, $\sum_\lambda m_n(\lambda)d_\lambda ||\lambda||^2$.

Our starting point is the observation that, since $\Sp_{2g,\Qp}$ is a totally split connected reductive group, its irreducible representations are absolutely irreducible; see~\cite[20.6]{MilneReductiveGroups}. The decomposition of $L_n$ into irreducibles in the sense of~\eqref{eq:ln-decomp} persists through scalar extension of the form $\Qp\to \mathbb{C}$. To be precise, write $V_{\lambda,\mathbb{C}} := V_\lambda\otimes \mathbb{C}$ for the associated rational representation of $\Sp_{2g,\mathbb{C}}$, then the decomposition of $L_{n,\mathbb{C}}$ into irreducibles has the form
\[
L_{n,\mathbb{C}} = \bigoplus_\lambda V_{\lambda,\mathbb{C}}^{\oplus m_n(\lambda)}.
\]
In particular, the $V_{\lambda,\mathbb{C}}$ remain irreducible representations of $\Sp_{2g,\mathbb{C}}$, and appear with the same multiplicity in $L_{n,\mathbb{C}}$ as $V_{\lambda}$ appears in $L_n$. We seek to restrict $L_{n,\mathbb{C}}$ to the compact real form of $\Sp_{2g,\mathbb{C}}$.

\begin{definition}[The unitary symplectic group; see~{\cite[\S7.2, p.~99, Exercise~7.4]{FULTON_HARRIS}}]
Write $\mathrm{USp}(2g) = \Sp_{2g}(\mathbb{C})\cap \mathrm{U}(2g)$ for the unitary symplectic group, cut out from $\Sp_{2g}(\mathbb{C})$ by its subgroup of unitary elements. Explicitly,
\[
\mathrm{USp}(2g) := \left\{M\in \Sp_{2g}(\mathbb{C}): M^TJM = J, \overline{M}^tM = I\right\},\qquad J =
\begin{pmatrix}
0& I_{g}\\
-I_{g}& 0
\end{pmatrix},
\]
where for a matrix $M$ in $M_{2g}(\mathbb{C})$, $\overline{M}$ denotes its complex conjugate.
\end{definition}
The advantage of restricting our representations from $\Sp_{2g,\mathbb{C}}$ to $\mathrm{USp}(2g)$ is that this functor preserves the finite-dimensional algebraic representation theory of $\Sp_{2g,\mathbb{C}}$, while allowing us to employ the character theory of compact Lie groups.

\begin{proposition}\label{prop:equiv-compact}
Restriction from $\Sp_{2g,\mathbb{C}}$ to $\mathrm{USp}(2g)$ preserves irreducible representations.
\end{proposition}

\begin{proof}
Write $\mathfrak{k} := \Lie(\mathrm{USp}(2g))$, and $\mathfrak{h}_\mathbb{C} := \Lie(\Sp_{2g,\mathbb{C}})$. Since $\mathrm{USp}(2g)$ is a compact real form of $\Sp_{2g,\mathbb{C}}$, we have
\[
\mathfrak{k}\otimes \mathbb{C} = \mathfrak{h}_\mathbb{C};
\]
see~\cite[Chapter~VIII, \S7, Theorem~12(a)]{SerreComplexSemisimpleLieAlgebras}. Now, let $W$ denote a finite dimensional irreducible representation of $\Sp_{2g,\mathbb{C}}$, and assume its restriction to $\mathrm{USp}(2g)$ splits into a direct sum of $W = W_1 \oplus W_2$, with the $W_i$ nonzero, stable $\mathrm{USp}(2g)$-modules. Differentiating the $\mathrm{USp}(2g)$-action, we find that each $W_i$ is $\mathfrak{k}$-stable, but then, by complex linearity, it is also $\mathfrak{h}_\mathbb{C}$-stable. Since $\Sp_{2g,\mathbb{C}}$ is connected, the decomposition $W = W_1\oplus W_2$ is also $\Sp_{2g,\mathbb{C}}$-stable, contradicting the $\Sp_{2g,\mathbb{C}}$-irreducibility of $W$.
\end{proof}

Write $\mathcal{L}_n$ and $\mathcal{V}_\lambda$ for the restrictions of $L_{n,\mathbb{C}}$ and $V_{\lambda,\mathbb{C}}$ to $\mathrm{USp}(2g)$. In view of Proposition~\ref{prop:equiv-compact}, we obtain an irreducible decomposition of $\mathcal{L}_n$, analogous to~\eqref{eq:ln-decomp}:
\begin{equation}\label{eq:ln-decomp-compact}
\mathcal{L}_n = \bigoplus_{\lambda}\mathcal{V}_{\lambda}^{\oplus m_n(\lambda)}.
\end{equation}

\begin{notation}
Where there is no cause of confusion, we shall write $L_n$ and $V_\lambda$ for the restrictions of the complexified representations to the unitary symplectic group.
\end{notation}

For a $\mathrm{USp}(2g)$-representation $W$, write $\chi_W: \mathrm{USp}(2g)\to \mathbb{C}$ for its character. Recall that, by definition, $\chi_W(g) := \mathrm{tr}(g|_{W})$. Let us abbreviate $\chi_n: \mathrm{USp}(2g)\to \mathbb{C}$ for the character of $\mathcal{L}_n$. Our first goal will be to give an asymptotic estimate for $b_n$. 

\begin{notation}
For functions $f,g:\mathbb{N}\to \mathbb{N}_{>0}$, the notation $f\asymp g$, means that 
\[
\lim_{n\to \infty}\dfrac{f(n)}{g(n)} = 1.
\]
\end{notation}

\begin{proposition}\label{prop:bn-growth}
With
\[
\alpha := g + \sqrt{g^2 - 1},
\]
one has, for all $n\ge 1$,
\[
b_n \asymp \dfrac{\alpha^n}{n}.
\]
\end{proposition}
\begin{remark}
Observe that under our assumption that $C$ is projective and hyperbolic, $g\ge 2$ necessarily. The asymptotics of Proposition~\ref{prop:bn-growth} do not make sense for an elliptic curve, since in this case the (geometric \'etale) fundamental group is abelian, hence its unipotent completion has no associated graded pieces beyond the abelianization. In this case $\alpha^n = 1$ for every $n$, while $b_n = 0$ for every $n\ge 2$.
\end{remark}

Throughout this subsection, we shall repeatedly employ the following presentation of $\chi_n$ due to Filip, expressed in terms of the character of the standard representation, denoted $\chi_V$.
\begin{theorem}[{\cite[Proposition 4.4, Eq.~(4.8)]{FilipSurfaceLie}}]\label{thm:filip_lie}
For all $u\in \mathrm{USp}(2g)$, one has
\[
\chi_n(u) = \dfrac{1}{n}\sum_{d|n}\mu(\dfrac{n}{d})P_d(\chi_V(u^{n/d})),
\]
where $\mu: \mathbb{N}\to \{-1,0,1\}$ is the M\"obius function, and $P_d: \mathbb{R}\to \mathbb{R}$ is the function given explicitly by
\[
P_d(x) = \left(\dfrac{x + \sqrt{x^2 - 4}}{2}\right)^d + \left(\dfrac{x - \sqrt{x^2 - 4}}{2}\right)^d.
\]
\end{theorem}

\begin{remark}\label{rem:domain}
Observe that it indeed suffices to define $P_d$ for real arguments. Indeed, since $u$ is symplectic and unitary, so is $u^{n/d}$, and therefore its eigenvalues come in conjugate pairs $(\gamma_i,\overline{\gamma_i})$. Therefore, $\chi_V(u^{n/d}) = \mathrm{tr}(u^{n/d}) = \sum_{i=1}^g 2\mathrm{Re}(\gamma_i) \in \mathbb{R}$. In fact, $\chi_V(u^{n/d})$ takes values strictly in $[-2g,2g]$.
\end{remark}

\begin{proof}[Proof of Theorem~\ref{prop:bn-growth}]
Writing $I := I_{2g}\in \mathrm{USp}(2g)$, we have
\[
\chi_n(I) = b_n.
\]
Plugging into Theorem~\ref{thm:filip_lie}, we get
\begin{equation}\label{eq:bn}
\chi_n(I) = \dfrac{1}{n}\sum_{d|n}\mu(\dfrac{n}{d})P_d(2g).
\end{equation}
By definition of $P_d$, we have
\[
P_d(2g) = \alpha^d + (g - \sqrt{g^2 - 1})^d,
\]
which gives rise to the trivial bound $|P_d(2g)|\le 2\alpha^d$, valid for every $d\ge 1$. As for $d=n$, we shall apply a slightly more careful analysis. Since
\[
0\le g - \sqrt{g^2 - 1} = (g - \sqrt{g^2 - 1})\cdot \dfrac{g + \sqrt{g^2 - 1}}{g + \sqrt{g^2 - 1}} = \dfrac{1}{g + \sqrt{g^2 - 1}} \le \dfrac{1}{g},
\]
we have
\begin{equation}\label{eq:p2g_alpha_bd}
|P_n(2g) - \alpha^n| \le \dfrac{1}{g^n} = o(1).
\end{equation}
Plugging back to~\eqref{eq:bn}, we obtain
\[
|\chi_n(I) - \dfrac{\alpha^n}{n}| \le o(\dfrac{1}{n}) + \dfrac{2}{n}\sum_{d|n, d\neq n}\alpha^d \ll \dfrac{1}{n}\sum_{d\le n/2}\alpha^d \ll \dfrac{\alpha^{n/2}}{n} = o(\dfrac{\alpha^n}{n}).
\]
Therefore,
\[
\chi_n(I) \asymp \dfrac{\alpha^n}{n},
\]
which completes the proof.
\end{proof}

Our next goal is to bound $m_n(\lambda)$. We shall prove
\begin{proposition}\label{prop:mn-bound}
For all $n\ge 1$ and for any highest weight vector $\lambda\in \mathbb{Z}^g$, one has
\[
m_n(\lambda) \ll d_\lambda b_n n^{-D/2}.
\]
\end{proposition}

\begin{proof}
By additivity of characters,~\eqref{eq:ln-decomp-compact} implies the identity
\[
\chi_n = \sum_{\lambda}m_n(\lambda)\chi_\lambda,
\]
where $\chi_\lambda := \chi_{V_\lambda}$. Now, fix a Haar probability measure $dk$ on $\mathrm{USp}(2g)$. By probability measure, we mean that
\[
\int_{\mathrm{USp}(2g)}dk = 1.
\]
Then, by orthogonality of characters for non-isomorphic $\mathrm{USp}(2g)$-representations, we have
\[
\int_{\mathrm{USp}(2g)}\chi_n(k)\overline{\chi_{\lambda}(k)}dk = m_n(\lambda).
\]
Weyl's unitary trick provides the bound $|\chi_{\lambda}(k)|\le d_\lambda$ for every $k\in \mathrm{USp}(2g)$, and therefore
\begin{equation}\label{eq:mn_1}
m_n(\lambda) \le d_\lambda\int_{\mathrm{USp}(2g)}|\chi_n(k)|dk.
\end{equation}
Our goal is therefore to bound $\int_{\mathrm{USp}(2g)}|\chi_n(k)|dk$. By Filip's formula, we have
\[
\chi_n(k) = \dfrac{1}{n}\sum_{d|n}\mu(n/d)P_d(\chi_V(k^{n/d})).
\]
For $d < n$, we bound $P_d(\chi_V(k^{n/d}))$ uniformly over $k\in \mathrm{USp}(2g)$ as follows. By definition,
\[
P_d(x) = \left(\dfrac{x + \sqrt{x^2 - 4}}{2}\right)^d + \left(\dfrac{x - \sqrt{x^2 - 4}}{2}\right)^d.
\]
A standard calculation shows that 
\begin{equation}\label{eq:small_x}
|x|\le 2\implies |P_d(x)| \le 2.
\end{equation}
While, for $2 < |x| \le 2g$,
\[
|P_d(x)| \le 2\cdot \left(\dfrac{|x| + \sqrt{x^2 - 4}}{2}\right)^d \le 2P_d(x)\le 2P_d(2g) \ll \alpha^d,
\]
where the first inequality follows from the fact that $\left(\dfrac{|x| + \sqrt{x^2 - 4}}{2}\right)^d$ bounds both summands comprising $|P_d(x)|$, the second inequality is due to the fact that both summands have the same sign in the range $2\le |x|\le 2g$ one of which is equal to $\left(\dfrac{|x| + \sqrt{x^2 - 4}}{2}\right)^d$ in absolute value, the third inequality follows from the monotonicity of $|P_d(x)|$ in that range, and the last inequality from~\eqref{eq:p2g_alpha_bd}. By Remark~\ref{rem:domain}, $\chi_V(k^{n/d}) \in [-2g,2g]$, and therefore
\[
|P_d(\chi_V(k^{n/d}))| \ll \alpha^d,\qquad k\in \mathrm{USp}(2g).
\]
Plugging back into our expression for $\chi_n(k)$, we obtain the bound:
\[
|\chi_n(k)| \le \dfrac{1}{n}|P_n(\chi_V(k))| + \dfrac{1}{n}\sum_{d|n,d\neq n}|P_d(\chi_V(k^{n/d}))|
\]
\[
\ll \dfrac{1}{n}|P_n(\chi_V(k))| + \dfrac{1}{n}\sum_{d\le n/2}\alpha^d \ll \dfrac{1}{n}|P_n(\chi_V(k))| + \dfrac{\alpha^{n/2}}{n},
\]
and~\eqref{eq:mn_1} now becomes
\begin{equation}\label{eq:mn_2}
m_n(\lambda) \ll \dfrac{d_\lambda}{n}\int_{\mathrm{USp}(2g)}|P_n(\chi_V(k))|dk + \dfrac{d_\lambda\alpha^{n/2}}{n}.
\end{equation}
Our next goal is to bound 
\[
\int_{\mathrm{USp}(2g)}|P_n(\chi_V(k))|dk.
\]
Let $\mathscr{S}\subseteq \mathrm{USp}(2g)$ denote the subset consisting of those $k$ for which $|\chi_V(k)| < 2$. For such $k$'s, \eqref{eq:small_x} implies that $|P_n(\chi_V(k))|\le 2$, hence, the contribution of $\dfrac{1}{n}\int_{\mathscr{S}}|P_n(\chi_V(k))|dk$ to $\dfrac{1}{n}\int_{\mathrm{USp}(2g)}|P_n(\chi_V(k))|dk$ is bounded by a constant multiple of $1/n$. A standard computation shows that for all $y\in \mathbb{R}$,
\begin{equation}\label{eq:pn-identity}
P_n\left(\pm(e^y + e^{-y})\right) = (\pm 1)^n\cdot\left(e^{ny} + e^{-ny}\right).
\end{equation}
What is special about the complement of $\mathscr{S}$, is that since for all $k\in \mathrm{USp}(2g)\setminus \mathscr{S}$ we have $|\chi_V(k)| \ge 2$, the equation
\[
e^y + e^{-y} = 2\cosh y = |\chi_V(k)|
\]
has a solution. Namely, it is $y = \arccosh \dfrac{|\chi_V(k)|}{2}$.~\eqref{eq:pn-identity} then gives for all such $k$,
\[
|P_n(\chi_V(k))| = e^{n\arccosh (|\chi_V(k)|/2)} + e^{-n\arccosh (|\chi_V(k)|/2)}.
\]
It follows that in $\mathrm{USp}(2g)\setminus \mathscr{S}$,
\[
\begin{aligned}
|P_n(\chi_V(k))|
&\le 2e^{n\arccosh (|\chi_V(k)|/2)}\\
&= 2e^{n\arccosh g}e^{n(\arccosh (|\chi_V(k)|/2) - \arccosh g)}\\
&\le 2(e^{n\arccosh g} + e^{-n\arccosh g})
   e^{n(\arccosh (|\chi_V(k)|/2) - \arccosh g)}\\
&= 2P_n(\chi_V(I))e^{n(\arccosh (|\chi_V(k)|/2) - \arccosh g)}\\
&\ll \alpha^ne^{n(\arccosh (|\chi_V(k)|/2) - \arccosh g)},
\end{aligned}
\]
where the last inequality used~\eqref{eq:p2g_alpha_bd}.

Write $\delta := 2g - |\chi_V(k)|$, then
\[
\arccosh (g - \delta/2) - \arccosh g = -\int_{g - \delta/2}^g \dfrac{dt}{\sqrt{t^2 - 1}} \le - \dfrac{\delta}{2\sqrt{g^2 - 1}}.
\]
Plugging back, we have
\begin{equation}\label{eq:pn_bound}
|P_n(\chi_V(k))| \ll \alpha^n\exp\left(-\dfrac{n(2g - |\chi_V(k)|)}{2\sqrt{g^2 - 1}}\right).
\end{equation}
As above, write $\mathfrak{k} := \Lie(\mathrm{USp}(2g))$, and consider the Hilbert--Schmidt inner product on $\mathfrak{k}$,
\[
\langle X,Y\rangle_{\mathrm{HS}} := \mathrm{Re}\left(\mathrm{tr}(X^*Y)\right) = -\mathrm{tr}(XY),
\]
where the last equality is since the elements of $\mathfrak{k}$ are skew hermitian, and therefore $X^* = -X$.

For each $X\in \mathfrak{k}$, since $X$ is skew hermitian and a logarithm of a symplectic matrix, its eigenvalues consist of purely imaginary pairs, $\pm \theta_1 i, \dots, \pm \theta_g i$, where the $\theta_j\in \mathbb{R}$. Consequently,
\[
\mathrm{tr}(\exp(X)) = 2\sum_{i= 1}^g\cos \theta_i,\qquad ||X||_{\mathrm{HS}}^2 = 2\sum_{i=1}^g \theta_i^2.
\]
For $0 < \theta  < 1$, we bound $1 - \cos \theta \gg \theta^2$ as follows. Since $\sin(\theta/2)$ is concave in this range, we have
\[
\sin(\theta/2) \ge \theta\sin (1/2).
\]
Therefore, 
\[
1 - \cos \theta = 2\sin^2(\theta/2) \ge 2\theta^2\sin^2(1/2) = (1-\cos 1)\cdot \theta^2 \gg \theta^2.
\]

Hence, for $X$ in a small ball around the origin, one has
\begin{equation}\label{eq:hs}
2g - \mathrm{tr}(\exp(X)) \gg 2\sum_{i=1}^g \theta_i^2 = ||X||_{\mathrm{HS}}^2.  
\end{equation}
Let $B_\epsilon\subseteq \mathfrak{k}$ denote the open $L^2$-ball of radius $\epsilon$ around the origin, so that $\exp: B_\epsilon \to \mathrm{USp}(2g)$ is a diffeomorphism onto $B_\epsilon$'s image. For $k\in \exp(B_{\epsilon})$, write $\log k \in \mathfrak{k}$ for $k$'s pre-image. Then, by~\eqref{eq:hs}
\[
\chi_V(I) - \chi_V(k) = 2g - \mathrm{tr}(k) \gg ||\log k||_{\mathrm{HS}}^2.
\]
Since multiplication by $-I$ is an isometry between $\exp(B_{\epsilon})$ and $-\exp(B_{\epsilon})$, which acts on $\chi_V(I) - \chi_V(k)$ by interchanging signs, setting
\[
C_\epsilon(I) := \exp(B_{\epsilon}),\qquad C_{\epsilon}(-I) :=- C_\epsilon(I),\qquad \text{and }C_\epsilon(\pm I) = C_\epsilon(I)\cup C_\epsilon(-I),
\]
we find that there is a constant $c_1\ge 0$, such that for every $k\in C_{\epsilon}(\pm I)$, one has
\[
2g - |\chi_V(k)| \ge c_1||\log k||_{\mathrm{HS}}^2.
\]
Plugging back into~\eqref{eq:pn_bound}, we get for every $k\in C_{\epsilon}(\pm I)$,
\[
|P_n(\chi_V(k))| \ll \alpha^n\exp\left(-\frac{nc_1||\log k||_{\mathrm{HS}}^2}{2\sqrt{g^2 - 1}}\right).
\]
As for $k\not\in C_{\epsilon}(\pm I) \cup \mathscr{S}$, write $\mathcal{D} := \mathrm{USp}(2g)\setminus \left(C_{\epsilon}(\pm I) \cup \mathscr{S}\right)$. Since $C_{\epsilon}(\pm I) \cup \mathscr{S}$ is open, its complement $\mathcal{D}$ in $\mathrm{USp}(2g)$ is a closed subset of a compact set, hence compact as well. Therefore, the non-negative continuous function $2g - |\chi_V(k)|\ge 0$ attains a minimum on this subset. Since $C_{\epsilon}(\pm I)$ contains open neighborhoods around the only points where $2g - |\chi_V(k)|$ vanishes, namely $\pm I$, $2g - |\chi_V(k)|$ is bounded below from $0$ in this region. That is, there exists some absolute constant $c_2 > 0$, independent of $n$, such that for every $k\in \mathcal{D}$, one has
\[
2g - |\chi_V(k)| \ge c_2 > 0.
\]
Plugging back to~\eqref{eq:pn_bound}, we bound 
\[
|P_n(\chi_V(k))| \ll \alpha^n\exp\left(-\frac{c_2n}{2\sqrt{g^2 - 1}}\right) = \alpha^ne^{-c_2'n},
\]
uniformly in $\mathcal{D}$, where $c_2' := c_2/2\sqrt{g^2 - 1} > 0$. In particular, the contribution of the domain $\mathcal{D}$ to our bound on $\dfrac{1}{n}\int_{\mathrm{USp}(2g)}|P_n(\chi_V(k))|dk$ is bounded above by
\[
\dfrac{1}{n}\int_{\mathcal{D}}|P_n(\chi_V(k))|dk \ll b_n\cdot e^{-c_2'n}.
\]
Finally, we are ready to deal with the main contribution to $\dfrac{1}{n}\int_{\mathrm{USp}(2g)}|P_n(\chi_V(k))|dk$, which occurs on the domain $C_{\epsilon}(\pm I)$. Our bound on $|P_n(\chi_V)|$ is invariant under the isometry $k\mapsto -k$, which allows us to bound
\[
\dfrac{1}{n}\int_{C_{\epsilon}(\pm I)}|P_n(\chi_V(k))|dk \ll b_n \int_{C_\epsilon(I)}\exp\left(-\frac{nc_1||\log k||_{\mathrm{HS}}^2}{2\sqrt{g^2 - 1}}\right)dk.
\]
Write $c_1' := c_1/(2\sqrt{g^2 - 1})$, our goal becomes to bound 
\[
\int_{C_\epsilon(I)}\exp(-nc_1'||\log k||_{\mathrm{HS}}^2)dk.
\]
Consider the change of variables $k = \exp(X)$, and parameterize $B_\epsilon$ as follows. Let $Y_1,\dots,Y_D$ denote an orthonormal basis of $\mathfrak{k}$ for the Hilbert--Schmidt inner product. For $X\in \mathfrak{k}$, write $X = \sum_{i=1}^D x_iY_i$ with the $x_i \in \mathbb{R}$. In these coordinates, $B_{\epsilon}$ corresponds to the ball of radius $\epsilon$ around the origin in $\mathbb{R}^D$, equipped with the standard Euclidean metric, $d\underline{x} = dx_1\cdot\dots\cdot dx_D$. Let $J(\underline{x}) = |\det((d\exp)_{\underline{x}})|$ denote the rescaling factor of this reparameterization, such that
\[
\int_{C_\epsilon(I)}\exp(-nc_1'||\log k||_{\mathrm{HS}}^2)dk = \int_{\underline{x}\in \mathbb{R}^D: ||\underline{x}|| < \epsilon}e^{-nc_1'\sum_{i=1}^D x_i^2}J(\underline{x})d\underline{x}.
\]
The rescaling factor $J(\underline{x})$ is continuous, hence bounded on any compact subset. Possibly shrinking $\epsilon > 0$ in hindsight, we may assume $J(\underline{x})$ extends continuously to the boundary of our radius $\epsilon$ ball, and that it is bounded in the domain $\{||\underline{x}|| < \epsilon\}$ by some absolute scalar independent of $n$, in the sense that
\[
\int_{\underline{x}\in \mathbb{R}^D: ||\underline{x}|| < \epsilon}e^{-nc_1'\sum_{i=1}^D x_i^2}J(\underline{x})d\underline{x} \ll 
\int_{\underline{x}\in \mathbb{R}^D: ||\underline{x}|| < \epsilon}e^{-nc_1'\sum_{i=1}^D x_i^2}d\underline{x}.
\]
Since the integrand is non-negative on $\mathbb{R}^D$, we may extend integration throughout, to obtain the bound
\[
\int_{\underline{x}\in \mathbb{R}^D: ||\underline{x}|| < \epsilon}e^{-nc_1'\sum_{i=1}^D x_i^2}d\underline{x}\le 
\int_{\underline{x}\in \mathbb{R}^D}e^{-nc_1'\sum_{i=1}^D x_i^2}d\underline{x}.
\]
In order to see the precise asymptotic dependence on $n$, we change variables once more. Setting $y_i = \sqrt{nc_1'}x_i$, so that $d\underline{y} = (nc_1')^{D/2}d\underline{x}$, we obtain the multivariate Gaussian integral
\[
\int_{\underline{x}\in \mathbb{R}^D}e^{-nc_1'\sum_{i=1}^D x_i^2}d\underline{x} = (nc_1')^{-D/2}\int_{\underline{y}\in \mathbb{R}^D}e^{-\sum_{i=1}^D y_i^2}d\underline{y} \ll n^{-D/2}.
\]

Our bound on $\dfrac{1}{n}\int_{\mathrm{USp}(2g)}|P_n(\chi_V(k))|dk$ becomes,
\[
\dfrac{1}{n}\int_{\mathrm{USp}(2g)}|P_n(\chi_V(k))|dk \ll \dfrac{1}{n} + b_n\cdot e^{-c'n} + b_n\cdot n^{-D/2} \ll b_n\cdot n^{-D/2}.
\]
Finally, plugging back into~\eqref{eq:mn_2}, our bound on $m_n(\lambda)$ becomes
\[
m_n(\lambda) \ll \dfrac{d_\lambda}{n}\int_{\mathrm{USp}(2g)}|P_n(\chi_V(k))|dk + \dfrac{d_\lambda\alpha^{n/2}}{n}\ll d_\lambda b_n n^{-D/2}.
\]
\end{proof}

Our final task in this subsection is to obtain our second-moment bound, i.e.~to bound the sum $\sum_{\lambda}m_n(\lambda)d_\lambda ||\lambda||^2$.

\begin{definition}[see~{\cite[\S6.3, Proposition~6.15]{KirillovLieGroups}}]
Let $(X_i)_{i=1}^D\in \mathfrak{k} = \Lie(\mathrm{USp}(2g))$ denote an orthonormal basis of $\mathfrak{k}$ with respect to the Hilbert--Schmidt inner product $\langle X,Y\rangle_{\mathrm{HS}} := \mathrm{tr}(-XY)$. The Casimir $\mathfrak{C} = \mathfrak{C}_{\mathrm{HS}}$ is defined as
\[
\mathfrak{C} := -\sum_{i=1}^D X_i^2 \in \mathcal{U}_\mathbb{C},
\]
where $\mathcal{U}_{\mathbb{C}} := \mathcal{U}(\mathfrak{h}_\mathbb{C}) = \mathcal{U}(\mathfrak{k}_{\mathbb{C}})$ is the universal enveloping algebra of the complexification of $\mathfrak{k}$, i.e.~of the complex symplectic Lie algebra $\mathfrak{h}_{\mathbb{C}}$.
\end{definition}

\begin{remark}[see~{\cite[\S6.3, Proposition~6.15]{KirillovLieGroups}}]\label{rem:casimir}
The Casimir element is independent of the choice of orthonormal basis, and belongs to the center $\mathcal{Z}(\mathcal{U}(\mathfrak{k}))$.
\end{remark}

The Lie algebra $\mathfrak{k}$ acts naturally on the ring of smooth functions on $\mathrm{USp}(2g)\to \mathbb{C}$ as follows. Let $X\in \mathfrak{k}$, and let $f: \mathrm{USp}(2g)\to \mathbb{C}$ be a smooth function. By $(Xf): \mathrm{USp}(2g)\to \mathbb{C}$ we mean the function given explicitly via:
\[
(Xf)(k) := \dfrac{d}{dt}f(k\exp(tX))|_{t=0}.
\]
The action of $\mathfrak{k}$ naturally extends to an action of the symplectic Lie algebra $\mathfrak{h}_\mathbb{C}$ by linearity, and to an action of its tensor algebra through composition. One verifies that this action descends to the universal enveloping algebra by observing that $(XY - YX)f = [X,Y]f$. We therefore view the Casimir as acting on the ring of smooth functions on $\mathrm{USp}(2g)$ through the action of the enveloping algebra.

\begin{claim}\label{claim:clambda}
For $\lambda$ a highest weight corresponding to an irreducible representation $V_\lambda$ of $\mathrm{USp}(2g)$, write $\rho_\lambda: \mathrm{USp}(2g)\to \GL(V_\lambda)$ for the representation homomorphism, and $d\rho_\lambda: \mathcal{U}_\mathbb{C}\to \End(V_\lambda)$ for the extension of its differential to the (complex) universal enveloping algebra. Let $\chi_\lambda: \mathrm{USp}(2g)\to \mathbb{C}$ denote $V_\lambda$'s character. Then $\chi_\lambda$ is a Casimir eigenfunction, in the sense that there is a constant $c_\lambda\in \mathbb{C}$ such that $\mathfrak{C}\chi_\lambda = c_\lambda\chi_\lambda$.
\end{claim}
\begin{proof}
For $k\in \mathrm{USp}(2g)$, and $X\in \mathfrak{k}$, we have
\[
(X\chi_\lambda)(k) := \dfrac{d}{dt}\chi_\lambda(k\exp(tX))|_{t=0} = \dfrac{d}{dt}\mathrm{tr}(\rho_\lambda(k\exp(tX)))|_{t=0}.
\]
By multiplicativity of $\rho_\lambda$ and the linearity of the trace, this becomes
\[
(X\chi_\lambda)(k) = \mathrm{tr}(\rho_\lambda(k)\dfrac{d}{dt}\rho_{\lambda}(\exp(tX))|_{t=0}) = \mathrm{tr}(\rho_\lambda(k)d\rho_{\lambda}(X)).
\]
Then $(\mathfrak{C}\chi_\lambda)(k) = \mathrm{tr}(\rho_\lambda(k)d\rho_{\lambda}(\mathfrak{C}))$. Since, by Remark~\ref{rem:casimir}, $\mathfrak{C}\in \mathcal{Z}(\mathcal{U}_{\mathbb{C}})$ is in the center of $\mathcal{U}_\mathbb{C}$, $d\rho_{\lambda}(\mathfrak{C})$ is a $\mathrm{USp}(2g)$-equivariant endomorphism of $V_\lambda$, hence, by Schur's lemma, is given by a multiplication by a scalar, which we denote by $c_\lambda \in \mathbb{C}$.
\end{proof}

A careful calculation of the Casimir's action on a highest weight vector for $V_\lambda$ allows one to determine the scalar $c_\lambda$ precisely.
\begin{theorem}[{\cite[\S25.4, Lemma~25.6]{EtingofLieGroupsLieAlgebras}}]\label{thm:etingof}
The scalar $c_\lambda$ of Claim~\ref{claim:clambda} may be expressed explicitly as follows
\[
c_\lambda = \dfrac{1}{2}\sum_{j=1}^g \lambda_j^2 + \sum_{j=1}^g(g-j+1)\lambda_j.
\]
\end{theorem}

We are ready to state and prove our second-moment estimate.
\begin{proposition}[second-moment estimate]\label{prop:second_moment}
For all $n\ge 1$, one has
\[
\sum_{\lambda}m_n(\lambda)d_\lambda ||\lambda||^2 \ll nb_n.
\]
\end{proposition}

\begin{proof}
It follows from Theorem~\ref{thm:etingof} that
\[
(\mathfrak{C}\chi_n)(I) = \sum_\lambda m_n(\lambda)(\mathfrak{C}\chi_\lambda)(I) = \sum_\lambda m_n(\lambda)c_\lambda \chi_\lambda(I).
\]
Observing that $c_\lambda \gg ||\lambda||^2$ and recalling that $\chi_\lambda(I) = \dim V_\lambda = d_\lambda$, we obtain
\[
(\mathfrak{C}\chi_n)(I) \gg \sum_\lambda m_n(\lambda)d_\lambda||\lambda||^2.
\]
On the other hand, we may bound $(\mathfrak{C}\chi_n)(I)$ from above using Filip's formula (see Theorem~\ref{thm:filip_lie}), which writes $\chi_n$ explicitly via
\[
\chi_n(k) = \dfrac{1}{n}\sum_{d|n}\mu(n/d)P_d(\chi_V(k^{n/d})).
\]
Our goal is therefore to understand $\mathfrak{C}P_d(\chi_V(k^{n/d}))$. Let $X\in \mathfrak{k}$. As a first step we compute $D_XP_d(\chi_V(k^{n/d}))$. By definition,
\[
D_XP_d(\chi_V(k^{n/d})) = \dfrac{d}{dt}P_d(\chi_V((k\exp(tX))^{n/d}))|_{t=0}.
\]
Write $h: \mathrm{USp}(2g)\to \mathrm{USp}(2g)$, where $h(k) = k^{n/d}$, so that $P_d(\chi_V(k^{n/d})) = (P_d\circ \chi_V\circ h)(k)$. By the chain rule,
\[
D_X(P_d\circ\chi_V\circ h)(k) = (dP_d)_{\chi_V(h(k))}\circ (d\chi_V)_{h(k)}\circ dh_k\circ \gamma_X'(0),
\]
where $\gamma_X(t) = k\exp(tX)$. Since $dh_k\circ \gamma_X'(0) = (D_Xh)(k)$, we may compute explicitly,
\[
(D_Xh)(k) = \dfrac{d}{dt}\left(\underbrace{k\exp(tX)\cdot \dots\cdot k\exp(tX)}_{n/d\text{ times}}\right)|_{t=0} = kXk^{n/d-1} + \dots + k^{n/d}X = \sum_{j=0}^{n/d-1}k^{n/d}\mathrm{Ad}(k^{-j})(X).
\]
Plugging back, we obtain
\[
D_X(P_d\circ\chi_V\circ h)(k) = (dP_d)_{\chi_V(h(k))}\circ (d\chi_V)_{h(k)}\circ\left(\sum_{j=0}^{n/d-1}k^{n/d}\mathrm{Ad}(k^{-j})(X)\right).
\]
Since $P_d$ is holomorphic at a neighborhood of $\chi_V(I^{n/d}) = \chi_V(I) = 2g$, at a sufficiently small neighborhood, $(dP_d)_{\chi_V(h(k))}$ is simply multiplication by $P_d'(\chi_V(k^d))$, hence we may simplify the above into
\[
D_X(P_d\circ\chi_V\circ h)(k) = P_d'(\chi_V(h(k)))\cdot (d\chi_V)_{h(k)}\circ\left(\sum_{j=0}^{n/d-1}h(k)\mathrm{Ad}(k^{-j})(X)\right).
\]
Furthermore, by definition, for each $0\le j\le n/d-1$ we have 
\[
(d\chi_V)_{h(k)}(h(k)\mathrm{Ad}(k^{-j})(X)) = (D_{\mathrm{Ad}(k^{-j})(X)}\chi_V)(h(k)) = (D_{\mathrm{Ad}(k^{-j})(X)}\chi_V)(k^{n/d}).
\]
Plugging back, we get
\[
D_X(P_d\circ\chi_V\circ h)(k) = P_d'(\chi_V(h(k)))\cdot\left(\sum_{j=0}^{n/d-1}(D_{\mathrm{Ad}(k^{-j})(X)}\chi_V)(k^{n/d})\right).
\]
Applying $D_X$ once more, following the same process, we obtain
\[
D_X^2(P_d\circ\chi_V\circ h)(k) = P_d''(\chi_V(h(k)))\cdot\left(\sum_{j=0}^{n/d-1}(D_{\mathrm{Ad}(k^{-j})(X)}\chi_V)(k^{n/d})\right)^2
\]
\[
+ P_d'(\chi_V(h(k)))\cdot \sum_{i,j=0}^{n/d-1}(D_{\mathrm{Ad}(k^{-i})(X)}D_{\mathrm{Ad}(k^{-j})(X)}\chi_V)(k^{n/d}).
\]
Specializing at $k=I$, the term $P_d''(\chi_V(h(k)))\cdot\left(\sum_{j=0}^{n/d-1}(D_{\mathrm{Ad}(k^{-j})(X)}\chi_V)(k^{n/d})\right)^2$ vanishes. Indeed, observing that for any $Y\in \mathfrak{k}$, and any $k$,
\[
(D_Y\chi_V)(k^{n/d}) = \mathrm{tr}(k^{n/d}Y),
\]
$(D_Y\chi_V)(I) = 0$, since $Y$ is traceless. In particular, $D_X^2(P_d\circ\chi_V\circ h)(I)$ specializes into
\[
D_X^2(P_d\circ\chi_V\circ h)(I) = \dfrac{n^2}{d^2}P_d'(2g)\cdot (D_X^2\chi_V)(I).
\]
It follows that
\[
\mathfrak{C}(P_d\circ\chi_V\circ h)(I) = \dfrac{n^2}{d^2}P_d'(2g)\cdot (\mathfrak{C}\chi_V)(I).
\]
By Claim~\ref{claim:clambda}, $\mathfrak{C}\chi_V = c_V \chi_V$, while, by Theorem~\ref{thm:etingof}, $c_V = g + 1/2$. Therefore,
\[
|\mathfrak{C}(P_d\circ\chi_V\circ h)(I)| \ll_g \dfrac{n^2}{d^2}|P_d'(2g)|.
\]
Next up, recall that
\[
P_d(x) = \left(\dfrac{x + \sqrt{x^2 - 4}}{2}\right)^d + \left(\dfrac{x - \sqrt{x^2 - 4}}{2}\right)^d.
\]
Differentiating, we get
\[
P_d'(x) = d\left[\left(\dfrac{x + \sqrt{x^2 - 4}}{2}\right)^{d-1}\left(\dfrac{1}{2} + \dfrac{x}{2\sqrt{x^2 - 4}}\right) + \left(\dfrac{x - \sqrt{x^2 - 4}}{2}\right)^{d-1}\left(\dfrac{1}{2} - \dfrac{x}{2\sqrt{x^2 - 4}}\right)\right].
\]
Substituting $x = 2g$, we obtain the bound:
\[
P_d'(2g) \ll d\alpha^{d-1} \ll d\alpha^d.
\]
Plugging back, our bound for $|\mathfrak{C}(P_d\circ\chi_V\circ h)(I)|$ becomes
\[
|\mathfrak{C}(P_d\circ\chi_V\circ h)(I)|\ll \dfrac{n^2\alpha^d}{d}.
\]
Substituting into Filip's formula, we obtain
\[
|\mathfrak{C}\chi_n(I)| \ll \dfrac{1}{n}\sum_{d|n}|\mathfrak{C}P_d(\chi_V(h(k)))| \ll n\sum_{d|n}\dfrac{\alpha^d}{d} \ll nb_n + n\sum_{d\le n/2}\dfrac{\alpha^d}{d} \ll nb_n,
\]
which completes the proof.
\end{proof}

We are ready to complete the proof of Theorem~\ref{thm:power_saving}, which states that 
\[
H^2(G_{\mathbb{Q},\Sigma}, L_n) \ll b_n\cdot n^{-\frac{1}{2(D+1)}}.
\]

\begin{proof}[Proof of Theorem~\ref{thm:power_saving}]
Since the bound in Theorem~\ref{thm:power_saving} is asymptotic, we may assume $n$ is large enough so that $n\not\in E\subsetneq\Z$, with $E$ the finite set of bad weights of Theorem~\ref{thm:h2-bound}. Decomposing $L_n$ as in~\eqref{eq:ln-decomp}, Theorem~\ref{thm:h2-bound} gives
\[
\dim_\Qp H^2(G_{\mathbb{Q},\Sigma}, L_n) = \sum_\lambda m_n(\lambda)\cdot\dim_\Qp H^2(G_{\mathbb{Q},\Sigma}, V_{\lambda,n}) \ll 
\sum_\lambda m_n(\lambda)\cdot\min\left\{d_\lambda, \dfrac{(1 + ||\lambda||)^{D-g-1}}{d_\lambda}\right\}.
\]
Let $R > 0$ be an optimization parameter, which may depend on our fixed $n$, to be decided upon in hindsight. We split the above sum with respect to how large $||\lambda||$ is. Writing
\[
\sum_\lambda m_n(\lambda)\cdot\min\left\{d_\lambda, \dfrac{(1 + ||\lambda||)^{D-g-1}}{d_\lambda}\right\} = 
\sum_{||\lambda|| < R\sqrt{n}} + \sum_{||\lambda|| \ge R\sqrt{n}} = (I) + (II).
\]
As for $(I)$, we shall bound each $\dim_\Qp H^2(G_{\mathbb{Q},\Sigma}, V_{\lambda,n})$ using the nontrivial bound, $\dfrac{(1 + ||\lambda||)^{D-g-1}}{d_\lambda}$. We obtain
\[
(I) \ll \sum_{||\lambda|| < R\sqrt{n}}m_n(\lambda)\dfrac{(1 + ||\lambda||)^{D-g-1}}{d_\lambda}.
\]
By Proposition~\ref{prop:mn-bound}, $m_n(\lambda) \ll d_\lambda b_n n^{-D/2}$, and plugging into the above bound on $(I)$, we obtain
\[
(I) \ll b_n n^{-D/2}\sum_{||\lambda|| < R\sqrt{n}}(1 + ||\lambda||)^{D-g-1}.
\]
Bounding $\sum_{||\lambda|| < R\sqrt{n}}(1 + ||\lambda||)^{D-g-1}$ amounts to a basic lattice point count problem. As there are at most $O(L^g)$ lattice points in $\mathbb{Z}^g$ of norm at most $L$, there are at most $O((R\sqrt{n})^g)$ summands in $\sum_{||\lambda|| < R\sqrt{n}}(1 + ||\lambda||)^{D-g-1}$. As each summand is bounded above by a constant multiple of $(R\sqrt{n})^{D-g-1}$, our bound on $(I)$ becomes
\[
(I) \ll b_n n^{-D/2}(R\sqrt{n})^{D-1} = b_n R^{D-1}/\sqrt{n}.
\]
We proceed to analyze $(II)$. Here, we apply the trivial bound $\dim_\Qp H^2(G_{\mathbb{Q},\Sigma}, V_{\lambda,n}) \ll d_\lambda$. $(II)$ becomes,
\[
(II) \ll \sum_{||\lambda|| \ge R\sqrt{n}}m_n(\lambda)d_\lambda \le \dfrac{1}{R^2n}\sum_{||\lambda|| \ge R\sqrt{n}}m_n(\lambda)d_\lambda||\lambda||^2,
\]
where the final inequality above is sometimes referred to as Rankin's trick. Extending the sum throughout all $\lambda$'s and applying the second-moment estimate of Proposition~\ref{prop:second_moment}, $(II)$ is bounded by
\[
(II) \ll \dfrac{1}{R^2n}\sum_{\lambda}m_n(\lambda)d_\lambda||\lambda||^2 \ll b_nR^{-2}.
\]
Fixing $R$ in hindsight to balance between our bounds on $(I)$ and $(II)$, i.e.~fixing $R$ to solve $R^{-2} = R^{D-1}/\sqrt{n}$, we have $R := n^{\frac{1}{2(D+1)}}$, and so
\[
(I) + (II) \ll b_n\cdot n^{-\frac{1}{D+1}},
\]
which completes the proof of our power saving bound.
\end{proof}

\subsection{Proof of Kim's conjecture}\label{ssec:kim-conj-proof}
We now complete our proof of Kim's conjecture, i.e.~the proof of Corollary~\ref{cor:kim-conj}, which states that the Chabauty--Kim condition,
\[
\dim H^1_f(G_{\mathbb{Q},\Sigma}, U_n^\et) < \dim U_n^\dR/F^0,
\]
holds for any smooth projective geometrically connected curve of genus at least $2$ over $\mathbb{Q}$, for $n$ large enough (depending on the curve). Our notations are as in the introduction, that is, here $U^\et_n$ and $U^\dR_n$ denote the $n$'th descending central series truncations of the pro-unipotent \'etale and de Rham fundamental groups of $C$ based at an auxiliary basepoint $b\in \mathcal{C}(\Z[1/S])$. We shall complete the proof of Corollary~\ref{cor:kim-conj} by analyzing the two sides of the desired dimension inequality. Most of the reduction steps we carry out here, reducing the Chabauty--Kim condition to an $H^2$-bound are due to Kim, and already appear in~\cite[\S3, proof of Proposition~2, pp.~122--124]{KIM}.

As for the de Rham side, we have:
\begin{proposition}\label{prop:de_rham_bound}
\[
\dim U_n^\dR/F^0 \gg \sum_{k\le n}b_k - g^n.
\]
\end{proposition}
\begin{proof}
Since $U_n^\dR$ is an algebraic group and $F^0 = F^0U_n^\dR$ is its subgroup, we have
\[
\dim U_n^\dR/F^0 = \dim U_n^\dR - \dim F^0U_n^\dR.
\]
Since this proof only deals with the de Rham side, we shall typically omit the superscript from $U^\dR$ and write $U$ within this proof, where there is no risk of confusion.

For each $k\ge 1$, write $\mathrm{gr}_kU$ for the $k$'th graded piece of the descending central series filtration of $U$. For each such $k$, we have a short exact sequence
\begin{equation}\label{eq:dcs_ses}
0\to \gr_{k+1} U\to U_{k+1} \to U_k \to 0,
\end{equation}
giving rise, by induction, to the dimension equality, $\dim U_n = \sum_{k\le n}\dim \mathrm{gr}_k U$. Olsson's \'etale-de Rham comparison isomorphism (see~\cite[Theorems~2.20-2.21]{CORWIN_ZEHAVI}) identifies the associated graded Lie algebras of the de Rham and \'etale unipotent fundamental groups after extension of scalars to $B_{\mathrm{dR}}$. Consequently, their graded pieces have the same dimensions, and in particular, $\dim \mathrm{gr}_k U = b_k$ for every $k$. Therefore,
\[
\dim U_n = \sum_{k\le n}b_k.
\]
As for $\dim F^0U_n$, since the descending central series is strict for the Hodge filtration (cf.~\cite[Corollary B.16]{CORWIN_ZEHAVI}), we similarly have
\[
\dim F^0U_n = \sum_{k\le n}\dim F^0\mathrm{gr}_kU.
\]
In~\cite[\S3, proof of Proposition~2, p.~123]{KIM}, Kim shows that for each $k$,
\[
\dim F^0\mathrm{gr}_kU \le g^k.
\]

Plugging back, our bound on $\dim F^0U_n$ becomes
\[
\dim F^0U_n \le \sum_{k\le n}g^k \ll g^n,
\]
which completes the proof.
\end{proof}

Next, we proceed to study the \'etale side. Since $H^1_f(G_{\mathbb{Q},\Sigma}, U_n^\et)\hookrightarrow H^1(G_{\mathbb{Q},\Sigma}, U_n^\et)$ is a closed embedding, we may trivially bound,
\begin{equation}\label{eq:mono_etale}
\dim H^1_f(G_{\mathbb{Q},\Sigma}, U_n^\et)\le \dim H^1(G_{\mathbb{Q},\Sigma}, U_n^\et).  
\end{equation}
The descending central series filtration gives a short exact sequence as in~\eqref{eq:dcs_ses}
\[
0\to \gr_{k+1} U^\et\to U_{k+1}^\et \to U_k^\et \to 0,
\]
which, in turn, gives rise to a long exact sequence in Galois cohomology,
\[
\dots\to H^1(G_{\mathbb{Q},\Sigma}, \gr_{k+1} U^\et)\to H^1(G_{\mathbb{Q},\Sigma}, U^\et_{k+1})\to H^1(G_{\mathbb{Q},\Sigma}, U^\et_k)\to \dots,
\]
which gives the upper bound
\[
\dim H^1(G_{\mathbb{Q},\Sigma}, U^\et_{k+1}) \le \dim H^1(G_{\mathbb{Q},\Sigma}, U^\et_{k}) + \dim H^1(G_{\mathbb{Q},\Sigma}, \gr_{k+1} U^\et).
\]
Here $\gr_{k+1} U^\et$ denotes a vector group, and $H^1(G_{\mathbb{Q},\Sigma}, \gr_{k+1} U^\et)$ is an affine space, whose dimension is exactly 
\[
\dim_\Qp H^1(G_{\mathbb{Q},\Sigma}, \gr_{k+1} U^\et(\Qp)) = \dim_\Qp H^1(G_{\mathbb{Q},\Sigma}, L_{k+1}),
\]
where $L_{k+1}$ is the $k+1$'th graded piece of the associated graded Lie algebra, and the equality of dimensions stems from the fact that the logarithm map identifies $\gr_{k+1} U^\et(\Qp)$ and $L_{k+1}$ as $p$-adic Galois representations. Applying this bound recursively, we obtain,
\[
\dim H^1(G_{\mathbb{Q},\Sigma}, U^\et_n) \le \sum_{k\le n}\dim H^1(G_{\mathbb{Q},\Sigma}, L_k).  
\]
By the Euler characteristic formula~\cite[Chapter~I, Theorem~5.1]{MilneADT}, one obtains for any finite dimensional $p$-adic Galois representation $W$,
\[
\dim H^1(G_{\mathbb{Q},\Sigma}, W) = \dim H^0(G_{\mathbb{Q},\Sigma}, W) + \dim H^2(G_{\mathbb{Q},\Sigma}, W) + \dim W^-,
\]
where $W^-$ denotes the $-1$-eigenspace of complex conjugation. We wish to apply this formula for $W = L_k$. We first argue that $L_k$ has no Galois invariants. Indeed, $L_k$ may be identified as a quotient of $L_1^{\otimes k}$, where $L_1 = H_1^\et(\overline{C}; \Qp)$. Since $C$ is smooth projective over $\mathbb{Q}$, its geometric \'etale homology is pure of weight $1$, hence $L_1^{\otimes k}$ is pure of weight $k$, and so is $L_k$. It follows that $\dim H^0(G_{\mathbb{Q},\Sigma}, L_k) = 0$.

Plugging back, we obtain
\begin{equation}\label{eq:h1-to-h2}
\dim H^1(G_{\mathbb{Q},\Sigma}, U^\et_n) \le \sum_{k\le n}\dim_\Qp H^2(G_{\mathbb{Q},\Sigma}, L_k) + \sum_{k\le n}\dim_\Qp L_k^-.  
\end{equation}

Therefore, our next goal is to bound $\dim L_k^-$.
\begin{proposition}\label{prop:lk-dim}
One has
\begin{equation}\label{eq:lk-_dim}
\dim_\Qp L_k^- = b_k/2 + O\left(\dfrac{\alpha^{n/2}}{n}\right).
\end{equation}
\end{proposition}
\begin{remark}
We note that in~\cite[\S3, proof of Proposition~2, p.~123]{KIM}, Kim shows that for all $k\ge 1$ odd, one has
\[
\dim_\Qp L_k^- = \dfrac{1}{2}\dim_\Qp L_k = b_k/2.
\]
The idea is that $L_k$ has a pure Hodge structure of weight $-k$, and therefore its complexification decomposes as
\[
L_k\otimes \mathbb{C} = \bigoplus_{a+b=-k}L_k^{a,b},
\]
with the $L_k^{a,b}$ the component of Hodge type $(a,b)$. Complex conjugation interchanges $L_k^{a,b}$ with $L_k^{b,a}$, and when $k$ is odd, these subspaces are necessarily different. In particular, denoting complex conjugation by $\sigma$, the Hodge pieces $L_k^{a,b}$ and $L_k^{b,a}$ contribute equally to the positive and negative $\sigma$-eigenspaces of $L_k$, as can be seen by associating an arbitrary basis $v_1,\dots,v_r$ of $L_k^{a,b}$, with the $\sigma$-eigenbasis
\[
(v_i + \sigma(v_i), v_i - \sigma(v_i))_{i=1}^r
\]
of $L_k^{a,b}\oplus L_k^{b,a}$. From the point of view of Kim's conjecture, the following argument, which recovers a uniform bound of $\dim_\Qp L_k^- = (1/2 + o(1))b_k$, valid for all $k$, even and odd, is redundant. We choose to include it here nevertheless, as we find it elegant, and believe it may be of independent interest.
\end{remark}

\begin{proof}[Proof of Proposition~\ref{prop:lk-dim}]
Write $\langle \sigma\rangle$ for the order $2$ group generated by complex conjugation, denoted $\sigma \in G_{\mathbb{Q},\Sigma}$. Then~\eqref{eq:lk-_dim} is equivalent to the statement that 
\[
|\mathrm{tr}(\sigma|L_k)| \ll \dfrac{\alpha^{n/2}}{n},
\]
which is what we aim to prove. With $L = \bigoplus_{k\ge 1}L_k$ the associated graded Lie algebra, write $\mathcal{U}(L)$ for the universal enveloping algebra of $L$. We start out by deriving a presentation of $\mathcal{U}(L)$ in terms of $T(V)$, where $T(V)$ is the tensor algebra generated by $V$. 

As before, we write $V := L_1$, and let $\mathbb{L}(V)$ denote the free Lie algebra on $V$. The Lie algebra $L$ admits a presentation as $\mathbb{L}(V)/(\rho)$, with $(\rho)$ the Lie ideal generated by a single element $\rho \in \mathrm{gr}_2\mathbb{L}(V)$; see~\S\ref{ssec:character-estimates}. Let $r\in T(V)$ denote the corresponding element of $T(V)$, obtained by replacing every Lie bracket $[u,v]$ appearing in $\rho$, with $u,v\in \mathrm{gr}_1\mathbb{L}(V) = V$, by the commutator
\[
X_u\otimes X_v - X_v\otimes X_u\in T(V).
\]
Then $\mathcal{U}(L)$ admits a presentation as $T(V)/(r)$, where $(r)$ denotes the two-sided ideal in $T(V)$ generated by $r$. The tensor algebra is naturally graded by tensor length, i.e.
\[
T(V) = \bigoplus_{n\ge 0}V^{\otimes n},
\]
and since $(r)$ is a graded ideal, the quotient, $\mathcal{U}(L) \cong T(V)/(r)$, is naturally graded as well. Write $A_n$ for its $n$'th  graded piece. For each $n\ge 2$, we have a short exact sequence of $\langle \sigma\rangle$-representations
\[
0\to A_{n-2}\otimes r\to A_{n-1}\otimes V\to A_n\to 0.
\]
It follows that
\[
\mathrm{tr}(\sigma|A_n) = \mathrm{tr}(\sigma|A_{n-1}\otimes V) - \mathrm{tr}(\sigma|A_{n-2}\otimes r).
\]
Now, recall that $\mathrm{tr}(\sigma|V) = 0$. Indeed, this is true since the Weil pairing $V^{\vee}\otimes V^\vee\to \Qp(-1)$ is perfect, and pairs $\sigma$-eigenvectors of the same eigenvalue to $0$, hence the positive and negative eigenspaces are equidimensional. Also, since $\Qp\cdot r \cong \Qp(1)$ as a Galois representation (see~\S\ref{ssec:character-estimates}), $\mathrm{tr}(\sigma|\Qp\cdot r) = -1$. It follows that
\[
\mathrm{tr}(\sigma|A_{n-1}\otimes V) = \mathrm{tr}(\sigma|A_{n-1})\cdot \mathrm{tr}(\sigma|V) = 0,
\]
and that
\[
\mathrm{tr}(\sigma|A_{n-2}\otimes r) = \mathrm{tr}(\sigma|A_{n-2})\cdot \mathrm{tr}(\Qp\cdot r) = -\mathrm{tr}(\sigma|A_{n-2}).
\]
Hence, for all $n\ge 2$,
\[
\mathrm{tr}(\sigma|A_{n}) = \mathrm{tr}(\sigma|A_{n-2}),
\]
and for every $n\ge 0$, one has
\[
\mathrm{tr}(\sigma|A_n) =
\begin{cases}
1& 2|n\\
0& 2\nmid n.
\end{cases}
\]
Let us denote by $H(t;\mathcal{U}(L))$ the generating series,
\[
H(t;\mathcal{U}(L)) := \sum_{n\ge 0}t^{2n} = \dfrac{1}{1-t^2}.
\]
We may express $H(t;\mathcal{U}(L))$ differently as follows. Since $L = \bigoplus_{k\ge 1}L_k$, we have
\[
\mathcal{U}(L) = \bigotimes_{k\ge 1}\mathcal{U}(L_k).
\]
By the equivariant PBW theorem (\cite[\S1, equation~(1.2), and \S4 p.~8]{FilipSurfaceLie}), 
\[
\mathrm{gr}\,\mathcal{U}(L_k) = S(L_k),
\]
that is, the associated graded algebra of $\mathcal{U}(L_k)$ is isomorphic to the symmetric algebra on $L_k$. The symmetric algebra on $L_k$ is naturally graded, and its $n$'th graded piece is isomorphic to $\mathrm{Sym}^nL_k$ not only as a $\Qp$-vector space, but also as a $\langle \sigma\rangle$-representation. $\mathcal{U}(L_k)$ embeds into $\mathcal{U}(L)$ as a $\langle\sigma\rangle$-subrepresentation, and its generators, which are identified with $L_k$, have degree $k$. For each $k\ge 1$, consider the generating series,
\[
H(t;\mathcal{U}(L_k)) := \sum_{n\ge 0}\mathrm{tr}(\sigma|\mathrm{Sym}^nL_k)t^{nk},
\]
then we have the identity
\begin{equation}\label{eq:hilbert-product}
H(t;\mathcal{U}(L)) = \prod_{k\ge 1}H(t;\mathcal{U}(L_k)).  
\end{equation}
In order to isolate $\mathrm{tr}(\sigma|L_k)$, we shall apply the operator $t\cdot \dfrac{d}{dt}\log$ to both sides. As for the left hand side, we get
\[
(t\dfrac{d}{dt})\log \dfrac{1}{1-t^2} = \dfrac{2t^2}{1 - t^2} = \sum_{n\ge 1}2t^{2n}.
\]

As for the right-hand side, fix $k\ge 1$ and denote the dimensions of the positive and negative eigenspaces of $L_k$ via:
\[
p_k := \dim L_k^+,\qquad n_k := \dim L_k^-.
\]
Then clearly,
\[
p_k+n_k = b_k,\qquad p_k-n_k = \mathrm{tr}(\sigma|L_k).
\]
Each $\sigma$-eigenvector $v\in L_k$ of eigenvalue $\epsilon\in\{\pm 1\}$ contributes a factor of $\sum_{n\ge 0}(\epsilon t^k)^n$ to $H(t;\mathcal{U}(L_k))$ and therefore,
\[
H(t;\mathcal{U}(L_k)) = \dfrac{1}{(1-t^k)^{p_k}}\cdot \dfrac{1}{(1+t^k)^{n_k}}.
\]
Taking $(t\dfrac{d}{dt})\log$, we get
\[
(t\dfrac{d}{dt})\log H(t;\mathcal{U}(L_k)) = -p_kt\dfrac{d}{dt}\log(1-t^k) - n_kt\dfrac{d}{dt}\log(1+t^k) = kt^k\left[\dfrac{p_k}{1-t^k} - \dfrac{n_k}{1+t^k}\right]
\]
\[
= kt^k\sum_{m\ge 0}(p_k - (-1)^m n_k)t^{mk} = \sum_{m\ge 1}k\cdot \mathrm{tr}(\sigma^m|L_k)t^{mk}.
\]
Summing over $k\ge 1$, the right-hand side of~\eqref{eq:hilbert-product} becomes
\[
\sum_{k\ge 1}\sum_{m\ge 1}k\cdot\mathrm{tr}(\sigma^m|L_k)t^{mk} = \sum_{n\ge 1}\left[\sum_{d|n}d\cdot \mathrm{tr}(\sigma^{n/d}|L_d)\right]t^n.
\]
Comparing coefficients of $t^n$, we get, for every $n\ge 1$,
\[
\sum_{d|n}d\cdot \mathrm{tr}(\sigma^{n/d}|L_d) = 2\delta_{2\nmid n} :=
\begin{cases}
2& 2\nmid n\\
0& 2\mid n.
\end{cases}
\]
For each $d|n$, since the eigenvalues of $\sigma^{n/d}$ all have absolute value $1$, we have the naive bound, $\mathrm{tr}(\sigma^{n/d}|L_d) \le \dim L_d = b_d$. Then, it follows that
\[
n\cdot|\mathrm{tr}(\sigma|L_n)| \le 2 + \sum_{d|n,d\neq n}|d\cdot \mathrm{tr}(\sigma^{n/d}|L_d)| \ll 1 + \sum_{d\le n/2}d\cdot b_d
\]
\[
\ll 1 + \sum_{d\le n/2}\alpha^d \ll \alpha^{n/2},
\]
where we used Proposition~\ref{prop:bn-growth} for the bound $k\cdot b_k\ll \alpha^k$. It follows that
\[
|\mathrm{tr}(\sigma|L_n)| \ll \dfrac{\alpha^{n/2}}{n},
\]
which completes the proof.
\end{proof}

We are finally ready to complete the proof of Kim's conjecture.
\begin{proof}[Proof of Corollary~\ref{cor:kim-conj}]
Plugging in our lower bound on $\dim U^\dR_n/F^0$ of Proposition~\ref{prop:de_rham_bound}, and our upper bound on $\dim H^1_f(G_{\mathbb{Q},\Sigma}, U^\et_n)$ obtained by combining~\eqref{eq:mono_etale} and~\eqref{eq:h1-to-h2}, the Chabauty--Kim's condition for $C/\mathbb{Q}$ (see~\eqref{eq:kim_dim_ineq}) reduces to
\[
\sum_{k\le n}\dim H^2(G_{\mathbb{Q},\Sigma}, L_k) + \sum_{k\le n}\dim L_k^- < \sum_{k\le n}b_k - O(g^n).
\]
Plugging in our bound on $\dim L_k^-$ of Proposition~\ref{prop:lk-dim} and our power saving bound on $\dim H^2(G_{\mathbb{Q},\Sigma}, L_k)$ of Theorem~\ref{thm:h2-bound}, the above inequality reduces to
\[
A\sum_{k\le n}b_k\cdot k^{-\frac{1}{D+1}} + \dfrac{1}{2}\sum_{k\le n}b_k + \sum_{k\le n}O\left(\dfrac{\alpha^{k/2}}{k}\right) + O(g^n) < \sum_{k\le n}b_k,
\]
where $A = A(C, S, p) > 0$ is some absolute constant. This is equivalent to
\begin{equation}\label{eq:reduction_final}
A\sum_{k\le n}b_k\cdot k^{-\frac{1}{D+1}} + \sum_{k\le n}O\left(\dfrac{\alpha^{k/2}}{k}\right) + O(g^n) < \dfrac{1}{2}\sum_{k\le n}b_k.  
\end{equation}
\[
\]
Observe that $\sum_{k\le n}\dfrac{\alpha^{k/2}}{k} = O(\alpha^{n/2})$, and that both $g^n$ and $\alpha^{n/2}$ are $o(\sum_{k\le n}b_k\cdot k^{-\frac{1}{D+1}})$. Hence, there is a larger constant $B > A > 0$, such that for all $n\ge 1$,
\[
A\sum_{k\le n}b_k\cdot k^{-\frac{1}{D+1}} + \sum_{k\le n}O\left(\dfrac{\alpha^{k/2}}{k}\right) + O(g^n) < B\sum_{k\le n}b_k\cdot k^{-\frac{1}{D+1}}.
\]
Thus,~\eqref{eq:reduction_final} reduces to
\[
B\sum_{k\le n}b_k\cdot k^{-\frac{1}{D+1}} < \dfrac{1}{2}\sum_{k\le n}b_k,
\]
which holds for all $n$ large enough.
\end{proof}

\printbibliography{}
\end{document}